\documentclass[reqno]{amsart}
\usepackage{amsmath}
\usepackage{amscd}
\usepackage{graphics}
\usepackage{latexsym}
\theoremstyle{plain}
\newtheorem{theorem}{Theorem}[section]
\newtheorem{thm}[theorem]{Theorem}
\newtheorem{cor}[theorem]{Corollary}
\newtheorem{lem}[theorem]{Lemma}
\newtheorem{prop}[theorem]{Proposition}

\theoremstyle{definition}
\newtheorem{defn}[theorem]{Definition}
\newtheorem{notat}[theorem]{Notation}

\newtheorem{hyp}[theorem]{Hypothesis}
\newtheorem{conj}[theorem]{Conjecture}
\newtheorem{ques}[theorem]{Question}
\newtheorem{rmk}[theorem]{Remark}

\theoremstyle{remark}

\newcommand{\QQ}{\mathbb{Q}}

\newcommand{\ZZ}{\mathbb{Z}}
\newcommand{\CC}{\mathbb{C}}

\newcommand{\AAA}{\mathbb{A}}
\newcommand{\PP}{\mathbb{P}}

\newcommand{\SP}{\text{Spec }}

\newcommand{\Kbm}[2]{\overline{\mathcal{M}}_{#1}(#2)}
\newcommand{\kbm}[2]{\overline{M}_{#1}(#2)}
\newcommand{\kbmo}[2]{{M}_{#1}(#2)}
\newcommand{\Kbmo}[2]{{\mathcal{M}}_{#1}(#2)}

\newcommand{\lt}{\left}
\newcommand{\rt}{\right}

\newcommand{\mc}{\mathcal}

\newcommand{\OO}{\mathcal O}

\def\PP{{\mathbb P}}

\begin{document}

\title[Kodaira dimension]{The Kodaira dimension of spaces of rational
curves on low degree hypersurfaces}

\author[Starr]{Jason Michael Starr}
\address{Department of Mathematics \\
  Massachusetts Institute of Technology \\ Cambridge MA 02139}
\email{jstarr@math.mit.edu} 
\date{\today}

\begin{abstract}
For a hypersurface in complex projective space $X\subset \PP^n$, we
investigate the singularities and Kodaira dimension of the Kontsevich
moduli spaces $\Kbm{0,0}{X,e}$ parametrizing rational curves of degree
$e$ on $X$.  If $d+e \leq n$ and $X$ is a general hypersurface of
degree $d$, we prove that $\Kbm{0,0}{X,e}$ has only canonical
singularities and we conjecture the same is true for the coarse moduli
space $\kbm{0,0}{X,e}$.  We prove that this conjecture is implied by
the ``inversion of adjunction'' conjecture of Koll\'{a}r and Shokurov.
Also we compute the canonical divisor of $\Kbm{0,0}{X,e}$ and show
that for most pairs $(d,e)$ with $n \leq d^2 \leq n^2$, the canonical
divisor is a \emph{big} divisor.  When combined with the above
conjecture, this implies that in many cases $\kbm{0,0}{X,e}$ is a
variety of general type.  This investigation is motivated by the
question of which Fano hypersurfaces are unirational.
\end{abstract}

\maketitle

\section{Introduction}~\label{intro}

\medskip\noindent 
Let $X\subset \PP^n$ be a general hypersurface of degree $d$ in
complex projective space $\PP^n$.  The Kontsevich moduli space
$\Kbm{0,0}{X,e}$ is a proper, Deligne-Mumford stack which contains as
an open substack the scheme parametrizing smooth, rational curves of
degree $e$ on $X$ (c.f. ~\cite{FP}).  Except when $d=1$, $d=2$ or
$e=1$, very little is known about the singularities and Kodaira
dimension of $\Kbm{0,0}{X,e}$.  When $d=1$ and $d=2$, the spaces
$\Kbm{0,0}{X,e}$ are all smooth ~\cite{Kv} and rational ~\cite{KP}.
When $e=1$, the space $\Kbm{0,0}{X,e}$ is just the space of lines on
$X$, which is completely understood ~\cite[Thm. V.4.3]{K}.

\medskip\noindent 
What is known for $d \geq 3$ and $e \geq 2$?  If $d < \frac{n+1}{2}$,
it is proved in ~\cite{HRS2} that $\Kbm{0,0}{X,e}$ is integral of the
expected dimension with only local complete intersection
singularities.  If also $d^2 + d +1 < n$, it is proved in ~\cite{HS2}
that $\Kbm{0,0}{X,e}$ has negative Kodaira dimension.  In fact
$\Kbm{0,0}{X,e}$ is \emph{rationally connected}.  If $d \geq n-1$, the
open substack of $\Kbm{0,0}{X,e}$ parametrizing smooth rational curves
is not Zariski dense: the locus of \emph{multiple covers of lines}
yields an irreducible component of $\Kbm{0,0}{X,e}$ not contained in
the closure of this open set.  However all evidence suggests that for
$d\leq n-2$ and for $X\subset \PP^n$ a general hypersurface of degree
$d$, the stack $\Kbm{0,0}{X,e}$ is irreducible for all $e\geq 1$.

\medskip\noindent
\begin{ques} \label{ques-0} 
For $d \leq n-2$ and $X\subset \PP^n$ a general hypersurface of degree
$d$, what type of singularities does $\Kbm{0,0}{X,e}$ have?  For which
$(n,d,e)$ are the singularities terminal, resp. canonical, log
canonical?
\end{ques}

\medskip\noindent 
Fix $e > 1$ and $d\geq 1$.  Let $\PP^N$ denote the projective space
parametrizing degree $d$ hypersurfaces $X\subset \PP^n$.  Let $C_d
\subset \PP^N \times \Kbm{0,0}{\PP^n,e}$ denote the closed substack
parametrizing pairs $([X],[f:D \rightarrow X])$ with $X\subset \PP^n$
a hypersurface and $f:D \rightarrow X$ a stable map in
$\Kbm{0,0}{X,e}$.  The main theorem of this paper is the following:

\medskip\noindent
\setcounter{section}{7}
\setcounter{theorem}{4}
\begin{thm} 
If $e\geq 2$ and if $d+e \leq n$, then $C_d$ is an integral, normal,
local complete intersection stack of the expected dimension which has
canonical singularities.
\end{thm}

\medskip\noindent
\setcounter{section}{7}
\setcounter{theorem}{5}
\begin{cor} 
If $e \geq 2$ and if $d+e \leq n$, then for a general hypersurface
$X\subset \PP(V)$ of degree $d$, the Kontsevich moduli space
$\Kbm{0,0}{X,e}$ is an integral, normal, local complete intersection
stack of the expected dimension $(n+1-d)e+(n-3)$ with only canonical
singularities.
\end{cor}
\setcounter{section}{1}
\setcounter{theorem}{1}

\medskip\noindent 
It seems reasonable to expect that Theorem~\ref{thm-main} is sharp.
For instance the inequality $d+e \leq n$ is consistent with the
inequality $d \leq n-1$ necessary for $C_d$ to be irreducible.  On the
other hand, Corollary~\ref{cor-main} is certainly not sharp: it fails
to account for the cases $d=1$ and $d=2$ where $\Kbm{0,0}{X,e}$ is
smooth for every $e$.  For $d \geq 3$ and $e\geq 2$, the space
$\Kbm{0,0}{X,e}$ will always be singular, but hopefully it is not too
singular.

\medskip\noindent
\begin{ques} \label{ques-1} 
For which integers $d$ and $n$ is it true that for a very general
hypersurface $X\subset \PP^n$ of degree $d$, $X$ is unirational?
\end{ques}

\medskip\noindent 
This is the question which motivates this paper, although no new
answers are given here.  It is known that one must have $d \leq n$ or
else the Kodaira dimension of $X$ is nonnegative.  For $d=1,2,3$ and
$n \geq d$, it is known that a general hypersurface $X \subset \PP^n$
of degree $d$ is unirational.  For each integer $d$ there is an
integer $\phi(d)$ such that for $n \geq \phi(d)$ and $X\subset \PP^n$
a general hypersurface, then $X$ is unirational (~\cite{HMP},
~\cite[Chapter 23]{K+92}).  It is conjectured that if $d \leq n$ but
large compared to $n$ -- e.g. if $d=n$ for $n \geq 4$ -- then a
general hypersurface $X\subset \PP^n$ is \emph{not} unirational.  But
no example of such a hypersurface has been proved to be
non-unirational.

\subsection{Koll\'{a}r's approach}
\medskip\noindent 
The connection between Question~\ref{ques-1} and this paper comes from
a suggestion by Koll\'{a}r in ~\cite{KollarSimple}, that a necessary
condition for a variety $X$ of dimension at least $2$ to be
unirational is that for a general point $p\in X$, there exists a
rational surface $S\subset X$ containing $p$.

\medskip\noindent
\begin{defn} \label{defn-manyRS} 
An irreducible, projective variety $X$ over a field $k$ is \emph{swept
by rational surfaces} (resp. \emph{separably swept by rational
surfaces}) if there exists an irreducible variety $Z$ and a rational
transformation $F: Z \times \PP^2 \rightarrow X$ such that
\begin{enumerate}
\item 
the rational transformation $F$ is dominant (resp. dominant
and separable), and
\item 
the rational transformation $(\text{pr}_Z,F): Z \times \PP^2
\rightarrow Z \times X$ is generically finite to its image
(resp. generically finite and separable to its image).
\end{enumerate}
\end{defn}

\medskip\noindent
\begin{rmk} \label{rmk-manyRS} 
There are several obvious remarks.
\begin{enumerate}
\item 
This definition makes sense for any field $k$, without assuming that
$k$ is algebraically closed or of characteristic $0$ (although this is
the case of interest in the rest of the paper).
\item 
In the definition above, we may replace $Z$ by a separable, dominant
cover and the conditions will still hold.
\item 
Let $X$ be a variety which is separably swept by rational surfaces.
Let $S \subset Z \times X$ denote the image of $(\text{pr}_Z,F)$.  By
~\cite[Thm. III.2.4]{K}, the base change $S \otimes_{K(Z)}
\overline{K}(Z)$ is a rational surface.  By ~\cite{Coombes}, in fact
there is a separable dominant morphism $Z'\rightarrow Z$ and a
birational transformation over $Z'$, $G: Z' \times \PP^2 \rightarrow
Z' \times_Z S$.  Up to replacing $Z$ by $Z'$ and replacing $F$ by the
composite $Z' \times \PP^2 \xrightarrow{G} Z'\times_Z S
\xrightarrow{\text{pr}_S} S \xrightarrow{\text{pr}_X} X$, we may
assume that $(\text{pr}_Z,F)$ is actually birational to its image.
\item 
If $X$ is (separably) unirational and $\text{dim}(X) \geq 2$, then $X$
is (separably) swept by rational surfaces.
\item 
If $X$ is swept by rational surfaces (resp. separably swept by
rational surfaces), then $X$ is uniruled (resp. separably uniruled).
\item 
Let $\text{dim}(X) = n$.  The definition above is equivalent to the
stronger condition where we demand that $\text{dim}(Z) = n-2$.
Moreover we may demand that $Z$ is smooth over $k$.  In fact, by de
Jong's alterations of singularities, up to replacing $Z$ by a
generically \'etale cover, we may even demand that $Z$ is smooth and
projective.
\item 
The condition that $X$ is swept by rational surfaces (resp. separably
swept by rational surfaces) is a birational property which is
equivalent to the condition that there exists a finitely-generated
field extension $L/k$ of transcendence degree $n-2$ and a finite
super-extension (resp. finite separable super-extension) of $K(X)/k$
of the form $L(t_1,t_2)/K(X)$ (resp. such that the compositum $L*K(X)$
equals $L(t_1,t_2)$).
\item 
If $X$ is swept by rational surfaces (resp. separably swept by
rational surfaces) and $f:X \rightarrow X'$ is a generically finite,
dominant rational transformation (resp. generically \'etale, dominant
rational transformation), then also $X'$ is swept by rational surfaces
(resp. separably swept by rational surfaces).
\item 
Unlike the analogous situation of rationally connected varieties,
given a family of smooth, projective varieties in characteristic zero,
it is unclear whether the condition of being swept by rational
surfaces is a closed condition, or even an open condition, on fibers
of the family.
\end{enumerate}
\end{rmk}

\medskip\noindent 
It is technically more convenient to work with pencils of rational
curves than to work with rational surfaces, and, replacing $\PP^2$ by
the birational surface $\PP^1 \times \PP^1$, we can rephrase the
condition above in terms of pencils of rational curves.

\medskip\noindent
\begin{defn} \label{defn-sweeping} 
Let $X$ be a projective variety and $e\geq 1$ an integer.  An
integral, closed substack $Y \subset \Kbm{0,0}{X,e}$ is
\emph{sweeping} (resp. \emph{separably sweeping}) if
\begin{enumerate}
\item 
For a general geometric point of $Y$, the associated stable map
$f:C\rightarrow X$ has irreducible domain and is birational to its
image.
\item 
The restriction over $Y$ of the universal morphism, $f:Y
\times_{\Kbm{0,0}{X,e}} \mc{C} \rightarrow X$ is surjective
(resp. surjective and separable).
\end{enumerate}
\end{defn}

\medskip\noindent
\begin{rmk} \label{rmk-sweeping} 
There are several obvious remarks.
\begin{enumerate}
\item
This definition makes sense even if $k$ is not algebraically closed or
of characteristic $0$; although in case of positive characteristic one
should keep in mind that $\Kbm{0,0}{X,e}$ may only be an Artin
algebraic stack with finite diagonal (not a Deligne-Mumford stack).
\item 
If $Y \subset \Kbm{0,0}{X,e}$ is sweeping, then there is a Zariski
dense open substack which is a scheme.  In particular, it makes sense
to ask whether $Y$ is uniruled (resp. separably uniruled).
\item 
If $Y \subset \Kbm{0,0}{X,e}$ is sweeping (resp. separably sweeping),
then for the irreducible component $M \subset \Kbm{0,0}{X,e}$ which
contains $Y$, also $M$ is sweeping (resp. separably sweeping).
\item
If $X$ is uniruled (resp. separably uniruled), then there is an
integer $e$ and an irreducible component $M \subset \Kbm{0,0}{X,e}$
which is sweeping (resp. separably sweeping).
\end{enumerate}
\end{rmk}

\medskip\noindent  
\begin{lem} \label{lem-sweeping} 
Let $X$ be a projective variety over a field $k$ (not necessarily
algebraically closed or of characteristic $0$).  The variety $X$ is
swept by rational surfaces (resp. separably swept by rational
surfaces) iff there exists an integer $e$ and a substack $Y \subset
\Kbm{0,0}{X,e}$ which is sweeping (resp. separably sweeping) such that
$Y$ is uniruled (resp. separably uniruled).
\end{lem}

\medskip\noindent
\begin{proof}
Let $X$ be a variety of dimension $n$ which is swept by rational
surfaces (resp. separably swept by rational surfaces).  Then there
exists a smooth quasi-projective variety $Z$ of dimension $n-2$ and a
rational transformation $F:Z \times \PP^1 \times \PP^1 \rightarrow X$
which is dominant (resp. dominant and separable) such that
$(\text{pr}_Z,F): Z \times \PP^1 \times \PP^1 \rightarrow Z \times X$
is generically finite (resp. birational to its image).  The
indeterminacy locus of $F$, say $I \subset Z \times \PP^1 \times
\PP^1$, is a subvariety which has codimension at least $2$ at every
point.  The projection $\text{pr}_{12}(I) \subset Z \times \PP^1$ has
codimension at least $1$ at every point, i.e. it is contained in a
divisor $D \subset Z \times \PP^1$.  Denoting $U = Z\times \PP^1 - D$,
the rational transformation $F:U \times \PP^1 \rightarrow X$ is a
regular morphism.  Up to shrinking $U$ further, the morphism
$(\text{pr}_{12},F): U \times \PP^1 \rightarrow U \times X$ is a
finite morphism (resp. a finite, birational morphism).  Let $C'
\subset U \times X$ denote the image and let $G:C \rightarrow U \times
X$ denote the normalization of $C'$.  In case $X$ is separably swept
by rational surfaces, $G$ is the same as $(\text{pr}_{12},F)$.

\medskip\noindent 
Up to shrinking $U$ further, the morphism $\text{pr}_U:C \rightarrow
U$ is a smooth, proper morphism of relative dimension $1$.  Moreover
every geometric fiber is dominated by $\PP^1$, and therefore the
geometric fibers are connected curves which are isomorphic to $\PP^1$.
So $\text{pr}_U:C \rightarrow U$ is a family of genus $0$ curves and
$G:C \rightarrow U \times X$ is a family of stable maps of genus $0$
to $X$ of some degree $e$.  There is an induced $1$-morphism $\zeta:U
\rightarrow \Kbm{0,0}{X,e}$.  Define $Y \subset \Kbm{0,0}{X,e}$ to be
the closed image substack of $\zeta$.

\medskip\noindent 
The claim is that $Y$ is sweeping (resp. separably sweeping) and that
the coarse moduli space of $Y$ is uniruled (resp. separably uniruled).
There is an open substack of $Y$ over which the geometric fibers of
the universal curve are irreducible.  By construction, this open
substack contains the image of $\zeta$, and so it is Zariski dense in
$Y$.  Similarly, the condition that the stable maps are birational to
their images is an open condition on $Y$.  Since this condition holds
on the image of $\zeta$, it holds on an open dense subset of $Y$.
Therefore $Y$ satisfies Item ($1$) of Definition ~\ref{defn-sweeping}.
Now $F: U \times \PP^1 \rightarrow X$ is dominant Therefore also $G:C
\rightarrow X$ is dominant, and so $f:Y \times_{\Kbm{0,0}{X,e}} \mc{C}
\rightarrow X$ is dominant.  So $Y$ satisfies Item ($2$) of
Definition~\ref{defn-sweeping}, i.e. $Y$ is sweeping.  Moreover if $X$
is separably swept by rational surfaces, then $F$ is dominant,
generically finite and separable which implies that also $f:Y
\times_{\Kbm{0,0}{X,e}} \mc{C} \rightarrow X$ is dominant, generically
finite and separable (since a sub-extension of a separable field
extension is separable), so $Y$ is separably sweeping.

\medskip\noindent
Consider $\zeta:U \rightarrow Y$.  The claim is that $\zeta$ is
dominant and generically finite (resp. dominant, generically finite
and separable).  There is a factorization of $F$ of the form
\begin{equation}
U\times \PP^1 \rightarrow U \times_Y
(Y \times_{\Kbm{0,0}{X,e}} \mc{C}) = C \rightarrow Y
\times_{\Kbm{0,0}{X,e}} \mc{C} \xrightarrow{f} X.
\end{equation}
Since $F$ is dominant and generically finite (resp. dominant,
generically finite and separable), each of these factors of $F$ is
dominant and generically finite (resp. dominant, generically finite
and separable).  In particular $U \times_Y (Y \times_{\Kbm{0,0}{X,e}}
\mc{C}) \rightarrow Y \times_{\Kbm{0,0}{X,e}} \mc{C}$ is dominant and
generically finite (resp. and separable).  But this is just the
base-change of $\zeta:U \rightarrow Y$ by the smooth surjective
morphism $\text{pr}_Y: Y \times_{\Kbm{0,0}{X,e}} \mc{C} \rightarrow
Y$.  Therefore $\zeta:U \rightarrow Y$ is dominant and generically
finite (resp. dominant, generically finite and separable).  Since $U$
is an open subset of $Z \times \PP^1$, we conclude that $Y$ is
uniruled (resp. separably uniruled).  This proves the forward
direction of the lemma.

\medskip\noindent
Conversely, suppose that $Y \subset \Kbm{0,0}{X,e}$ is sweeping
(resp. separably sweeping) and $Y$ is uniruled (resp. separably
uniruled).  Then we can find a rational transformation $\zeta:Z \times
\PP^1 \rightarrow Y$ which is dominant and generically finite
(resp. dominant, generically finite and separable).

\medskip\noindent
Let $\pi:C \rightarrow Z \times \PP^1$ be a projective completion of
the pullback of the universal family of curves over $\Kbm{0,0}{X,e}$,
and let $H:C \rightarrow X$ be the pullback of the universal stable
map.  Blowing up $C$ if necessary, we may suppose that $C$ is normal
and that $H$ is a regular morphism.  Moreover, over a dense open
subset of $Z\times \PP^1$, the projection $\pi$ is a smooth, proper
morphism whose geometric fibers are connected curves of genus $0$.  So
the geometric generic fiber of $\text{pr}_Z\circ \pi:C \rightarrow Z$
is a conic bundle over $\PP^1$.  By Tsen's theorem
(c.f. ~\cite[Cor. IV.6.6.2]{K}), the base-change $C \otimes_{K(Z)}
\overline{K}(Z)$ is a rational surface.  By ~\cite{Coombes}, in fact
there is a dominant, generically finite and separable morphism
$Z'\rightarrow Z$ and a birational transformation over $Z'$, $G:Z'
\times \PP^2 \rightarrow Z' \times_Z C$.  Observe that the composition
$Z' \times \PP^1 \rightarrow Z'\times \PP^1 \rightarrow Y$ is still
dominant and generically finite (resp. dominant, generically finite
and separable).  Therefore we may replace $Z$ by $Z'$ so that $C$ is
birational to $Z \times \PP^2$ over $Z$.

\medskip\noindent
Denote by $F$ the rational transformation $H\circ G: Z \times \PP^2
\rightarrow X$.  There is a factorization of $F$ of the form
\begin{equation}
Z \times \PP^2 \xrightarrow{G} C \rightarrow Y \times_{\Kbm{0,0}{X,e}}
\mc{C} \xrightarrow{f} X.
\end{equation}
By assumption, each of these factors is dominant (resp. dominant and
separable).  Therefore $F$ is dominant (resp. dominant and separable).
Consider $(\text{pr}_Z \circ \pi,F)\circ G:Z \times \PP^2 \rightarrow
Z\times X$.  By the hypotheses on $Y$, this morphism is generically
finite (resp. generically finite and separable to its image).
Therefore $X$ is swept by rational surfaces (resp. separably swept by
rational surfaces).  This finishes the proof of the lemma.
\end{proof}

\medskip\noindent
Because of the lemma, it is natural to try to understand the sweeping
substacks of $\Kbm{0,0}{X,e}$, and in particular to try to understand
the Kodaira dimension of these substacks (recall that a standard
conjecture from the minimal model program predicts that an algebraic
variety is uniruled iff its Kodaira dimension is negative).  If
$X\subset \PP^n$ is a hypersurface of degree $d \leq n$, then for $e
>> 0$, an irreducible component of $\Kbm{0,0}{X,e}$ will itself be
sweeping.  So the first step is to determine the Kodaira dimensions of
the irreducible components of $\Kbm{0,0}{X,e}$ which are sweeping.

\medskip\noindent
As mentioned above, for $X$ a general hypersurface of degree $d <
\frac{n+1}{2}$ the stacks $\Kbm{0,0}{X,e}$ are irreducible and
reduced, and the same is conjectured whenever $d < n-1$.  Under the
hypothesis that each $\Kbm{0,0}{X,e}$ is an integral, normal stack of
the expected dimension, one can compute the canonical divisor of
$\Kbm{0,0}{X,e}$, i.e. one can compute the \emph{expected canonical
divisor}.  This is carried out in Section~\ref{sec-canKBM}, and is a
straightforward extension of the derivation in ~\cite{Pand97}.  The
exciting observation is that when $d$ satisfies the inequalities
$d\leq n-4$ and $d^2 \geq n+2$, then for almost all values of $e$, the
canonical divisor is \emph{big}.  And when $d\leq n-7$ and $d(d+1)
\geq 2(n+1)$, then the canonical divisor is big for every $e\geq 1$.

\medskip\noindent
Recall that a sufficient condition for a variety $M$ to be of general
type is that the canonical divisor $K_M$ is big and $M$ has only
canonical singularities.  This raises the hope that in the degree
range above the spaces $\kbm{0,0}{X,e}$ are all of general type.  The
missing ingredient is an analysis of the singularities of
$\kbm{0,0}{X,e}$.  This paper is the result of a ``initial
investigation'' of the singularities of $\kbm{0,0}{X,e}$.  Obviously
much work is still needed to prove that $\kbm{0,0}{X,e}$ contains no
uniruled sweeping subvariety.

\subsection{Detailed summary} \label{subsec-summary} 

\medskip\noindent
The proof of Theorem~\ref{thm-main} is a
deformation-and-specialization argument. The stack
$\Kbm{0,0}{\PP^n,e}$ is smooth, therefore the singularities of $C_d$
come from loci in $\Kbm{0,0}{\PP^n,e}$ over which the fiber dimension
of $\pi_d:C_d \rightarrow \Kbm{0,0}{\PP^n,e}$ jumps.  This defines a
stratification of $\Kbm{0,0}{\PP^n,e}$, and the ``deepest'' stratum
corresponds to the locus $Y \subset \Kbm{0,0}{\PP^n,e}$ parametrizing
multiple covers of lines.  This stratum is a smooth variety, and the
normal bundle of $Y \subset \Kbm{0,0}{\PP^n,e}$ is a certain sheaf of
$(e-1) \times (n-1)$ matrices.

\medskip\noindent
The normal cone of $\pi_d^{-1}(Y) \subset C_d$ is a projective cone
over the normal bundle of $Y \subset \Kbm{0,0}{\PP^n,e}$.  When $d+e
\leq n$ this projective cone is even a projective Abelian cone
associated to a torsion-free sheaf on the normal bundle of $Y \subset
\Kbm{0,0}{\PP^n,e}$.  This torsion-free sheaf is essentially the
direct sum of $d$ copies of the quotient of a bundle of rank $(n-1)$
by the universal $(e-1) \times (n-1)$ matrix mentioned above.  By an
explicit resolution of singularities, one concludes that this
projective Abelian cone is canonical. Then \emph{deformation to the
normal cone} produces a family over $\PP^1$ whose fibers over $\AAA^1$
are all isomorphic to $C_d$ and whose fiber over $\infty$ is the
normal cone of $\pi_d^{-1}(Y)$.

\medskip\noindent
Applying inversion-of-adjunction type results to this family, one
concludes that there exists an open substack $U \subset
\Kbm{0,0}{\PP^n,e}$ containing $Y$ such that $\pi_d^{-1}(U)$ is
canonical.  Moreover the action of the group $\text{GL}_{n+1}$ on
$\PP^n$ induces an action of $\text{GL}_{n+1}$ on
$\Kbm{0,0}{\PP^n,e}$.  The open substack $U$ is
$\text{GL}_{n+1}$-invariant, but also $Y$ intersects the closure of
every orbit of $\text{GL}_{n+1}$ on $\Kbm{0,0}{\PP^n,e}$.  Therefore
$U$ is all of $\Kbm{0,0}{\PP^n,e}$, which proves that $C_d$ is
canonical.

\medskip\noindent
Section ~\ref{sec-discrep}, Section ~\ref{sec-deter}, Section
~\ref{sec-spec} and Section ~\ref{sec-dnc} are all of a foundational
nature, proving basic results about the singularities of the relative
Grassmannian, or \emph{Grassmannian cone}, associated to a
torsion-free coherent sheaf $\mc{E}$ which has (local) projective
dimension $1$.  The main result of Section~\ref{sec-discrep} is
Proposition~\ref{prop-relcan}, which relates the singularities of a
Grassmannian cone $C$ parametrizing rank $r$ locally free quotients of
$\mc{E}$ to the singularities of the pair $(B,r\cdot B_{g-1})$ where
$B_{g-1}$ is the closed subscheme determined by the \emph{Fitting
ideal} of $\mc{E}$.  The main result of Section~\ref{sec-deter} is
Proposition~\ref{prop-cone} which computes the singularities of the
Grassmannian cone of a direct sum of $a$ copies of the cokernel of the
universal $g\times f$ matrix on the affine space of $g \times f$
matrices.  The main result of Section~\ref{sec-spec} is
Corollary~\ref{cor-spec} which applies known results related to
inversion-of-adjunction to prove adjunction type results for a pair
$(B, r\cdot B_{g-1})$ as above.  Section~\ref{sec-dnc} is a review of
the construction of deformation to the normal cone in preparation for
the proof of the main theorem.

\medskip\noindent
In Section~\ref{sec-covers} and Section~\ref{sec-proof} the proof of
the main theorem is given.  Section~\ref{sec-covers} introduces the
closed substack $Y \subset \Kbm{0,0}{\PP^n,e}$ parametrizing multiple
covers of lines; the main result is Proposition~\ref{prop-phid} which
is an analysis of the coherent sheaves used to define $C_d$ when
restricted to a first-order neighborhood of $Y$ in
$\Kbm{0,0}{\PP^n,e}$.  Section~\ref{sec-proof} gives the proof of
Theorem~\ref{thm-main} along the lines discussed above.

\medskip\noindent
In Section~\ref{sec-RSBT}, we use the Reid--Shepherd-Barron--Tai
criterion to prove that the coarse moduli space $\kbm{0,0}{\PP^n,e}$
has only canonical singularities (and in most cases it is even
terminal).  Combining this analysis with the proof of
Theorem~\ref{thm-main}, we show in Section~\ref{sec-conj} that the
inversion-of-adjunction conjecture of Koll\'ar and Shokurov implies
that the coarse moduli space of the stack $C_d$ has only canonical
singularities when either $e\geq 3$ and $d+e \leq n$ or $e=2$ and $d+3
\leq n$.

\medskip\noindent
Finally, in Section~\ref{sec-canKBM}, we compute the expected
canonical divisor on the stack $\Kbm{0,0}{X,e}$, which is the same as
the canonical divisor on the coarse moduli space $\kbm{0,0}{X,e}$ in
many cases.  We show that when $n +1 < d^2 < (n-3)^2$, for most
choices of $e$ the expected canonical divisor of $\kbm{0,0}{X,e}$ is
big.

\medskip\noindent
\textbf{Acknowledgments}  
This paper is a continuation of ~\cite{HRS2} and ~\cite{HS2}.  My
greatest debt is to my coauthors Joe Harris and Mike Roth.  I am also
grateful for useful conversations with Jiun-Cheng Chen, A. Johan de
Jong, Mircea Musta\c{t}\v{a}, and especially J\'{a}nos Koll\'{a}r.  I
was supported by NSF Grant DMS-0201423.

\section{Discrepancies of a Grassmannian cone} \label{sec-discrep}

\medskip\noindent
Let $B$ be a Noetherian scheme which is connected, normal, and
$\QQ$-Gorenstein of pure dimension $b$.  Let $\phi:\mc{G} \rightarrow
\mc{F}$ be a morphism of locally free $\OO_B$-modules of rank $g$ and
$f$ respectively such that the cokernel $\mc{E} = \text{Coker}(\phi)$
has generic rank $e=f-g$.  In this section all results are of a local
nature on $B$.  So the results apply equally as well to a coherent
sheaf $\mc{E}$ which has local projective dimension $1$ (in the sense
of ~\cite[p. 280]{Ma}).

\medskip\noindent
\begin{notat} \label{notat-C} 
Denote by $\text{det}(\mc{E})$ the invertible sheaf
$\text{det}(\mc{F}) \otimes_{\OO_B} \text{det}(\mc{G})^\vee$.  Let $r$
be an integer, $1\leq r \leq e$, and denote by the pair $(\pi:C
\rightarrow B, \alpha:\pi^* \mc{E} \rightarrow \mc{Q})$ the relative
Grassmannian cone over $B$ parametrizing rank $r$ locally free
quotients of $\mc{E}$.  Denote by $\OO_C(1)$ the invertible sheaf on
$C$, $\text{det}(\mc{Q})$.  Denote by the pair $(\rho:C' \rightarrow
B, \beta:\rho^* \mc{F} \rightarrow \mc{Q}')$ the relative Grassmannian
bundle over $B$ parametrizing rank $r$ locally free quotients of
$\mc{F}$.  Denote by $\OO_{C'}(1)$ the invertible sheaf on $C'$,
$\text{det}(\mc{Q}')$.  The surjection $\pi^* \mc{F} \rightarrow \pi^*
\mc{E} \xrightarrow{\alpha} \mc{Q}$ induces a morphism of $B$-schemes
which we denote $h:C \rightarrow C'$.
\end{notat}

\medskip\noindent
Comparing universal properties, it is clear that $h$ is a closed
immersion whose ideal sheaf is the image of the composite morphism
\begin{equation}
\rho^* \mc{G} \otimes (Q')^\vee \xrightarrow{\phi
  \otimes \beta^\dagger}  \pi^*\mc{F} \otimes_{\OO_{C'}} \rho^*\mc{F}^\vee
  \xrightarrow{\text{Trace}} \OO_{C'}
\end{equation}
which we denote by $\gamma$.

\medskip\noindent
Observe that if $\mc{E}$ is locally free, then $\pi:C \rightarrow B$
is Zariski locally a Grassmannian bundle, but in general the fiber
dimension is not necessarily constant.  The case of interest in the
remainder of the paper is $r=1$, in which case $\pi:C \rightarrow B$
is a \emph{projective Abelian cone} (paraphrasing the notation of
\cite{F} and \cite{BF}).  But the following results hold for arbitrary
$r$.  We begin by stating an obvious condition for $C$ to be
irreducible.

\medskip\noindent
\begin{notat} \label{notat-Bk} 
Let $k$ be an integer $k=0,\dots,g$.  Denote by $B_k \subset B$ the
closed subscheme whose ideal sheaf is generated by the $(k+1)\times
(k+1)$-minors of $\phi$, i.e. $B_k$ is the locus where $\phi$ has rank
at most $k$.
\end{notat}  

\medskip\noindent
\begin{lem} \label{lem-codims} 
The coherent sheaf $\mc{E}$ is torsion-free iff
$\text{codim}_B(B_{g-1}) \geq 2$.  The scheme $C$ is irreducible iff
$\text{codim}_B(B_k - B_{k-1}) \geq r(g-k) +1$ for all
$k=0,\dots,g-1$, in which case $C$ has dimension $c=b+r(e-r)$.
Furthermore, the scheme $C$ is regular in codimension $1$ points if
$\text{codim}_B(B_k - B_{k-1}) \geq r(g-k)+2$ for all $k=0,\dots,g-1$.
If $g=0$, all of these conditions are vacuously satisfied.
\end{lem}

\begin{proof}
Torsion sections of the sheaf $\mc{E}$ correspond locally on $S$ to
sections of $\mc{F}$ which are generically in the image of $\mc{G}$.
Since $S$ is normal and since $\mc{G}$ is locally free, the image of
$\mc{G}$ in $\mc{F}$ equals the intersection of its localization at
all codimension $1$ points of $S$.  Therefore a section of $\mc{F}$
corresponds to a torsion-free section of $\mc{E}$ iff its localization
at all codimension $1$ points of $S$ is torsion-free in $\mc{E}$.
Therefore $\mc{E}$ is torsion-free iff $\text{codim}_B(B_{g-1}) \geq
2$.

\medskip\noindent
Of course $C'$ has pure dimension $b+r(f-r)$ at every point.  Consider
the morphism $\gamma$.  Since the rank of $\rho^* \mc{G} \otimes
(\mc{Q}')^\vee$ is $g\cdot r$, the dimension of $C$ at every point is
at least $c:= b + r(f-r) - g\cdot r = b + r(e-r)$.

\medskip\noindent
The restriction of $\pi$ over the locally closed subscheme
$B_k-B_{k-1}$ is proper and smooth of relative dimension $r\lt( (f-k)
- r \rt)$ and has geometrically irreducible (and nonempty) fibers.  In
particular, the preimage of $B-B_{g-1}$ is normal and irreducible of
dimension $b+r(e-r)$, i.e. $c$.  Therefore to prove that $C$ is
irreducible, it suffices to prove that for each $k=0,\dots, g-1$, the
dimension of $\pi^{-1}(B_k - B_{k-1})$ is at most $c-1$.  Conversely,
if any of these sets has dimension $c$ or greater, then it is not in
the closure of $\pi^{-1}(B-B_{g-1})$ which proves that $C$ is
reducible.  Therefore $C$ is irreducible iff $\text{dim} \pi^{-1}(B_k
- B_{k-1}) \leq c-1$.  But the dimension of this set is clearly $c -
\lt[ \text{codim}_B \lt(B_k - B_{k-1} \rt) -r(g-k) \rt]$.  So $C$ is
irreducible iff $\text{codim}_B \lt(B_k - B_{k-1} \rt) \geq r(g-k) +
1$ for all $k=0,\dots, g-1$.

\medskip\noindent
Now suppose that in fact $\text{codim}_B \lt(B_k - B_{k-1} \rt) \geq
r(g-k) + 2$ for all $k=0,\dots,g-1$.  The preimage $\pi^{-1}(B -
B_{g-1})$ is normal.  So to prove that $C$ is regular in codimension
$1$ points, it suffices to prove that $C$ is regular in codimension
$1$ points which are contained in one of the subsets $\pi^{-1}(B_k -
B_{k-1})$ for $k=0,\dots, g-1$.  But our inequality guarantees that
each of the sets $\pi^{-1}(B_k - B_{k-1})$ has codimension at least
$2$ in $C$, therefore there are no such codimension $1$ points.
\end{proof}

\begin{rmk} \label{rmk-codims} 
It can happen that $\phi$ satisfies the second inequality so that $C$
is irreducible, and yet $C$ is not regular in codimension $1$, e.g. on
$\AAA^2$ consider the morphism $\phi:\OO_{\AAA^2} \rightarrow
\OO_{\AAA^2}^{\oplus 2}$ with matrix $(x^2,y^2)^\dagger$.  On the
other hand, the third inequality is not a necessary condition,
e.g. see Proposition~\ref{prop-cone} below.
\end{rmk}

\medskip\noindent
\begin{hyp} \label{hyp-tf} 
Unless explicitly stated otherwise, the coherent sheaf $\mc{E}$ is
torsion-free, i.e. $\text{codim}_B(B_{g-1}) \geq 2$.
\end{hyp}

\medskip\noindent
\begin{lem} \label{lem-discrep} 
Suppose that $C$ has pure dimension $c=b+r(e-r)$.
\begin{enumerate}
\item 
If $B$ is Cohen-Macaulay, then $C$ is Cohen-Macaulay.  If also $C$ is
regular in codimension $1$, then $C$ is normal.
\item 
If $B$ is Gorenstein, then $C$ is Gorenstein.
\item 
The morphism $\pi$ admits a relative dualizing complex of the form
$\omega_\pi[r(e-r)]$ where $\omega_\pi$ is
the invertible sheaf $\pi^*\text{det}(\mc{E})^{\otimes r}
\otimes_{\OO_C} \OO_C(-e)$.  
\item 
If $C$ is normal, then $C$ is $\QQ$-Gorenstein and the $\QQ$-Cartier
divisor class $K_C$ equals $\pi^*K_B + K_\pi$ where $K_\pi$ is the
divisor class of $\omega_\pi$.
\end{enumerate}
\end{lem}

\medskip\noindent
\begin{proof}
By assumption $h(C)$ has pure codimension $g\cdot r$ in $C$, which
equals the rank of $\mc{G}\otimes_{\OO_{C'}} (Q')^\vee$.  By
~\cite[Thm. 17.3 and Thm. 17.4]{Ma}, if $B$ is Cohen-Macaulay then
$h(C)$ is Cohen-Macaulay.  In this case it follows from Serre's
criterion ~\cite[Thm. 23.8]{Ma} that if $C$ is regular in codimension
$1$, then $C$ is normal.  Using ~\cite[Exer. 18.1]{Ma}, if $B$ is
Gorenstein then $h(C)$ is Gorenstein.

\medskip\noindent
The morphism $\rho$ is smooth and has a relative dualizing complex
$\omega_\rho[r(f-r)]$ where $\omega_\rho$ is the invertible sheaf
$\rho^*\text{det}(\mc{F})^{\otimes r} \otimes_{\OO_{C'}}
\OO_{C'}(-f)$.  The morphism $h$ is a regular embedding and has a
relative dualizing complex $\omega_h[-rg]$ where $\omega_h$ is the
pullback of the invertible sheaf $\textit{Ext}^{rg}_{\OO_{C'}}
(h_*\OO_C, \OO_{C'})$.  Forming the Koszul complex associated to the
sheaf map $\gamma$, we conclude that $\omega_h \cong
\pi^*\text{det}\lt( \mc{G}^\vee \rt)^{\otimes r} \otimes_{\OO_C}
\OO_{C}(g)$.  Therefore the composite $\pi = \rho\circ h$ has a
relative dualizing complex $\omega_\pi[r(e-r)]$ where $\omega_\pi$ is
an invertible sheaf isomorphic to $\pi^*\text{det}(\mc{E})^{\otimes r}
\otimes \OO_C(-e)$.

\medskip\noindent
Suppose that $C$ is normal.  Let $U \subset C$ be the smooth locus of
$C$ and let $V\subset C'$ denote the smooth locus of $C'$.  Since $h$
is a regular embedding, we have that $U \subset h^{-1}(V)$.  Also,
since $\rho$ is smooth, we have that $V$ is just $\rho^{-1}(W)$ where
$W\subset B$ is the smooth locus.  Now $\omega_B|_W$ is isomorphic to
$\OO_B(K_B)|_W$ and $\omega_{C'}|_{\rho^{-1}(W)}$ is isomorphic to
$\rho^*\OO_B(K_B) \otimes_{\OO_{C'}} \omega_\rho$.  By the same
reasoning as above, we have that $\omega_C|_V \cong
h^*\rho^*\OO_B(K_B) \otimes \omega_\pi$.  But of course $\omega_C|_V$
is isomorphic to $\OO_C(K_C)|_V$.  Therefore the $\QQ$-Weil divisor
class $K_C$ is equal to $\rho^*K_C + K_\pi$ where $K_\pi$ is the
divisor class of $\omega_\pi$.  Since this is a $\QQ$-Cartier divisor
class, $C$ is $\QQ$-Gorenstein.
\end{proof}

\medskip\noindent
\begin{cor} \label{cor-discrep} 
Let $Y \subset B$ be a closed subscheme which is a regular embedding
of pure codimension $\text{codim}_B(Y)$, i.e. for every closed point
$p\in Y$, the ideal sheaf $\mc{I}_p \subset \OO_{B,p}$ is generated by
a regular sequence of length $\text{codim}_B(Y)$.
\begin{enumerate}
\item
If $C \times_B Y$ has the expected dimension $b-\text{codim}_B(Y) +
r(e-r)$, then there exists an open subset $U\subset B$ containing $Y$
such that $C \times_B U$ has the expected dimension $b + r(e-r)$.
\item
If also $C\times_B Y$ is irreducible, then we can choose $U$ so that
$C \times_B U$ is irreducible.
\item
If also $C\times_B Y$ is normal and $B$ is Cohen-Macaulay, then we can
choose $U$ so that $C\times_B U$ is normal.
\end{enumerate}
\end{cor}

\medskip\noindent
\begin{proof}
Let $C_i \subset C$ be an irreducible component of $C$ which has
nonempty intersection with $C\times_B Y$.  The dimension of $C$ is at
least $b+r(e-r)$, and we are trying to prove that it is exactly
$b+r(e-r)$.  Since $Y \subset B$ is a regular embedding locally
defined by a regular sequence of length $\text{codim}_B(Y)$, it
follows by Krull's Hauptidealsatz that $\text{dim}(C_i \times_B Y)
\geq \text{dim}(C_i) - \text{codim}_B(Y)$.  On the other hand, since
$C_i \times_B Y$ is a closed subscheme of $C \times_B Y$,
$\text{dim}(C_i \times_B Y) \leq \text{dim}(C \times_B Y) = b+r(e-r) -
\text{codim}_B(Y)$.  Therefore we conclude that $\text{dim}(C_i) = b +
r(e-r)$.  So for any irreducible component $C_i \subset C$ whose
dimension is larger than $b+r(e-r)$, we have that $\pi(C_i) \cap Y =
\emptyset$.  We define $U$ to be the complement of the finitely many
closed sets $\pi(C_i)$ as above.  This proves Item $(1)$.

\medskip\noindent
Suppose now that also $C \times_B Y$ is irreducible.  By
Lemma~\ref{lem-codims}, for each $k=0,\dots,g-1$, we have that
$\text{codim}_Y(Y_k - Y_{k-1}) \geq r(g-k) + 1$.  Now $Y_k = B_k \cap
Y$.  So again by Krull's Hauptidealsatz, for every irreducible
component $(B_k)_i$ of $B_k$ which intersects $Y$, we have that
$\text{dim}((B_k)_i) \leq \text{dim}(Y_k) + \text{codim}_B(Y)$,
i.e. $\text{codim}_B((B_k)_i) \geq \text{codim}_Y(Y_k)$.  Now we
shrink the $U$ from the last paragraph by taking the complement of the
finitely many irreducible components $(B_k)_i$ which don't intersect
$Y$ and which have the wrong codimension.  Then for every $(B_k)_i$
which intersects $U$, we have that $\text{codim}_B((B_k)_i) \geq
\text{codim}_Y(Y_k) \geq r(g-k) + 1$.  Therefore by
Lemma~\ref{lem-codims}, we have that $C \times_B U$ is irreducible.

\medskip\noindent
Finally, suppose that also $C\times_B Y$ is normal and $B$ is
Cohen-Macaulay.  By Item $(1)$ of Lemma~\ref{lem-discrep}, to prove
that $C \times_B U$ is normal it suffices to prove that $C \times_B U$
is regular in codimension one.  Let $(C\times_B U)_{\text{sing}}$ be
the singular locus of $C\times_B U$.  Since $C\times_B Y \subset
C\times_B U$ is a Cartier divisor, every regular point of $C\times_B
Y$ is also a regular point of $C\times_B U$.  Therefore the
intersection of $(C\times_B U)_{\text{sing}}$ with $C\times_B Y$ is
contained in $(C\times_B Y)_{\text{sing}}$.  Since $C\times_B Y$ is
normal, $(C\times_B Y)_{\text{sing}}$ has codimension at least two in
$C\times_B Y$.  So again by Krull's Hauptidealsatz we have that every
irreducible component of $(C\times_B U)_{\text{sing}}$ which
intersects $C\times_B Y$ has codimension at least $2$ in $C\times_B
U$.  So, after shrinking $U$ more, we may assume that $C\times_B U$ is
normal.
\end{proof}

\medskip\noindent
\begin{hyp} \label{hyp-1} 
From now on we will assume that $C$ is irreducible of the expected
dimension $c = b + r(e-r)$.
\end{hyp}

\medskip\noindent
\begin{defn} \label{defn-res} 
A morphism of schemes $u:\widetilde{B} \rightarrow B$ is a
\emph{resolution of} $\mc{E}$ if
\begin{enumerate}
\item $u$ is a birational, proper morphism
\item $\widetilde{B}$ is smooth 
\item the exceptional locus of $u$ is a simple normal crossings
  divisor $E_1 \cup \dots \cup E_k$, and
\item the coherent sheaf $\widetilde{\mc{E}} :=
  u^*\mc{E}/\text{torsion}$ is a locally free
  $\OO_{\widetilde{B}}$-module of rank $e$.
\end{enumerate}
\end{defn}

\medskip\noindent
\begin{notat} \label{notat-res} 
Let $u:\widetilde{B} \rightarrow B$ be a resolution of $\mc{E}$.
Denote by $\widetilde{\mc{G}}$ the kernel of the induced surjective
sheaf map $u^*\mc{F} \rightarrow \widetilde{\mc{E}}$ and denote by
$\widetilde{\phi}: \widetilde{\mc{G}} \rightarrow u^* \mc{F}$ the
induced injection, i.e. $\widetilde{\mc{E}}$ is the cokernel of
$\widetilde{\phi}$.  Denote by the pair $(\text{pr}_1:\widetilde{B}
\times_B C' \rightarrow \widetilde{B}, \text{pr}_2^*\beta:
\text{pr}_2^*\rho^*\mc{F} \rightarrow \text{pr}_2^*\mc{Q}')$ the
base-change of $(\rho:C' \rightarrow B, \beta)$.  Denote by the pair
$(\widetilde{\pi}:\widetilde{C} \rightarrow \widetilde{B},
\widetilde{\alpha}:\widetilde{\pi}^* \widetilde{\mc{E}} \rightarrow
\widetilde{Q})$ the Grassmannian bundle parametrizing rank $r$ locally
free quotients of $\widetilde{\mc{E}}$.  Denote by
$\OO_{\widetilde{C}}(1)$ the invertible sheaf on $\widetilde{C}$,
$\text{det}(\widetilde{Q})$.  The surjection $u^*\mc{F} \rightarrow
\widetilde{\mc{E}}$ induces a closed immersion which we denote
$\widetilde{h}:\widetilde{C} \rightarrow \widetilde{B} \times_B C'$.
\end{notat}

\medskip\noindent
Because the morphism $u^*\mc{F} \rightarrow \widetilde{\mc{E}}$
factors through the pullback $u^*\mc{F} \rightarrow u^*\mc{E}$, the
composition $\text{pr}_2\circ h: \widetilde{C} \rightarrow C'$ factors
through $h$, i.e. there is an induced morphism $v:\widetilde{C}
\rightarrow C$.  Of course $v^* \mc{Q} \cong \widetilde{\mc{Q}}$ and
$\pi\circ v = u\circ \widetilde{\pi}$.

\medskip\noindent
\begin{lem} \label{lem-res} 
The morphism $v: \widetilde{C} \rightarrow C$ is a weak resolution of
singularities, that is
\begin{enumerate}
\item
$v$ is a proper, birational morphism, and
\item $\widetilde{C}$ is nonsingular.
\end{enumerate} 
Moreover, the exceptional locus of $v$ is contained in the divisor
$\widetilde{\pi}^{-1}(E_1 \cup \dots \cup E_k)$.
\end{lem}

\medskip\noindent
\begin{proof}
This is obvious.
\end{proof}

\medskip\noindent
\begin{rmk} \label{rmk-discrep} 
If $B$ is a finite type scheme over an algebraically closed field of
characteristic $0$, then a resolution of $\mc{E}$ exists.  Form the
Grassmannian bundle $(\pi_e:G\rightarrow B,\alpha_e: \pi_e^* \mc{E}
\rightarrow Q)$ parametrizing rank $e$ locally free quotients of
$\mc{E}$.  Since $\mc{E}$ generically has rank $e$, there is an
irreducible component $G_0$ of $G$ such that $\pi_e:G_0 \rightarrow B$
is birational.  Let $Z\subset G_0$ denote the fundamental locus of
$\pi_e^{-1}$.  By ~\cite{Hir}, we can find a log resolution
$\widetilde{B} \rightarrow B$ of the pair $(G_0,Z)$, and the induced
morphism $u:\widetilde{B} \rightarrow B$ will be as above.
\end{rmk}

\medskip\noindent
\begin{lem} \label{lem-logres} 
The morphism $u:\widetilde{B} \rightarrow B$ is a log resolution of
the pair $(B,B_{g-1})$.
\end{lem}

\medskip\noindent
\begin{proof}
The morphism $\phi:\mc{G}\rightarrow \mc{F}$ induces an element
$\bigwedge^g\phi \in \textit{Hom}_{\OO_B}(\bigwedge^g
\mc{G},\bigwedge^g \mc{F})$, i.e. an element in $\bigwedge^g \mc{F}
\otimes_{\OO_B} \lt(\bigwedge^g \mc{G} \rt)^\vee$.  There is an
induced map
\begin{equation}
\text{Id}\otimes \bigwedge^g \phi: \bigwedge^e \mc{F} \rightarrow
\bigwedge^e \mc{F} \otimes_{\OO_B} \bigwedge^g \mc{F} \otimes_{\OO_B}
\lt( \bigwedge^g \mc{G} \rt)^\vee.
\end{equation}
We can compose this map with the wedge product map on $\mc{F}$ to get
a map $\bigwedge^e \mc{F} \rightarrow \text{det}(\mc{E})$.  Consider
the restriction of this map to the torsion-free subsheaf $\mc{G}\cdot
\bigwedge^{e-1} \mc{F} \subset \bigwedge^e \mc{F}$.  On the generic
point of $B$, it is clear that this map is the zero map.  But a
morphism between torsion-free sheaves on $B$ which is zero at the
generic point of $B$ is the zero map.  So the restriction of the map
to $\mc{G} \cdot \bigwedge^{e-1} \mc{F}$ is zero, which proves that
the map factors through a map $\psi: \bigwedge^e \mc{E} \rightarrow
\text{det}(\mc{E})$.

\medskip\noindent 
Since $\mc{E}$ is torsion-free, $B_{g-1}$ has codimension at least $2$
in $B$.  So $\psi$ is an isomorphism in codimension $1$.  Denote by
$\mc{I} \subset \OO_B$ the unique ideal sheaf so that
$\text{Image}(\psi)$ equals $\mc{I} \cdot \text{det}(\mc{E})$.  In
other words, $\mc{I}$ is the $e^{\text{th}}$ \emph{Fitting ideal} of
$\phi$, i.e. the ideal sheaf generated by the $g\times g$ minors of
$\phi$.  This is precisely the ideal sheaf of the subscheme $B_{g-1}$.
It follows that $\lt( u^*\bigwedge^e \mc{E} \rt)/\text{torsion}$ is
just $u^{-1}\mc{I}\cdot u^*\text{det}(\mc{E})$.

\medskip\noindent
The surjection $u^*\mc{E} \rightarrow \widetilde{\mc{E}}$ induces a
surjection $\lt(u^*\bigwedge^e \mc{E}\rt)/\text{torsion} \rightarrow
\bigwedge^e \widetilde{\mc{E}}$.  It is easy to see that this
surjection is in fact an isomorphism.  Therefore we have a canonical
isomorphism $u^{-1}\mc{I}\cdot u^*\text{det}(\mc{E}) \cong \text{det}
(\widetilde{\mc{E}})$.  In particular, the pullback ideal sheaf
$u^{-1}\mc{I}\cdot \OO_{\widetilde{B}}$ is a Cartier divisor.
Moreover, this divisor is a subdivisor of the simple normal crossings
divisor $E_1 \cup \dots \cup E_k$, and so it is also a simple normal
crossings divisor.  So we have proved the lemma.
\end{proof}

\medskip\noindent
By definition, the log discrepancies $a(E_i;B,r\cdot B_{g-1})$ of
$(B,r\cdot B_{g-1})$ along the divisors $E_1,\dots,E_k \subset
\widetilde{B}$ are defined by
\begin{equation}
K_{\widetilde{B}} - u^*K_B - r\cdot u^{-1}(B_{g-1}) = \sum_{i=1}^k \lt(
a(E_i;B,r\cdot B_{g-1})
-1 \rt) E_i
\end{equation}
where $u^{-1}(B_{g-1}) \subset \widetilde{B}$ is defined to be the
closed subscheme corresponding to $u^{-1}\mc{I} \cdot
\OO_{\widetilde{B}}$.  By the isomorphism above, $u^{-1}\mc{I} \cong
\text{det} (\widetilde{\mc{E}}) \otimes_{\OO_{\widetilde{B}}} u^*
\text{det}(\mc{E})^\vee$.  On the other hand, applying
Lemma~\ref{lem-discrep} to both $C$ and $\widetilde{C}$, we conclude
that
\begin{equation}
K_{\widetilde{C}} - v^* K_C = \widetilde{\pi}^* \lt( K_{\widetilde{B}} -
u^* K_B + r\cdot C_1(\text{det}(\widetilde{\mc{E}}) ) -
r\cdot u^*C_1(\text{det}(\mc{E}) ) \rt).
\end{equation}
So we have the following lemma.

\medskip\noindent
\begin{lem}  \label{lem-relcan} 
The relative canonical divisor of $v:\widetilde{C} \rightarrow C$ is
  equal to the following divisor
\begin{equation}
K_{\widetilde{C}} - v^*K_C = \sum_{i=1}^k \lt( a(E_i;B,r\cdot B_{g-1}) - 1 \rt)
\widetilde{\pi}^* E_i.
\end{equation}
\end{lem}

\medskip\noindent
\begin{prop} \label{prop-relcan} 
The pair $(C,\emptyset)$ is log canonical (resp. Kawamata log
terminal, canonical) iff the pair $(B,r\cdot B_{g-1})$ is log
canonical (resp. Kawamata log terminal, canonical).
\end{prop}

\medskip\noindent
\begin{proof}
Of course the total discrepancy of $(C,\emptyset)$,
$\text{totaldiscrep}(C,\emptyset)$, is the minimum of $0$ and the
discrepancy of $C$, $\text{discrep}(C,\emptyset)$.  So $(C,\emptyset)$
is log canonical (resp. Kawamata log terminal, canonical) iff
$\text{totaldiscrep}(C,\emptyset)$ is $=0$ (resp. $> -1$, $\geq -1$).
By a standard argument (c.f. ~\cite[Cor. 2.32]{KM} and
~\cite[Prop. 1.3(iv)]{EMY}), $\text{totaldiscrep}(C,\emptyset) =
\text{totaldiscrep}(\widetilde{C}, -K_{\widetilde{C}/C})$.  Since
$\widetilde{C}$ is smooth and $\widetilde{\pi}^*(E_1 \cup \dots \cup
E_k)$ is a simple normal crossings divisor, the total discrepancy
(c.f. ~\cite[Defn. 2.28]{KM} and ~\cite[Defn.  2.34]{KM}) is given by
a combinatorial formula in the coefficients $a(\widetilde{\pi}^*(E_i);
\widetilde{C}, -K_{\widetilde{C}/C})$.  But these are the same as the
coefficients $a(E_i;B,r\cdot B_{g-1})$.  So the total discrepancy of
$(C,\emptyset)$ equals the minimum of $0$ and the integers $a(E_i;B,
r\cdot B_{g-1})-1$.  Moreover, by assumption all of the divisors $E_i$
are exceptional for $u$.  Therefore the minimum of the integers
$a(E_i;B, r\cdot B_{g-1})-1$ is the discrepancy of $(B, r\cdot
B_{g-1})$.  Therefore $(C,\emptyset)$ is log canonical, etc. iff
$(B,r\cdot B_{g-1})$ is log canonical, etc.
\end{proof}

\medskip\noindent
\begin{rmk} \label{rmk-res} 
If each log discrepancy $a(E_i;B,r\cdot B_{g-1})$ is different than
$1$, then it follows that the exceptional locus of $v:\widetilde{C}
\rightarrow C$ is all of $\widetilde{\pi}^{-1}(E_1 \cup \dots \cup
E_k)$ and $v$ is a strong resolution of singularities.  In this case
one can also conclude that $C$ is terminal iff $(B,r\cdot B_{g-1})$ is
terminal.  But if any log discrepancies equal $1$ this can fail;
e.g. if $B=\AAA^2$ and $\phi:\OO_{\AAA^2} \rightarrow
\OO_{\AAA^2}^{\oplus 2}$ is the map with matrix $(x,y)^\dagger$, then
for $r=1$ the cone $C$ is the blowing up of $\AAA^2$ in the origin.
So $C$ is smooth, and thus terminal.  But the pair $(\AAA^2,\{0\})$ is
only canonical.
\end{rmk}

\medskip\noindent
\subsection{Further manipulations} \label{subsec-further}

\medskip\noindent 
The following results are straightforward and aren't actually used in
the rest of the paper, but it is natural to state them here.  In this
subsection we will denote the Grassmannian cone $C$ by $C_r$ to
emphasize the integer $r$.  Proofs are left to the reader.

\medskip\noindent
\begin{lem} \label{lem-further1} 
   Let $1\leq s < r \leq e$.  If $C_r$ is normal of pure dimension
  $b+r(e-r)$, then also $C_s$ is normal of pure dimension $b+s(e-s)$.
  Moreover, if $C_r$ is log canonical (resp. Kawamata log terminal,
  canonical, terminal) then also $C_s$ is log canonical
  (resp. Kawamata log terminal, canonical, terminal).
\end{lem}

\medskip\noindent 
Let $\phi^{(i)}:\mc{G}^{(i)} \rightarrow \mc{F}^{(i)}$, $i=1,\dots,N$
be a sequence of injective morphisms of locally free sheaves such that
for each $i=1,\dots, N$ the cokernel $\mc{E}^{(i)}$ has generic rank
$e^{(i)} = f^{(i)} - g^{(i)}$.  Let $r^{(1)},\dots,r^{(N)}$ be a
sequence of integers with $1\leq r^{(i)} \leq e^{(i)}$.  Let
$\pi^{(i)}:C^{(i)} \rightarrow B$ denote the Grassmannian cone of rank
$r^{(i)}$ locally free quotients of $\mc{E}^{(i)}$ and let $\pi:C
\rightarrow B$ denote the fiber product $C^{(1)} \times_B \dots
\times_B C^{(N)}$.

\medskip\noindent
\begin{lem} \label{lem-further2} 
Suppose that $C$ has pure dimension $c=b+\sum_i
  r^{(i)}(e^{(i)}-r^{(i)})$.  Then Lemma~\ref{lem-discrep} applies to
  $\pi:C \rightarrow B$ where now the dualizing complex is
  $\omega_\pi[\sum_i r^{(i)}(e^{(i)}-r^{(i)})]$ with $\omega_\pi$
  equal to the tensor product of $\pi^*\lt[
  \text{det}(E^{(1)})^{\otimes r^{(1)}} \otimes \dots \otimes
  \text{det}(E^{(N)})^{\otimes r^{(N)}} \rt]$ with
  $\OO_{C^{(1)}}(-e^{(1)})\otimes \dots \otimes
  \OO_{C^{(N)}}(-e^{(N)})$.
\end{lem}

\medskip\noindent
\begin{lem} \label{lem-further3} 
Suppose that $C$ is normal of pure dimension $c$.  For each
  $i=1,\dots,N$ let $Z^{(i)}$ denote the closed subscheme associated
  to the $e^{(i)}$ Fitting ideal of $\phi^{(i)}$.  Then for each
  divisor $E$ of $K(B)$, the log discrepancy
  $a(\pi^{-1}(E);C,\emptyset)$ equals the log discrepancy
  $a(E;B,\sum_i r^{(i)} Z^{(i)})$.  And $(C,\emptyset)$ is log
  canonical (resp. Kawamata log terminal, canonical) iff $(B,\sum_i
  r^{(i)} Z^{(i)})$ is log canonical (resp. Kawamata log terminal,
  canonical).  Also, for every subset $I \subset \{1, \dots, N \}$,
  the fiber product $C^{(I)} = \prod_B (C^{(i)} | i \in I )$ is normal
  of pure dimension $b + \sum_{i\in I} r^{(i)} Z^{(i)}$, and if $C$ is
  log canonical (resp. Kawamata log terminal, canonical), then
  $C^{(I)}$ is log canonical (resp. Kawamata log terminal, canonical).
\end{lem}

\medskip\noindent
Define $\mc{G}:=\oplus_i \mc{G}_i$, define $\mc{F}:= \oplus_i
\mc{F}_i$, define $\phi:\mc{G} \rightarrow \mc{F}$ to be the direct
sum over $i=1,\dots, N$ of $\phi^{(i)}$ and define $\mc{E}$ to be the
cokernel of $\phi$, i.e. $\mc{E} \cong \oplus_i \mc{E}^{(i)}$.  Denote
by $e$ the sum $\sum_i e^{(i)}$ and let $r$ be an integer $1\leq r
\leq e$.  Define $\pi':C' \rightarrow B$ to be the Grassmannian cone
parametrizing rank $r$ locally free quotients of $\mc{E}$.

\medskip\noindent
\begin{lem} \label{lem-further4} 
\begin{enumerate}
\item The $e^{\text{th}}$ Fitting ideal of $\phi$ is the product
  $\mc{I}^{(1)}\cdot \dots \cdot \mc{I}^{(N)}$.
\item  If $C'$ is normal of pure dimension $b+r(e-r)$, then for every
  subset $I \subset \{1,\dots,N\}$, the Grassmannian bundle $C'_I
  \rightarrow B$ parametrizing rank $r$ locally free quotients of
  $\oplus_{i\in I} \mc{E}^{(i)}$ is normal of pure dimension
  $b+r(\sum_{i\in I} e^{(i)} -r)$. 
\item If, moreover, $C'$ is log canonical (resp. Kawamta log terminal,
  canonical), then also $C'_I$ is log canonical
  (resp. Kawamata log terminal, canonical).
\end{enumerate}
\end{lem}

\section{Log discrepancies of generic determinantal varieties}
\label{sec-deter} 

\medskip\noindent
Let $K$ be a field, not necessarily algebraically closed or of
characteristic zero.  In this section, for convenience, we will work
in the category of $K$-schemes.  The interested reader will see how to
prove all the analogous results over $\SP(\ZZ)$, and thus over an
arbitrary base scheme.

\medskip\noindent
Let $S$ be a $K$-scheme and let $\mc{G}$, $\mc{F}$ be locally free
$\OO_S$-modules of finite rank $g$ and $f$ respectively with $g\leq
f$.  Define $\pi^{(0)}:M^{(0)}(S,\mc{G},\mc{F})\rightarrow S$ to be
the affine bundle $\underline{\text{Spec}}_S
\text{Sym}^*\lt(\textit{Hom}_{\OO_S}(\mc{G},\mc{F})^\vee \rt)$.  When
there is no risk of confusion, we denote $M^{(0)}(S,\mc{G},\mc{F})$ by
$M^{(0)}$.  There is a tautological sheaf map
\begin{equation}
\phi: \mc{G}\otimes_{\OO_S} \OO_{M^{(0)}} \rightarrow \mc{F}
\otimes_{\OO_S} \OO_{M^{(0)}}.
\end{equation}

\medskip\noindent
\begin{defn} \label{defn-gdv} 
For $k=0,\dots, g$, the $k^{\text{th}}$ \emph{generic determinantal
variety} $M^{(0)}_k \subset M^{(0)}$ is defined to be the closed
subscheme of $M^{(0)}$ whose ideal sheaf is generated by the
$(k+1)\times (k+1)$-minors of $\phi$ just as in
Notation~\ref{notat-Bk} (c.f. also ~\cite[Sec. II.2]{ACGH}).  For
notation's sake, we define $M^{(0)}_{-1}$ to be the empty set.  In
particular, $M^{(0)}_0$ is just the zero section of $\pi^{(0)}:M^{(0)}
\rightarrow S$.
\end{defn}

\medskip\noindent
In this section we will compute the log discrepancies of the pair
$(M^{(0)}, M^{(0)}_k)$.  This is straightforward once we have a log
resolution.  The log resolution is obtained by first blowing up
$M^{(0)}_0$, then blowing up the strict transform of $M^{(0)}_1$, the
blowing up the strict transform of $M^{(0)}_2$, etc.

\subsection{The log resolution} \label{subsec-logres} 
\medskip\noindent 
The log resolution of $(M^{(0)}, M^{(0)}_k)$ which we use is the
obvious one: we succesively blow up the strict transforms of the
schemes $M^{(0)}_0, M^{(0)}_1,\dots, M^{(0)}_k$.  Using the action of
the group $\text{GL}(\mc{F}) \times_S \text{GL}(\mc{G})$, it is easy
to prove this does give a log resolution.  For completeness, we go
through the proof in (somewhat tedious) detail.

\medskip\noindent
\begin{rmk} \label{rmk-wonderful} 
The set of subschemes $M^{(0)}_k \subset M^{(0)}$ form a
stratification, and in case the ground field is $\CC$, this is a
\emph{conical stratification} in the sense of ~\cite{MP98} and the
blowing up we construct coincides with the \emph{minimal wonderful
compactification}.  We choose not to follow ~\cite{MP98} for two
reasons: First of all, the log resolution we construct is completely
obvious and exists over a ground field $K$ of arbitrary
characteristic, not just over $\CC$ (and in fact over $\SP(\ZZ)$,
though we don't prove this).  More importantly, for our purposes it is
crucial that the log resolution we construct have the \emph{additional
property} that for each $M^{(0)}_k$ the inverse image of the ideal
sheaf of $M^{(0)}_k$ is an invertible sheaf, i.e. it contains no
embedded points.  This typically fails for the minimal wonderful
compactification associated to a conical stratification, e.g. if one
considers the nodal plane cubic $p\in C$ sitting in $\PP^2$ sitting as
a linear subvariety of $\PP^3$, then $(\{p\}, C - \{p\}, \PP^3 - C)$
is a conical stratification of $\PP^3$ and the inverse image of the
ideal sheaf of $C$ in the minimal wonderful compactification has an
``embedded line'' on the exceptional divisor over $p$.  It would be
interesting to know if there are extra hypotheses of a general nature
which can be added to the definition of a conical stratification so
that the minimal wonderful compactification has the additional
property.
\end{rmk}

\medskip\noindent
In the case that $f=g$, the log resolution we construct is identical
to that in ~\cite{Kausz}.  Moreover in ~\cite{Kausz} it is proved that
the inverse image of the ideal sheaf of $M^{(0)}_k$ is an invertible
sheaf.  However we are mostly interested in the case $f\neq g$, so we
give the full description of the log resolution and proofs of the
basic properties of the log resolution.  We begin by giving a more
precise definition of the sequence of blowing ups mentioned above.

\medskip\noindent
\begin{lem} \label{lem-exist} 
There exists a sequence of schemes $M^{(r)}$ for $r=0,\dots,g$ and
morphisms $u^{(s,r)}:M^{(r)} \rightarrow M^{(s)}$ for each $0\leq s
\leq r \leq g$ with the following properties
\begin{enumerate}
\item For $0\leq t \leq s \leq r \leq g$, we have $u^{(t,s)}\circ
  u^{(s,r)} = u^{(t,r)}$.
\item For $r=0,\dots, g$, the morphism $u^{(0,r)}:M^{(r)}\rightarrow
  M^{(0)}$ is an isomorphism over the open subscheme $M^{(0)} -
  M^{(0)}_{r-1}$.
\item For each $0\leq r \leq k \leq g$, define $M^{(r)}_k \subset
  M^{(r)}$ to be the closure of the pullback by $u^{(0,r)}$ of
  $M^{(0)}_r - M^{(0)}_{k-1}$.  Then, for $r=0,\dots,g-1$ the morphism
  $u^{(r,r+1)}:M^{(r+1)} \rightarrow M^{(r)}$ is the blowing up of
  $M^{(r)}$ along $M^{(r)}_r$.
\end{enumerate}
\end{lem}
  
\begin{proof}  
This is almost tautological.  The one thing that needs to be checked
is that when we construct $M^{(r+1)}$ as the blowing up of $M^{(r)}$
along $M^{(r)}_r$, that the induced map $u^{0,r+1}:M^{(r+1)}
\rightarrow M^{(0)}$ is an isomorphism over $M^{(0)} - M^{(0)}_{r}$.
But this follows immediately from the fact that $u^{r,0}$ is an
isomorphism over $M^{(0)} - M^{(0)}_{r-1}$ and the fact that
$u^{r,r+1}$ is an isomorphism over the preimage under $u^{r,0}$ of
$M^{(0)} - M^{(0)}_r$.
\end{proof}

\medskip\noindent
\begin{notat} \label{notat-Erk} 
Let $(S,\mc{G},\mc{F})$ be a datum with $\text{rank}(\mc{G}) =g$,
$\text{rank}(\mc{F})=f$ (and of course $g\leq f$). For each
$r=1,\dots,k$, denote by $E^{(r)}_{r-1}(S,\mc{G},\mc{F}) \subset
M^{(r)}(S,\mc{G},\mc{F})$ the exceptional divisor of the blowing up
$u^{r-1,r}$.  For $r=2,\dots,k$, and for $k=0,\dots,r-2$, denote by
$E^{(r)}_k\subset M^{(r)}$ the strict transform of $E^{(k+1)}_k$ under
the morphism $u^{k+1,r}:M^{(r)} \rightarrow M^{(k+1)}$.  Clearly for
each $0\leq r< s \leq g$, the exceptional locus of $u^{r,s}:M^{(s)}
\rightarrow M^{(r)}$ is $E^{(s)}_{r} \cup \dots \cup E^{(s)}_{s-1}$.
\end{notat}

\medskip\noindent
\begin{defn} \label{defn-mor} 
Let $(S,\mc{G}',\mc{F}')$ and $(S,\mc{G},\mc{F})$ be data of pairs of
locally free sheaves on $S$ with $g' \leq f'$ and $g\leq f$.  A
\emph{morphism} between the data is a triple $\zeta = (p,q,T)$ where
\begin{enumerate}
\item
$p:\mc{G} \rightarrow \mc{G}'$ is a surjective morphism of
$\OO_S$-modules,
\item 
$q:\mc{F}'\rightarrow \mc{F}$ is an injective morphism of
$\OO_S$-modules whose cokernel is locally free, and
\item
$T:\mc{G} \rightarrow \mc{F}$ is a morphism of $\OO_S$-modules,
\end{enumerate}
and
such that we have direct sum decompositions 
\begin{equation}
\mc{G} = \text{Ker}(T) \oplus \text{Ker}(p),\ \ \mc{F} =
\text{Image}(T) \oplus \text{Image}(q).
\end{equation}
In particular, we have $g-g' = f-f'$.  The \emph{rank} of $\zeta$ is
the common integer $g-g'=f-f'$.
\end{defn}

\medskip\noindent
\begin{lem} \label{lem-mor} 
Let $\zeta = (p,q,T)$ be a morphism $(S,\mc{G}',\mc{F}') \rightarrow
(S,\mc{G},\mc{F})$ of rank $l$.  Then for each $r=0,\dots,g'$ there
exists a morphism of $S$-schemes
\begin{equation}
\tau^{r}(\zeta):\text{GL}(\mc{F}) \times_S
\text{GL}(\mc{G}) \times_S M^{(r)}(S,\mc{G}',\mc{F}') \rightarrow
M^{(r+l)}(S,\mc{G},\mc{F}) 
\end{equation}
and which satisfy the following conditions
\begin{enumerate}
\item 
The composition $u^{0,l}\circ \tau^0(\zeta)$ is a morphism whose image
is contained in $M^{(0)}(S,\mc{G},\mc{F}) -
M^{(0)}_{l-1}(S,\mc{G},\mc{F})$.
\item 
The morphism $\tau^0(\zeta)$ is the unique morphism such that the
composition $u^{0,l}\circ \tau^0(\zeta)$ is the morphism whose
restriction to the fiber over a point $x\in S$ maps a triple
$(\alpha,\beta,L)$ in $\text{GL}(\mc{F}_x) \times \text{GL}(\mc{G}_x)
\times M^{(0)}(x,\mc{G}'_x, \mc{F}'_x)$ to the element $\alpha\circ
(T_x + q_x\circ L \circ p_x)\circ \beta^{-1}$.
\item 
For $0\leq r \leq s \leq g'$, we have $u^{l+r,l+s}\circ
\tau^{s}(\zeta)$ equals $\tau^r(\zeta) \circ \lt(\text{Id} \times
\text{Id} \times u^{r,s}\rt)$; moreover the commutative diagram is a
Cartesian square.
\item 
For each $0\leq r \leq g'$, the morphism $\tau^r(\zeta)$ is
quasi-compact, separated and smooth and the image is equal to the
preimage under $u^{0,l+r}$ of $M^{(0)}-M^{(0)}_{l-1}$.
\item 
For each $0\leq r\leq k \leq g'$ the preimage under $\tau^r(\zeta)$ of
$M^{(r+l)}_{k+l}(S,\mc{G},\mc{F})$ is the closed subscheme
$\text{GL}(\mc{F})\times_S \text{GL}(\mc{G}) \times_S
M^{(r)}_k(S,\mc{G}', \mc{F}')$ of $\text{GL}(\mc{F})\times_S
\text{GL}(\mc{G}) \times_S M^{(r)}(S,\mc{G}',\mc{F}')$.
\item 
For each $0\leq r \leq g'$ and each $0\leq i \leq r-1$ the preimage
under $\tau^r(\zeta)$ of $E^{(r+l)}_{i+l}(S,\mc{G},\mc{F})$ is the
Cartier divisor $\text{GL}(\mc{F})\times_S \text{GL}(\mc{G}) \times_S
E^{(r)}_{i}(S,\mc{G}',\mc{F})$.
\end{enumerate}
\end{lem}

\medskip\noindent
\begin{proof}
We will prove by induction on $r$ that for each $0 \leq r \leq g'$,
there exists a sequence of morphisms $\tau^0(\zeta), \dots,
\tau^r(\zeta)$ satisfying Items $(1)$ through $(6)$.  We begin with
$r=0$.  First we give a more precise definition of the morphism in
Item $(2)$.  We will define a morphism
\begin{equation}
\tau(\zeta): \text{GL}(\mc{F}) \times_S \text{GL}(\mc{G}) \times_S
M^{(0)}(S,\mc{G}',\mc{F}') \rightarrow M^{(0)}(S,\mc{G},\mc{F})
\end{equation}
by giving a natural transformation of the obvious functors represented
by the two schemes and invoking Yoneda's lemma.

\medskip\noindent
Suppose $T$ is any $k$-scheme.  By the universal properties of the
three factors, a morphism of $T$ to $\text{GL}(\mc{F}) \times_S
\text{GL}(\mc{G}) \times_S M^{(0)}(S,\mc{G}',\mc{F}')$ is equivalent
to a morphism $f:T\rightarrow S$ together with a triple
$(\alpha,\beta,L)$ where $\alpha:f^*\mc{F} \rightarrow f^*\mc{F}$ is
an automorphism of $\OO_T$-modules, where $\beta:f^*\mc{G} \rightarrow
f^*\mc{G}$ is an automorphism of $\OO_T$-modules, and where
$L:f^*\mc{G}' \rightarrow f^*\mc{F}'$ is a morphism of
$\OO_T$-modules.  There is an associated morphism of $\OO_T$-modules
$f^*\mc{G} \rightarrow f^*\mc{F}$ by $L'=\alpha\circ \lt( f^*T +
f^*q\circ L f^*p\rt) \circ \beta^{-1}$.  By the universal property of
$M^{(0)}(S,\mc{G},\mc{F})$ the pair $(f,L')$ is equivalent to a
morphism $T\rightarrow M^{(0)}(S,\mc{G},\mc{F})$.  The association
$(f,\alpha,\beta,L)\mapsto (f,L')$ is clearly a natural transformation
of the Yondea functors and so determines a morphism $\tau$.  Also, it
is clear that $\tau$ is a morphism of $S$-schemes.

\medskip\noindent
Because of the direct sum decompositions of $\mc{G}$ and $\mc{F}$
induced by $(p,q,T)$, we have for any point $x\in T$ that the rank of
$L'_x$ equals $\text{rank}(T_{f(x)}) + \text{rank}(L_x)$,
i.e. $l+\text{rank}(L_x)$ which is bigger than $l-1$.  Therefore the
image of $\tau$ is contained in the complement of
$M^{(0)}_{l-1}(S,\mc{G},\mc{F})$.  Since
$u^{0,l}:M^{(l)}(S,\mc{G},\mc{F}) \rightarrow
M^{(0)}(S,\mc{G},\mc{F})$ is an isomorphism over the complement of
$M^{(0)}_{l-1}(S,\mc{G},\mc{F})$, there is a unique morphism
$\tau^0(\zeta)$ such that $u^{0,l}\circ \tau^0(\zeta)$ equals $\tau$.
This proves Item $(1)$ and Item $(2)$.

\medskip\noindent
Next we establish Item $(4)$ for $\tau^0(\zeta)$.  This is equivalent
to the claim that $\tau$ is quasi-compact, smooth and separated.  It
is clear that both $\text{GL}(\mc{F})\times_S \text{GL}(\mc{G})
\times_S M^{(0)}(S,\mc{G}',\mc{F}')$ and $M^{(0)}(S,\mc{G},\mc{F})$
are quasi-compact, smooth and separated over $S$.  Therefore $\tau$ is
quasi-compact, finitely-presented and separated and to show that
$\tau$ is smooth it suffices to check the Jacobian criterion for
fibers of $\tau$ over geometric points of $S$.

\medskip\noindent
Suppose that $x\in S$ is a geometric point and that $(\alpha,\beta,L)$
is a point of $\text{GL}(\mc{F}) \times_S \text{GL}(\mc{G}) \times_S
M^{(0)}(S,\mc{G}',\mc{F}')$ lying above $x$.  Now the smooth group
scheme $\text{GL}(\mc{F}) \times_S \text{GL}(\mc{G})$ acts on both the
domain and target of $\tau$ and $\tau$ is equivariant for this action.
Therefore it suffices to check the Jacobian criterion at one
representative point of each orbit, i.e. we may assume that
$\alpha=\text{id}_{\mc{F}}$ and that $\beta = \text{id}_{\mc{G}}$.
The Zariski tangent space to the fiber of $\text{GL}(\mc{F}) \times_S
\text{GL}(\mc{G}) \times_S M^{(0)}(S,\mc{G}',\mc{F}')$ at any point is
canonically identified with the $\kappa(x)$-vector space
$\text{Hom}(\mc{F}_x,\mc{F}_x) \times \text{Hom}(\mc{G}_x,\mc{G}_x)
\times \text{Hom}(\mc{G}'_x,\mc{F}'_x)$.  Similarly the Zariski
tangent space to the fiber of $M^{(0)}(S,\mc{G},\mc{F})$ is
canonically identified with $\text{Hom}(\mc{G}_x,\mc{F}_x)$.  And
$d\tau$ maps a triple $(\alpha_1, \beta_1, L_1)$ to the element
$\alpha_1\circ(T_x+q_x\circ L\circ p_x) + q_x\circ L_1 \circ p_x +
(T_x+q_x\circ L\circ p_x) \circ \beta_1$.

\medskip\noindent
Let $d=\text{rank}(L)$.  We can choose ordered bases for $\mc{G}_x$
and $\mc{F}_x$ with respect to which $T_x$ has the matrix
representation
\begin{equation}
\lt[ \begin{array}{r|r|r}
I_{l,l} & 0_{l,d} & 0_{l,g'-d} \\
\hline 0_{d,l} & 0_{d,d} & 0_{d,g'-d} \\
\hline 0_{f'-d,l} & 0_{f'-d,d} & 0_{f'-d,g'-d}
\end{array} \rt]
\end{equation}
and with respect to which $q_x\circ L \circ p_x$ has the matrix
representation
\begin{equation}
\lt[ \begin{array}{r|r|r}
0_{l,l} & 0_{l,d} & 0_{l,g'-d} \\
\hline 0_{d,l} & I_{d,d} & 0_{d,g'-d} \\
\hline 0_{f'-d,l} & 0_{f'-d,d} & 0_{f'-d,g'-d}
\end{array} \rt]
\end{equation}
For any linear operator $L'\in \text{Hom}(\mc{G}_x,\mc{F}_x)$ consider
the matrix representation of $L'$ with respect to the ordered bases
above:
\begin{equation}
\lt[ \begin{array}{r|r|r}
L'_1 & L'_2 & L'_3 \\
\hline L'_4 & L'_5 & L'_6 \\
\hline L'_7 & L'_8 & L'_9
\end{array} \rt]
\end{equation}
where the block submatrices $L'_i$ have the same dimensions as the
blocks in the matrices of $T$ and $L$.  Then if we denote
\begin{equation}
\alpha_1 = \lt[ \begin{array}{r|r|r}
L'_1 & L'_2 & 0 \\
\hline L'_4 & L'_5 & 0 \\
\hline L'_7 & L'_8 & 0 
\end{array}\rt]
\end{equation}
\begin{equation}
\beta_1 = \lt[ \begin{array}{r|r|r}
0_{l,l} & 0_{l,d} & L'_3 \\
\hline 0_{d,l} & 0_{d,d} & L'_6 \\
\hline 0_{f'-d,l} & 0_{f'-d,d} & 0_{f'-d,g'-d}
\end{array} \rt]
\end{equation}
\begin{equation}
q_x\circ L_1 \circ p_x = \lt[ \begin{array}{r|r|r}
0_{l,l} & 0_{l,d} & 0_{l,g'-d} \\
\hline 0_{d,l} & 0_{d,d} & 0_{d,g'-d} \\
\hline 0_{f'-d,l} & 0_{f'-d,d} & L'_9
\end{array} \rt]
\end{equation}
we have that the pair $(\alpha_1,\beta_1,L_1)$ maps to $L'$ under
$d\tau$.  So $\tau$ satisfies the Jacobian criterion, therefore $\tau$
and $\tau^{0}(\zeta)$ are smooth.  This proves Item $(4)$ for
$\tau^0(\zeta)$.  

\medskip\noindent
Since for $L' = \tau(\alpha,\beta,L)$ we have $\text{rank}(L') = l +
\text{rank}(L)$, it is clear that the preimage under $\tau$ of
$M^{(0)}_{l+k}(S,\mc{G},\mc{F})$ is the closed subscheme
$\text{GL}(\mc{F}) \times_S \text{GL}(\mc{G}) \times_S M^{(0)}_{k}(S,
\mc{G}', \mc{F}')$.  Finally Item $(3)$ and Item $(6)$ are vacuous for
the ``sequence'' of morphisms $\tau^0(\zeta)$.  So we have proved that
for $r=0$, there is a sequence of morphisms $\tau^0(\zeta),\dots,
\tau^r(\zeta)$ satisfying Items $(1)$ through $(6)$.

\medskip\noindent
Now comes the induction step.  For some $r=1,\dots,g'$, suppose that
morphisms $\tau^{0}(\zeta),\dots,\tau^{r-1}(\zeta)$ have been
constructed satisfying the conditions.  Since $\tau^{r-1}(\zeta)$ is
smooth, the fiber product of $\tau^{r-1}(\zeta)$ with the blowing up
$u^{l+r-1,l+r}: M^{(l+r)}(S,\mc{G},\mc{F}) \rightarrow
M^{(l+r-1)}(S,\mc{G},\mc{F})$ of $M^{(l+r-1)}(S,\mc{G},\mc{F})$ along
the closed subscheme $M^{(l+r-1)}_{l+r-1}(S,\mc{G},\mc{F})$ is
canonically isomorphic to the blowing up of $\text{GL}(\mc{F})\times_S
\text{GL}(\mc{G})\times_S M^{(r-1)}(S,\mc{G}',\mc{F}')$ along the
preimage of $M^{(l+r-1)}_{l+r-1}(S,\mc{G},\mc{F})$.  By the induction
hypothesis, the preimage is precisely $\text{GL}(\mc{F}) \times_S
\text{GL}(\mc{G}) \times_S M^{(r-1)}_{r-1}(S,\mc{G}',\mc{F}')$.  So
the base-change of $u^{l+r-1,l+r}$ by $\tau^{r-1}(\zeta)$ is just
$\text{Id} \times \text{Id}\times u^{r-1,r}$.  We define
$\tau^r(\zeta):\text{GL}(\mc{F})\times_S \text{GL}(\mc{G}) \times_S
M^{(r)}(S,\mc{G}',\mc{F}') \rightarrow M^{(l+r)}(S,\mc{G},\mc{F})$ to
be the base-change of $\tau^{r-1}(\zeta)$ by $u^{l+r-1,l+r}$.

\medskip\noindent
By construction of $\tau^r(\zeta)$ and the induction hypothesis, it is
clear that $\tau^0(\zeta),\dots,\tau^r(\zeta)$ satisfy Item $(3)$ of
the lemma.  Since $\tau^r(\zeta)$ is the base-change of a
quasi-compact, separated and smooth morphism, $\tau^r(\zeta)$ is also
quasi-compact, separated and smooth, Item $(4)$ is true.

\medskip\noindent
Next we prove Item $(5)$.  Since $\tau^r(\zeta)$ and
$\tau^{r-1}(\zeta)$ are smooth, and by Item $(3)$ the process of
forming the strict transform by $u^{l+r,l+r-1}$ of a closed subscheme
and then forming the preimage under $\tau^r(\zeta)$ is the same as the
process of first forming the preimage under $\tau^{r-1}(\zeta)$ and
then forming the strict transform under $\text{Id} \times \text{Id}
\times u^{r,r-1}$.  By the induction hypothesis and Item $(5)$, the
preimage under $\tau^{r-1}(\zeta)$ of $M^{(r+l-1)}_k(S,\mc{G},\mc{F})$
equals $ \text{GL}(\mc{F}) \times_S \text{GL}(\mc{G}) \times_S
M^{(r-1)}_k(S,\mc{G}',\mc{F}')$.  The strict transform of this
subscheme under $\text{Id} \times \text{Id} \times u^{r-1,r}$ is
$\text{GL}(\mc{F}) \times_S \text{GL}(\mc{G}) \times_S
M^{(r)}_k(S,\mc{G}',\mc{F}')$.  So Item $(5)$ is satisfied for
$\tau^r(\zeta)$.

\medskip\noindent
Finally we prove Item $(6)$.  As in the last paragraph, for
$i=0,\dots,r-2$ the pullback by $\tau^r(\zeta)$ of the strict
transform by $u^{l+r,l+r-1}$ of $E^{(l+r-1)}_{i+l}(S,\mc{G},\mc{F})$
equals the strict transform of the pullback by $\tau^{r-1}(\zeta)$.
By the induction hypothesis and Item $(6)$, this pullback is
$\text{GL}(\mc{F}) \times_S \text{GL}(\mc{G}) \times_S
E^{(r-1)}_i(S,\mc{G}',\mc{F}')$.  The strict transform of this
subscheme under $\text{Id} \times \text{Id} \times u^{r,r-1}$ is
$\text{GL}(\mc{F}) \times_S \text{GL}(\mc{G}) \times_S
E^{(r)}_i(S,\mc{G}',\mc{F}')$.  Finally, using Item $(3)$, the
preimage under $\tau^r(\zeta)$ of the exceptional divisor of
$u^{r+l,r+l-1}$, i.e. $E^{(r+l)}_{r+l-1}(S,\mc{G},\mc{F})$, equals the
exceptional divisor of $\text{Id} \times \text{Id} \times u^{r,r-1}$,
i.e. $\text{GL}(\mc{F}) \times_S \text{GL}(\mc{G}) \times_S
E^{(r)}_{r-1}(S,\mc{G}', \mc{F}')$.  This proves Item $(6)$ and
finishes the proof of the lemma by induction on $r$.
\end{proof}

\medskip\noindent
\begin{lem} \label{lem-snc} 
Let $V$ be a smooth $K$-scheme and let $D\subset V$ be a simple normal
crossings divisor.  The $K$-scheme $\AAA^1 \times V$ is smooth and the
divisor $D' = (\AAA^1\times D) \cup (\{0\} \times V)$ is a simple
normal crossings divisor in $\AAA^1 \times V$.
\end{lem}

\medskip\noindent
\begin{proof}  This follows immediately from the definition of
  \emph{simple normal crossings divisor}.
\end{proof}

\medskip\noindent
\begin{notat} \label{notat-Ur} 
Denote by $U^{(0)} \subset M^{(0)}(S,\mc{G},\mc{F})$ the complement of
the zero section $M^{(0)}_0(S,\mc{G},\mc{F})$.  For each $r=0,\dots,
g$ denote by $U^{(r)} \subset M^{(r)}$ the open subscheme $U^{(r)} =
\lt(u^{0,r}\rt)^{-1}(U^{(0)})$ and denote by $v^{r,s}:U^{(s)}
\rightarrow U^{(r)}$ the restriction of $u^{r,s}$ to $U^{(s)}$.  For
each $r=0,\dots,g$ and each $i=r,\dots,g$, denote by $\mc{I}^{(r)}_k$
the ideal sheaf of the closed subscheme $M^{(r)}_k \subset M^{(r)}$.
Finally denote by $f:\AAA^1 \times U^{(0)} \rightarrow M^{(0)}$ the
morphism which sends a pair $(\lambda, L)$ in $\AAA^1 \times U^{(0)}$
to the point $\lambda\cdot L \in M^{(0)}$.
\end{notat}  

\medskip\noindent
The preimage under $f$ of $M^{(0)}_0$ is precisely $\{0\} \times
U^{(0)}$, which is a Cartier divisor in $\AAA^1 \times U^{(0)}$.
Therefore, by the universal property of blowing up, there is a unique
morphism $f^{1}:\AAA^1 \times U^{(0)} \rightarrow M^{(1)}$ such that
$u^{0,1}\circ f^{1} = f$.  It is easy to check that $f^{1}:\AAA^1
\times U^{(1)} \rightarrow M^{(1)}$ is a $\mathbb{G}_m$-torsor where
$\mathbb{G}_m$ acts on $\AAA^1 \times U^{(1)}$ by $\mu\cdot
(\lambda,L) = (\mu\cdot \lambda, \mu^{-1}\cdot L)$.  In particular,
$f^{1}$ is smooth and surjective.  The preimage under $f^{1}$ of
$E^{(1)}_0$ is $\{0\} \times U^{(1)}$.  The next lemma proves that for
$k =1,\dots, g-1$, the preimage under $f^{1}$ of $M^{(1)}_k$ is
$\AAA^1 \times \lt(U^{(1)}\cap M^{(1)}_k \rt)$.

\medskip\noindent
\begin{lem} \label{lem-pref} 
Define $\mc{G}^{(1)}$ to be the locally free sheaf
$(u^{0,1})^*\mc{G}(E^{(1)}_0)$ and let $(u^{0,1})^*\mc{G} \rightarrow
\mc{G}^{(1)}$ denote the canonical sheaf map.
\begin{enumerate}
\item 
There is a factorization $\phi^{(1)}:\mc{G}^{(1)} \rightarrow
(u^{0,1})^*\mc{F}$ of $(u^{0,1})^*\phi$.
\item 
The pullback of $\phi^{(1)}$ by $f^1$ is canonically isomorphic to the
pullback $\text{pr}_2^* \phi$ of the restriction of $\phi$ to $U^{(0)}
= U^{(1)}$.
\item 
For every geometric point $x$ of $M^{(1)}$,
$\text{rank}(\phi^{(1)}|_x) \geq 1$.
\item 
For each $k=1,\dots, g-1$, the inverse image ideal sheaf
$(u^{0,1})^{-1}(\mc{I}^{(0)}_k)$ equals $\mc{I}^{(1)}_k\cdot
\OO_{M^{(1)}}\lt( -(k+1)E^{(1)}_0 \rt)$.
\item 
For each $k=1,\dots, g-1$, the preimage under $f^{(1)}$ of $M^{(1)}_k$
equals $\AAA^1 \times \lt(U^{(1)} \cap M^{(1)}_k \rt)$.
\end{enumerate}
\end{lem}

\begin{proof}
The restriction of $(u^{0,1})^*\phi$ to $E^{(1)}_0$ is the zero map,
therefore it factors through the elementary transform up of
$(u^{0,1})^*\mc{G}$, i.e. it factors through a morphism
$\phi^{(1)}:\mc{G}^{(1)} \rightarrow (u^{0,1})^*\mc{F}$.  This
establishes Item $(1)$.

\medskip\noindent
Now $(f^1)^* \mc{G}^{(1)}$ equals $\text{pr}_2^*\mc{G}({0}\times
U^{(1)})$, which is canonically isomorphic to $\text{pr}_2^*\mc{G}$.
Via this isomorphism, the pullback $(f^1)^* \phi^{(1)}$ equals
$\text{pr}_2^* \phi$.  Since $\phi|_{U^{(0)}}$ has rank at least $1$
at all geometric points, the same is true of $\text{pr}_2^* \phi$.
Therefore $\phi^{(1)}$ has rank at least $1$ at all geometric points,
i.e. Item $(3)$ is true.

\medskip\noindent
One can check that two ideal sheaves are equal after a faithfully flat
base change.  Since $f^1$ is faithfully flat, it follows that to prove
both Item $(4)$ and Item $(5)$, it suffices to prove that the inverse
image ideal sheaf of $\mc{I}^{(0)}_k$ in $\AAA^1 \times U^{(1)}$
equals
\begin{equation}
\text{pr}_2^{-1}(\mc{I}^{(0)}) \cdot \OO_{\AAA^1 \times
  U^{(1)}}(-(k+1)\{0\} \times U^{(1)}).
\end{equation}
Let $t$ denote the coordinate on $\AAA^1$.  Then the preimage under
$f^1$ of $(u^{0,1})^*\phi$ is precisely the matrix $t\cdot
\text{pr}_2^*\phi$.  Therefore the ideal sheaf generated by the
$(k+1)\times (k+1)$-minors of this matrix is just $t^{k+1}$ times the
ideal generated by $(k+1) \times (k+1)$-minors of $\text{pr}_2^*\phi$,
i.e. the inverse image ideal sheaf of $\mc{I}^{(0)}_k$ under
$u^{0,1}\circ f^1$ is as above.  This finishes the proof of the lemma.
\end{proof}

\medskip\noindent
\begin{lem} \label{lem-f} 
For each $r=1,\dots,g$ there exists a morphism of $S$-schemes $f^{r}:
\AAA^1 \times U^{(r)} \rightarrow M^{(r)}$ and which satisfy the
following conditions
\begin{enumerate}
\item 
The morphism $u^{0,1}\circ f^{1}$ is a morphism whose image is
contained $U^{(0)}$.
\item 
The morphism $f^{1}$ is the unique morphism such that $u^{1,0}\circ
f^{1} = f$.
\item 
For $1\leq r \leq s \leq g$ we have $u^{r,s} \circ f^{s} = f^{r} \circ
\lt( \text{Id} \times v^{r,s} \rt)$; moreover the commutative diagram
is a Cartesian square.
\item 
For each $1 \leq r \leq g$, the morphism $f^r$ is a
$\mathbb{G}_m$-torsor, in particular it is surjective and smooth.
\item 
For each $1 \leq r \leq k \leq g$, the preimage under $f^r$ of
$M^{(r)}_k$ is the closed subscheme $\AAA^1 \times \lt(U^{(r)} \cap
M^{(r)}_k \rt)$.
\item 
For each $2 \leq r \leq g$ and $1 \leq k \leq r-1$, the preimage under
$f^r$ of $E^{(r)}_k$ is $\AAA^1 \times \lt(U^{(r)} \cap E^{(r)}_k
\rt)$.  And for each $1\leq r \leq g$, the preimage under $f^r$ of
$E^{(r)}_0$ is $\{0\} \times U^{(r)}$.
\end{enumerate}
\end{lem}

\begin{proof}
Item $(1)$ is trivial and is only included to maintain symmetry with
Lemma~\ref{lem-mor}.  Item $(2)$ follows from the construction of
$f^{(1)}$ above.  We will prove by induction on $r$ that for each
$1\leq r \leq g$ there exists a sequence of morphisms $f^{1},\dots,
f^{r}$ satisfying Items $(1)$ through $(6)$.  This is already
established for $r=1$, where Item $(5)$ follows from Item $(5)$ of
Lemma~\ref{lem-pref}.

\medskip\noindent
Now comes the induction step.  For some $r=1,\dots, g$, suppose that
morphisms $f^1,\dots,f^{r-1}$ have been constructed satisfying the
conditions.  Since $f^{r-1}$ is smooth, the fiber product of $f^{r-1}$
with the blowing up $u^{r-1,r}: M^{(r)} \rightarrow M^{(r-1)}$ along
the closed subscheme $M^{(r-1)}_{r-1}$ is canonically isomorphic to
the blowing up of $\AAA^1 \times U^{(r-1)}$ along the preimage of
$M^{(r-1)}_{r-1}$.  By the induction hypothesis and Item $(5)$, the
preimage is precisely $\AAA^1 \times (U^{(r-1)} \cap
M^{(r-1)}_{r-1})$.  So the base-change of $u^{r-1,r}$ by $f^{r-1}$ is
just $\text{Id} \times v^{r-1,r}$.  We define $f^r:\AAA^1 \times
U^{(r)} \rightarrow M^{(r)}$ to be the base-change of $f^{r-1}$ by
$u^{r-1,r}$.

\medskip\noindent
By construction of $f^r$ and the induction hypothesis, it is clear
that $f^1,\dots,f^r$ satisfy Item $(3)$ of the lemma.  Since $f^r$ is
the base-change of the $\mathbb{G}_m$-torsor $f^{r-1}$, also $f^r$ is
a $\mathbb{G}_m$-torsor.  So Item $(4)$ is true.

\medskip\noindent
Next we prove Item $(5)$.  Using the Cartesian property of Item $(3)$
and using that $f^{r-1}$ and $f^r$ are smooth, the process of forming
the strict transform under $u^{r-1,r}$ and then forming the preimage
under $f^r$ is the same as the process of first forming the preimage
under $f^{r-1}$ and then forming the strict transform under $\text{Id}
\times v^{r-1,r}$.  By the induction hypothesis and Item $(5)$, the
preimage under $f^{r-1}$ of $M^{(r-1)}_k$ equals $\AAA^1 \times
(U^{(r-1)} \cap M^{(r-1)}_k)$.  The strict transform of this subscheme
under $\text{Id} \times v^{r-1,r}$ is $\AAA^1 \times (U^{(r)} \cap
M^{(r)}_k)$.

\medskip\noindent
Finally we prove Item $(6)$.  As in the last paragraph, for
$k=0,\dots,r-2$ the pullback by $f^r$ of the strict transform of
$E^{(r-1)}_k$ equals the strict transform of the pullback by
$f^{r-1}$.  By the induction hypothesis and Item $(6)$, the pullback
of $E^{(r-1)}_0$ equals $\{0\}\times U^{(r-1)}$ and the pullback of
$E^{(r-1)}_k$ equals $\AAA^1 \times (U^{(r-1)} \cap E^{(r-1)}_k)$ for
$k=1,\dots,r-2$.  The strict transforms of these subschemes are $\{0\}
\times U^{(r)}$ and $\AAA^1 \times (U^{(r)} \cap E^{(r)}_k)$
respectively.  Finally, using Item $(3)$, the preimage under $f^r$ of
the exceptional divisor of $u^{r-1,r}$, i.e. $E^{(r)}_{r-1}$, equals
the exceptional divisor of $\text{Id} \times v^{r-1,r}$, i.e. $\AAA^1
\times (U^{(r)} \cap E^{(r)}_{r-1})$.  This proves Item $(6)$ and
finishes the proof of the lemma by induction on $r$.
\end{proof}

\medskip\noindent
\begin{prop} \label{prop-smooth} 
This is the main result of this section.
\begin{enumerate}
\item
For each $r=0,\dots,g$, the scheme $M^{(r)}$ is smooth over $S$.  
\item
For each $r=1,\dots,g$, the closed subscheme $E^{(r)}_{0} \cup \dots \cup
E^{(r)}_{r-1}$ is a simple normal crossings divisor in $M^{(r)}$;
moreover the intersection with every geometric fiber over $S$ is a
simple normal crossings divisor.
\item
For each $r=1,\dots,g$, the scheme $M^{(r)}_r$ is smooth over $S$, and
therefore $M^{(r)}_r \rightarrow M^{(r)}$ is a regular embedding. 
\item
For each $r=0,\dots,g$ there exists a locally free sheaf of rank $g$,
$\mc{G}^{(r)}$ on $M^{(r)}$ and a morphism of sheaves
$\phi^{(r)}:\mc{G}^{(r)} \rightarrow (u^{0,r})^* \mc{F}$ such that
$\mc{G}^{(0)} = \mc{G}$ and $\phi^{(0)}=\phi$, such that $\mc{G}^{(1)}
= (u^{0,1})^* \mc{G}^{(0)}(E^{(1)}_0)$ and $\phi^{(1)}$ is as in Lemma
~\ref{lem-pref}, and which satisfy the following condition: for each
$r=1,\dots,g$ there is a factorization $(u^{r-1,r})^*\mc{G}^{(r-1)}
\xrightarrow{\psi^{(r)}} \mc{G}^{(r)} \xrightarrow{\phi^{(r)}}
(u^{0,r})^*\mc{F}$ such that the cokernel of $\psi^{(r)}$ is the
pushforward from $E^{(r)}_{r-1}$ of a locally free sheaf of rank
$g+1-r$ and such that $\phi^{(r)}$ has rank at least $r$ at all
geometric points, and has rank $g$ at the generic point of
$E^{(r)}_{r-1}$.
\item  
The morphisms $\psi^{(r)}$ and $\phi^{(r)}$ above are unique up to
unique isomorphism (if we think of the sheaves $\mc{G}^{(r)}$ as being
subsheaves of $\mc{G}\otimes_{\OO_{M^{(0)}}} K(M^{(0)})$, then the
morphisms are honestly unique).  Moreover they are equivariant for the
obvious action of the group scheme $\text{GL}(\mc{F}) \times
\text{GL}(\mc{G})$.
\item
For each $0 \leq r < s \leq g$ and each $k\geq s$, the inverse image
ideal sheaf $(u^{r,s})^*\mc{I}^{(r)}_k$ equals
\begin{equation}
\mc{I}^{(s)}_k \cdot
\OO_{M^{(s)}} \lt(-\lt(\sum_{i=r}^{s-1} (k+1-i) E^{(s)}_i \rt) \rt).
\end{equation}
\end{enumerate}
\end{prop}

\medskip\noindent
\begin{proof}
The proposition can be checked Zariski locally over $S$.  And Zariski
locally over $S$, the locally free sheaves $\mc{G}$ and $\mc{F}$ are
free.  Therefore we can reduce to the case that $S = \SP(K)$.  We
prove the result by induction on $g$.  For $g=0$, there is nothing to
prove.  Therefore, by way of induction, we may suppose that $g>0$ and
the result has been proved whenever $\text{rank}(\mc{G}) < g$.

\medskip\noindent
The idea of the induction step is the following.  First of all, for
$M^{(0)}$ Items $(1)$ through $(6)$ are all obvious, so we only need
to consider $M^{(r)}$ with $r=1,\dots,g$.  To begin with we restrict
over $U^{(1)}, \dots, U^{(g)}$ and check Items $(1)$ through $(6)$
when restricted over these open sets.  We may check this after making
a smooth surjective base-change on the sets $U^{(r)}$.
Lemma~\ref{lem-mor} gives us a sequence of such base-changes
$\tau^{r-1}(\zeta)$ for any rank $1$ morphism
$\zeta:(S,\mc{G}',\mc{F}') \rightarrow (S,\mc{G},\mc{F})$.  After
base-change by $\tau^{r-1}(\zeta)$, Items $(1)$ through $(6)$ over
$U^{(r)}$ reduce to Items $(1)$ through $(6)$ over
$M^{(r-1)}(S,\mc{G}',\mc{F}')$.  Since $\text{rank}(\mc{G}') <
\text{rank}(\mc{G})$, these follow from the induction hypothesis.

\medskip\noindent
Next, to establish the proposition over all of $M^{(r)}$, we use the
sequence of morphisms $f^r:\AAA^1 \times U^{(r)} \rightarrow M^{(r)}$
from Lemma~\ref{lem-f}.  These morphisms are smooth surjective, so
Items $(1)$ through $(6)$ may be checked after base-change by $f^r$.
And these reduce to Items $(1)$ through $(6)$ over $U^{(r)}$.  Since
we have already established Items $(1)$ through $(6)$ over $U^{(r)}$,
this finishes the induction step.  The proof follows by induction.

\medskip\noindent
First we establish Items $(1)$ through $(6)$ when restricted over
$U^{(1)},\dots , U^{(g)}$.  Let $T \in U^{(0)}$ be a geometric point
with $\text{rank}(T)=1$.  Now we can find $\mc{G}'$ of rank $g-1$,
$\mc{F}'$ of rank $f-1$, $p:\mc{G} \rightarrow \mc{G}'$, and
$q:\mc{F}' \rightarrow \mc{F}$ so that $\zeta=(p,q,T)$ is a morphism
$(\mc{G}',\mc{F}') \rightarrow (\mc{G},\mc{F})$.  By
Lemma~\ref{lem-mor}, for each $r=0,\dots,g-1$, we have a
quasi-compact, separated, smooth morphism
$\tau^r(\zeta):\text{GL}(\mc{F})\times_S \text{GL}(\mc{G})\times_S
M^{(r)}(S,\mc{G}', \mc{F}') \rightarrow M^{(r+1)}(S,\mc{G},\mc{F})$
whose image is $U^{(r+1)}$.

\medskip\noindent
Since $\tau^{r-1}(\zeta): \text{GL}(\mc{F}) \times_S \text{GL}(\mc{G})
\times_S M^{(r-1)}(S,\mc{G}',\mc{F}') \rightarrow U^{(r)}$ is smooth
and surjective, to prove that $U^{(r)}$ is smooth over $S$ it suffices
to prove that $\text{GL}(\mc{F}) \times_S \text{GL}(\mc{G}) \times_S
M^{(r-1)}(S,\mc{G}',\mc{F}')$ is smooth over $S$.  By the induction
hypothesis, $M^{(r-1)}(S,\mc{G}',\mc{F}')$ is smooth over $S$.  And
$\text{GL}(\mc{F})$ and $\text{GL}(\mc{G})$ are obviously smooth over
$S$.  Therefore the fiber product is smooth over $S$.  This
establishes Item $(1)$ over $U^{(r)}$.

\medskip\noindent
Similarly, to show that $U^{(r)} \cap \lt( E^{(r)}_0 \cup \dots \cup
E^{(r)}_{r-1} \rt)$ is a simple normal crossings divisor, it suffices
to prove that the pullback under $\tau^{r-1}(\zeta)$ is a simple
normal crossings divisor.  By Item $(6)$ of Lemma~\ref{lem-mor}, this
pullback is of the form $\text{GL}(\mc{F}) \times_S \text{GL}(\mc{G})
\times_S E$ where $E = E^{(r-1)}_0 \cup \dots E^{(r-1)}_{r-2}$ (the
``missing'' divisor is due to the fact that the pullback of
$E^{(r)}_0$ is the empty set).  By the induction hypothesis, $E$ is a
simple normal crossings divisor in $M^{(r-1)}(S,\mc{G}',\mc{F}')$.
Therefore $\text{GL}(\mc{F}) \times_S \text{GL}(\mc{G}) \times_S E$ is
a simple normal crossings divisor in $\text{GL}(\mc{F}) \times_S
\text{GL}(\mc{G}) \times_S M^{(r-1)}(S,\mc{G}',\mc{F}')$.  It follows
that $U^{(r)} \cap \lt( E^{(r)}_0 \cup \dots \cup E^{(r)}_{r-1} \rt)$
is a simple normal crossings divisor in $U^{(r)}$.  This establishes
Item $(2)$ over $U^{(r)}$.

\medskip\noindent
To show that $U^{(r)}\cap M^{(r)}_r$ is smooth over $S$, it suffices
to prove that the pullback under $\tau^{r-1}(\zeta)$ is smooth over
$S$.  By Item $(5)$ of Lemma~\ref{lem-mor}, the pullback is
$\text{GL}(\mc{F}) \times_S \text{GL}(\mc{G}) \times_S
M^{(r-1)}_{r-1}(S,\mc{G}', \mc{F}')$.  By the induction hypothesis
$M^{(r-1)}_{r-1}(S,\mc{G}',\mc{F}')$ is smooth over $S$, so also
$\text{GL}(\mc{F}) \times_S \text{GL}(\mc{G}) \times_S
M^{(r-1)}_{r-1}(S,\mc{G}',\mc{F}')$ is smooth over $S$.  It follows
that $U^{(r)} \cap M^{(r)}_r$ is smooth over $S$.  This establishes
Item $(3)$ over $U^{(r)}$.

\medskip\noindent 
Item $(4)$ is quite a bit more involved.  By the induction hypothesis,
the maps $(\psi')^{(r)}:(u^{r-1,r})^*\lt(\mc{G}' \rt)^{(r-1)}
\rightarrow \lt(\mc{G}'\rt)^{(r)}$ and the maps $(\phi')^{(r)}: \lt(
\mc{G}' \rt)^{(r)} \rightarrow (u^{0,r})^* \mc{F}'$ on
$M^{(r)}(S,\mc{G}',\mc{F}')$ are all defined and satisfy the
conditions in Item $(4)$ and Item $(5)$.  First we deal with Item
$(5)$ on $U^{(r)}$.  Observe that if the sequence of maps
$\psi^{(s)},\phi^{(s)}$ exists for $s=0,\dots,r-1$, then there is at
most one pair $\psi^{(r)},\phi^{(r)}$ which satisfies the hypothesis.
This is because the restriction of $(u^{r-1,r})^*\phi^{(r-1)}$ to
$E^{(r)}_{r-1}$ has rank at least $r-1$.  So the kernel has rank at
least $g+1-r$.  If there exists a pair $\psi^{(r)},\phi^{(r)}$, then
we must have that the restriction of $(u^{r-1,r})^*\phi^{(r-1)}$ to
$E^{(r)}_{r-1}$ has constant rank $r-1$ (i.e. the cokernel is locally
free of rank $g+1-r$) and $\psi^{(r)}: (u^{r-1,r})^*\mc{G}^{(r-1)}
\rightarrow \mc{G}^{(r)}$ must be the elementary transform up along
$E^{(r)}_{r-1}$ whose kernel equals the kernel of
$(u^{r-1,r})^*\phi^{(r-1)}$.  And then $\phi^{(r)}$ is the unique
morphism through which $(u^{r-1,r})^*\phi^{(r-1)}$ factors.  This
establishes Item $(4)$ (in fact without restricting over $U^{(r)}$).
Equivariance with respect to $\text{GL}(\mc{F}) \times_S
\text{GL}(\mc{G})$ follows by induction on $r$ and the uniqueness just
mentioned.

\medskip\noindent
Now we show the existence of $\psi^{(r)},\phi^{(r)}$ when restricted
over $U^{(r)}$.  We will prove this by faithfully flat (in fact
smooth) descent, i.e. we will construct a descent datum for the
faithfully flat cover $\tau^{r-1}(\zeta): \text{GL}(\mc{F}) \times_S
\text{GL}(\mc{G}) \times_S M^{(r-1)}(S,\mc{G}',\mc{F}') \rightarrow
U^{(r)}$.  The uniqueness in Item $(4)$ and equivariance with respect
to $\text{GL}(\mc{F}) \times_S \text{GL}(\mc{G})$ is what will give
the cocycle condition.

\medskip\noindent
For each $r=1,\dots,g$ define $\mc{G}^{(r)}_{\text{pre}}$ on
$M^{(r-1)}(S,\mc{G}',\mc{F}')$ to be the direct sum of $\lt( \mc{G}'
\rt)^{(r-1)}$ and $\text{Ker}(p) \otimes_{\OO_S} \OO_{M^{(r-1)}}$.  In
particular, $\mc{G}^{(1)}$ is simply $\mc{G} \otimes_{\OO_S}
\OO_{M^{(0)}}$.  For each $r=2,\dots, g$ define
$\psi^{(r)}_{\text{pre}}: (u^{r-2,r-1})^*\mc{G}^{(r-1)}_{\text{pre}}
\rightarrow \mc{G}^{(r)}_{\text{pre}}$ to be the direct sum of
$(\psi')^{(r-1)}:(u^{r-2,r-1})^* \lt( \mc{G}' \rt)^{(r-2)} \rightarrow
\lt( \mc{G}' \rt)^{(r-1)}$ and the identity map on $\text{Ker}(p)
\otimes_{\OO_S} \OO_{M^{(r-1)}}$.  For each $r=1,\dots, g$ define
$\phi^{(r)}_{\text{pre}}: \mc{G}^{(r)}_{\text{pre}} \rightarrow
\mc{F}\otimes_{\OO_S} \OO_{M^{(r-1)}}$ to be the sum of the map
\begin{equation}
q\circ
(\phi')^{(r-1)}: \lt( \mc{G}' \rt)^{(r-1)} \rightarrow \mc{F}'
\otimes_{\OO_S} \OO_{M^{(r-1)}} \rightarrow \mc{F} \otimes_{\OO_S}
\OO_{M^{(r-1)}}
\end{equation}
with the map $T: \text{Ker}(p)\otimes_{\OO_S} \OO_{M^{(r-1)}}
\rightarrow \mc{F} \otimes_{\OO_S} \OO_{M^{(r-1)}}$.  

\medskip\noindent
On $\text{GL}(\mc{F})$ there is a universal automorphism
$\alpha:\mc{F} \otimes_{\OO_S} \OO \rightarrow \mc{F} \otimes_{\OO_S}
\OO$ and on $\text{GL}(\mc{F})$ there is a universal automorphism
$\beta:\mc{G} \otimes_{\OO_S} \OO \rightarrow \mc{G} \otimes_{\OO_S}
\OO$.  Also denote by $\alpha$ and $\beta$ the pullbacks of these
automorphisms to $\text{GL}(\mc{F}) \times_S \text{GL}(\mc{G})
\times_S M^{(r-1)}(S,\mc{G}', \mc{F}')$.  By definition, the pullback
by $\tau^{r-1}(\zeta)$ of $(u^{0,r})^*\phi$ equals $\alpha \circ
(u^{0,r-1})^*\text{pr}_3^*\phi^{(1)}_{\text{pre}} \circ \beta^{-1}$.
Now we specify the part of the descent data on $\text{GL}(\mc{F})
\times_S \text{GL}(\mc{G}) \times_S M^{(r-1)}(S,\mc{G}',\mc{F}')$
which define $\mc{G}^{(r)}$, $\psi^{(r)}$ and $\phi^{(r)}$.  For each
$r=1,\dots, g$ we define $\tau^{r}(\zeta)^* \mc{G}^{(r)}$ to be
$\text{pr}_3^* \mc{G}^{(r)}_{\text{pre}}$.  We define
$\tau^{0}(\zeta)^* \psi^{(1)}$ to be $\beta^{-1}$, where the domain of
$\beta^{-1}$ is identified with $\mc{G}\otimes_{\OO_S} \OO \cong
\tau^{0}(\zeta)^* (u^{0,1})^* \mc{G}^{(0)}$ and where the range of
$\beta^{-1}$ is identified with $\mc{G}\otimes_{\OO_S} \OO \cong
\text{pr}_3^* \mc{G}^{(1)}_{\text{pre}}$.  For $r=2,\dots,g$, we
define $\tau^{r-1}(\zeta)^* \psi^{(r)}$ to be $\text{pr}_3^*
\psi^{(r)}_{\text{pre}}$.  For all $r=1,\dots,g$, we define
$\tau^{r-1}(\zeta)^* \phi^{(r)}$ to be $\alpha \circ \text{pr}_3^*
\phi^{(r)}_{\text{pre}}$.

\medskip\noindent
Now to finish specifying the descent data, we have to give patching
morphisms on the fiber product of $\tau^{r-1}(\zeta)$ with itself.
There are canonical descent data associated to the sheaves $\mc{G}
\otimes_{\OO_S} \OO_{M^{(r)}}$ and $\mc{F} \otimes_{\OO_S}
\OO_{M^{(r)}}$ on $M^{(r)}(S,\mc{G},\mc{F})$.  And, up to unique
isomorphism, there is at most one way of extending the descent data
for $\mc{G}^{(r)}$, $\psi^{(r)}$ and $\phi^{(r)}$ so that the descent
data giving $\psi^{(r)}$ and $\phi^{(r)}$ are morphisms from the
descent datum for $\mc{G} \otimes_{\OO_S} \OO_{M^{(r)}}$ to the
descent datum for $\mc{G}^{(r)}$ and from the descent datum for
$\mc{G}^{(r)}$ to the descent datum for $\mc{F} \otimes_{\OO_S}
\OO_{M^{(r)}}$ respectively.  This doesn't prove that such a descent
datum exists.  Proving existence is an exercise in the compatibilities
of all the sheaves and morphisms defined so far, and is left to the
interested reader.  The key point, as always, is that on
$M^{(r)}(S,\mc{G},\mc{F})$, on the base-change by $\tau^{r-1}(\zeta)$,
and on the double base-change by $\tau^{r-1}(\zeta)$, the morphism
$\psi^{(r)}$ is, up to unique isomorphism, the elementary transform up
determined by the kernel of $(u^{r-1,r})^*\phi^{(r-1)}$ restricted to
$E^{(r)}_{r-1}$.

\medskip\noindent
The upshot is that the sheaves $\mc{G}^{(r)}$ and sheaf maps
$\psi^{(r)}$, $\phi^{(r)}$ exist when restricted over $U^{(r)}$.  To
check the properties in Item $(4)$, i.e. that the cokernel of
$\psi^{(r)}$ is as specified and that the rank of $\phi^{(r)}$ is as
specified, we can check after base-change by $\tau^{r-1}(\zeta)$.  And
then it follows from the construction of our descent datum, and by the
induction hypothesis applied to $M^{(r-1)}(S,\mc{G}',\mc{F}')$.
Again, the details are left to the interested reader.

\medskip\noindent
Next we show that Item $(6)$ holds when we restrict over $U^{(s)}$.
The map $v^{1,0}:U^{(1)} \rightarrow U^{(0)}$ is an isomorphism, so we
may suppose that $1 \leq r < s \leq g$.  To check that two ideal
sheaves are equal, we can check after faithfully flat base-change, so
we base-change by $\tau^{s-1}(\zeta)$.  By Item $(3)$ of
Lemma~\ref{lem-mor}, the inverse image under $\tau^{s-1}(\zeta)$ of
the inverse image under $u^{r,s}$ equals the inverse image under
$\text{Id} \times \text{Id} \times u^{r-1,s-1}$ of
$\tau^{r-1}(\zeta)$.  By Item $(5)$ of Lemma~\ref{lem-mor}, the
inverse image under $\tau^{r-1}(\zeta)$ of $\mc{I}^{(r)}_k$ equals the
ideal sheaf of the closed subscheme $\text{GL}(\mc{F}) \times_S
\text{GL}(\mc{G}) \times_S M^{(r-1)}_{k-1}(S,\mc{G}',\mc{F}')$ of
$\text{GL}(\mc{F}) \times_S \text{GL}(\mc{G}) \times_S
M^{(r-1)}(S,\mc{G}',\mc{F}')$.  By the induction hypothesis and Item
$(5)$, the inverse image of $\mc{I}^{(r-1)}_{k-1}(S,\mc{G}',\mc{F}')$
under $u^{r-1,s-1}$ is the product of
$\mc{I}^{(s-1)}_{k-1}(S,\mc{G}',\mc{F}')$ with the invertible ideal
sheaf associated to the Cartier divisor $\sum_{j=r-1}^{s-2}
((k-1)+1-j) E^{(s-1)}_j$.  Making the substitution $i=j+1$, the
Cartier divisor is $\sum_{i=r}^{s-1}(k+1-i) E^{(s-1)}_{i-1}$.  Taking
the inverse image of this ideal sheaf under $\text{pr}_3$, and using
Item $(5)$ and Item $(6)$ of Lemma~\ref{lem-mor}, we get the same
ideal sheaf as the inverse image under $\tau^{s-1}(\zeta)$ of the
ideal sheaf in Item $(6)$ above.  This establishes Item $(6)$ over
$U^{(r)}$.  So the proposition is proved ``over $U^{(r)}$''.

\medskip\noindent
To finish the induction step, we have to prove that Items $(1)$
through $(6)$ hold over all of $M^{(r)}(S,\mc{G},\mc{F})$.  To do this
we use the morphisms $f^r:\AAA^1 \times U^{(r)} \rightarrow M^{(r)}$
from Lemma~\ref{lem-f}.  By Item $(4)$ of Lemma~\ref{lem-f}, the
morphism $f^r$ is smooth and surjective, so the check the target of
$f^r$ is smooth over $S$, it suffices to check the domain of $f^r$ is
smooth over $S$.  As established above, $U^{(r)}$ is smooth over $S$
so that $\AAA^1 \times U^{(r)}$ is smooth over $S$.  Therefore
$M^{(r)}$ is smooth over $S$.  This establishes Item $(1)$ over
$M^{(r)}$.

\medskip\noindent
Similarly, to show that $E^{(r)}_0 \cup \dots \cup E^{(r)}_{r-1}$ is a
simple normal crossings divisor, it suffices to prove that the
pullback by $f^{r}$ is a simple normal crossings divisor.  By Item
$(6)$ of Lemma~\ref{lem-f}, the preimage of $E^{(r)}_0$ is $\{0\}
\times U^{(r)}$ and the preimage of $E^{(r)}_1 \cup \dots \cup
E^{(r)}_{r-1}$ is $\AAA^1 \times \lt(U^{(r)} \cap \lt( E^{(r)}_1 \cup
\dots \cup E^{(r)}_{r-1} \rt) \rt)$.  As established above, $U^{(r)}
\cap \lt( E^{(r)}_1 \cup \dots \cup E^{(r)}_{r-1} \rt)$ is a simple
normal crossings divisor in $U^{(r)}$.  So by Lemma~\ref{lem-snc}, the
divisor $(f^r)^{-1} \lt( E^{(r)}_0 \cup \dots \cup E^{(r)}_{r-1} \rt)$
is a simple normal crossings divisor.  This establishes Item $(2)$
over $M^{(r)}$.

\medskip\noindent
To show that $M^{(r)}_r$ is smooth over $S$, it suffices to show that
the preimage under $f^r$ is smooth over $S$.  By Item $(5)$ of
Lemma~\ref{lem-f}, the preimage of $M^{(r)}_r$ is $\AAA^1 \times \lt(
U^{(r)} \cap M^{(r)}_r \rt)$.  As established above, $U^{(r)} \cap
M^{(r)}_r$ is smooth.  So $\AAA^1 \times \lt( U^{(r)} \cap M^{(r)}_r
\rt)$ is smooth, and therefore $M^{(r)}_r$ is smooth.  This
establishes Item $(3)$ over $M^{(r)}$.

\medskip\noindent
As over $U^{(r)}$, Item $(4)$ and Item $(5)$ are a bit more involved
(although we've already done most of the work).  As before, Item $(5)$
is automatic provided we can prove the existence of $\mc{G}^{(r)}$,
$\psi^{(r)}$ and $\phi^{(r)}$ satisfying the hypotheses of Item $(4)$.
We prove existence by faithfully flat descent with respect to the
faithfully flat (in fact smooth) morphism $f^r:\AAA^1 \times U^{(r)}
\rightarrow M^{(r)}$.  Now we specify the part of the descent data on
$\AAA^1 \times U^{(r)}$ which define $\mc{G}^{(r)}$, $\psi^{(r)}$ and
$\phi^{(r)}$.  For each $r=1,\dots, g$, we define $(f^r)^*
\mc{G}^{(r)}$ to be $\text{pr}_2^* (\mc{G}^{(r)}|_{U^{(r)}})\otimes
\OO_{\AAA^1 \times U^{(r)}}(\{ 0 \} \times U^{(r)})$ (of course this
definition looks circular, but recall that we have constructed
$\mc{G}^{(r)}|_{U^{(r)}}$, $\psi^{(r)}|_{U^{(r)}}$ and
$\phi^{(r)}|_{U^{(r)}}$ above).  Let $t$ be the coordinate on $\AAA^1$
considered as a global section of the invertible sheaf $\OO_{\AAA^1
\times U^{(r)}}(\{ 0 \} \times U^{(r)})$ whose vanishing locus is
precisely $\{ 0 \} \times U^{(r)}$; this is a bit at odds with the
usual terminology which would call this section $1$.  The point is
that there is a canonical everywhere nonzero global section of this
invertible sheaf (which in the usual terminology would be denoted by
$\frac{1}{t}$, but which we prefer to denote by $1$), and with respect
to this trivialization the regular function $t$ corresponds to a
section whose vanishing locus is $\{ 0 \} \times U^{(r)}$.

\medskip\noindent 
We define $(f^1)^*(\psi^{(1)})$ to be the map
\begin{equation}
\text{Id} \otimes t:
\mc{G} \otimes_{\OO_S} \OO_{\AAA^1 \times U^{(1)}} \rightarrow \mc{G}
\otimes_{\OO_S} \OO_{\AAA^1 \times U^{(1)}} ( \{0\} \times U^{(1)} ).
\end{equation}
Of course we identify the domain of this map with the pullback by
$f^1$ of $\mc{G}^{(0)} = \mc{G} \otimes_{\OO_S} \OO_{M^{(1)}}$ and we
identify the target with $(f^1)^* \mc{G}^{(1)}$ defined above.  For
$r=2,\dots, g$, we define $(f^r)^*(\psi^{(r)})$ to be the map
$\text{pr}_2^*(\psi^{(r)}|_{U^{(r)}})\otimes \text{Id}$.  We define
$(f^1)^*(\phi^{(1)})$ to be the composition of the canonical
isomorphism $\text{Id} \otimes 1: \text{pr}_2^*\mc{G}^{(r)}\otimes
\OO_{\AAA^1 \times U^{(r)}}(\{ 0 \} \times U^{(r)}) \rightarrow
\text{pr}_2^* \mc{G}^{(r)}$ with
$\text{pr}_2^*(\phi^{(r)}|_{U^{(r)}})$.  Observe that $(f^1)^*
\phi^{(1)}$ is the same map constructed in Lemma~\ref{lem-pref}.

\medskip\noindent
We have to check that these definitions of $(f^r)^*\psi^{(r)}$ and
$(f^r)^*\phi^{(r)}$ have the properties from Item $(4)$.  For $r=1$,
this is precisely Lemma~\ref{lem-pref}.  Suppose that $r \geq 2$.  As
established above, $(\phi^{(r)}|_{U^{(r)}}) \circ
(\psi^{(r)}|_{U^{(r)}})$ equals
$(v^{r-1,r})^*(\phi^{(r-1)}|_{U^{(r-1)}})$.  Pulling back by
$\text{pr}_2^*$ and using Item $(3)$ from Lemma~\ref{lem-f}, we
conclude that $(f^r)^*\phi^{(r)}\circ (f^r)^* \psi^{(r)}$ equals
$(\text{Id} \times v^{r-1,r})^*(f^{r-1})^*\phi^{(r-1)} = (f^r)^*
(u^{r-1,r})^* \phi^{(r-1)}$.  This is the first necessary property.
As established above, the cokernel of $\psi^{(r)}|_{U^{(r)}}$ is the
push forward of a locally free sheaf of rank $g+1-r$ from the divisor
$U^{(r)}\cap E^{(r)}_{r-1}$ in $U^{(r)}$.  Therefore the cokernel of
$(f^r)^* \psi^{(r)}$, i.e. the cokernel of $\text{pr}_2^*
(\psi^{(r)}|_{U^{(r)}})$, is the push forward of a locally free sheaf
of rank $g+1-r$ from the divisor $\AAA^1 \times (U^{(r)} \cap
E^{(r)}_{r-1})$, i.e. the divisor $(f^r)^{-1}(E^{(r)}_{r-1})$.  This
is the second necessary property.  And as established above,
$\phi^{(r)}|_{U^{(r)}}$ has rank at least $r$ at all geometric points.
Therefore $(f^r)^*\phi^{(r)}$,
i.e. $\text{pr}_2^*(\phi^{(r)}|_{U^{(r)}})$, has rank at least $r$ at
all geometric points.  This is the last necessary property.

\medskip\noindent
To finish the proof of Item $(4)$, we need to define the part of the
descent data coming from the double base-change by $f^r$.  As above,
there is at most one way of completing the descent data so that
$\psi^{(r)}$ and $\phi^{(r)}$ give morphisms with the given descent
data for $\mc{G} \otimes_{\OO_S} \OO_{M^{(r)}}$ and $\mc{F}
\otimes_{\OO_S} \OO_{M^{(r)}}$.  Checking that one can complete the
descent data is an exercise left to the interested reader, where the
key point, as above, is that $\psi^{(r)}$ is the elementary transform
up determined by the kernel of $(u^{r-1,r})^*\phi^{(r-1)}$ restricted
to $E^{(r)}_{r-1}$.  This finishes the proof that Item $(4)$ and Item
$(5)$ are satisfied on $M^{(r)}$.

\medskip\noindent
Finally we establish Item $(6)$.  First suppose that $r=0$.  By Item
$(4)$ of Lemma~\ref{lem-pref}, $(u^{0,1})^{-1}(\mc{I}^{(0)}_k)$ equals
$\mc{I}^{(1)}_k \cdot \OO_{M^{(1)}}(-(k+1)E^{(1)}_0)$.  By Item $(5)$
of Lemma~\ref{lem-f}, $(f^1)^{-1}(\mc{I}^{(1)}_k)$ is the ideal sheaf
$\text{pr}_2^{-1}(\mc{I}^{(1)}_k|_{U^{(1)}})$.  By Item $(6)$ of
Lemma~\ref{lem-f}, $(f^1)^{-1}\OO_{M^{(1)}}(-E^{(1)}_0)$ is the ideal
sheaf of $\{0\} \times U^{(1)}$.  The process of forming the preimage
of an ideal sheaf by $\text{pr}_2$ and then forming the preimage of
that ideal sheaf by $\text{Id} \times v^{1,s}$ is the same as the
process of first forming the preimage of the ideal sheaf by $v^{1,s}$
and then forming the preimage by $\text{pr}_2$.  As established above,
$(v^{1,s})^{-1}(\mc{I}^{(1)}_k|_{U^{(1)}})$ equals the restriction to
$U^{(s)}$ of the ideal sheaf
\begin{equation}
\mc{I}^{(s)}_k \cdot \OO_{U^{(s)}}\lt(- \lt( \sum_{i=1}^{s-1} (k+1-i)
E^{(s)}_i \rt) \rt).
\end{equation}
And the inverse image under $\text{Id} \times v^{1,s}$ of the ideal sheaf of
$\{0\} \times U^{(1)}$ is the ideal sheaf of $\{ 0 \} \times
U^{(s)}$.  Putting the pieces together, the inverse image under $\text{Id}
\times v^{1,s}$ of the inverse image under $f^1$ of the inverse image
under $u^{0,1}$ of $\mc{I}^{(0)}_k$ equals the ideal sheaf
\begin{equation}
\text{pr}_2^{-1}(\mc{I}^{(s)}_k|_{U^{(s)}}) \cdot \OO_{\AAA^1 \times
  U^{(s)}} \lt( -(k+1)\{0\} \times U^{(s)} \rt) \otimes
\OO_{\AAA^1 \times U^{(s)}}
  \lt( - \lt( \sum_{i=1}^{s-1}(k+1-i) \AAA^1 \times \lt(
  U^{(s)} \cap E^{(s)}_i \rt) \rt) \rt).
\end{equation}
By Item $(5)$ and Item $(6)$ of Lemma~\ref{lem-f}, this is precisely
the inverse image under $f^s$ of the ideal sheaf which appears in Item
$(6)$ of the proposition.  Using Item $(3)$ of Lemma~\ref{lem-f} one
last time, and using that one can check equality of ideal sheaves
after faithfully flat base-change, we conclude that Item $(6)$ holds
when $r=1$.

\medskip\noindent
Checking Item $(6)$ when $r > 1$ is even easier and follows by the
same sort of argument as above; the details are left to the interested
reader.  This finishes the proof that Item $(6)$ holds for $M^{(r)}$,
and thus finishes the proof that the proposition holds for $M^{(r)}$.
So the proposition is proved by induction on the rank of $\mc{G}$.
\end{proof}

\medskip\noindent
\subsection{Computation of the log discrepancies} \label{subsec-comp}
Using the log resolution from the last section, we compute the log
discrepancies of the pairs $(M^{(0)},M^{(0)}_k)$.  We also have reason
to compute the log discrepancies of the pairs $(M^{(0)},q\cdot
M^{(0)}_k)$.  Then, combined with the results of
Section~\ref{sec-discrep}, we use these computations to find the log
discrepancies of some projective cones.

\medskip\noindent
\begin{lem} \label{lem-mults} 
Let $(S,\mc{G},\mc{F})$ be a datum with $\text{rank}(\mc{G}) = g$ and
$\text{rank}(\mc{F}) = f$.  For each $0 \leq r < s \leq g$, the
relative canonical divisor of $u^{r,s}:M^{(s)} \rightarrow M^{(r)}$
equals
\begin{equation}
K_{M^{(s)}} - (u^{r,s})^*K_{M^{(r)}} = \sum_{i=r}^{s-1} \lt(
(f-i)(g-i) -1 \rt) E^{(s)}_{i}.
\end{equation}
For $r< s \leq k < g$ and each positive integer $q$, 
the inverse image under $u^{r,s}$ of the ideal sheaf
$\lt( \mc{I}^{(r)}_k \rt)^q$ equals
\begin{equation}
\lt( \mc{I}^{(s)}_k \rt)^q \cdot \OO_{M^{(s)}} \lt( - q\lt(
\sum_{i=r}^{s-1} (k+1-i) E^{(s)}_i \rt) \rt).
\end{equation}
and the associated cycle is
\begin{equation}
q \lt[ M^{(s)}_{k} \rt] + \sum_{i=r}^{s-1} q(k+1-i)E^{(s)}_{i}
\end{equation}
Finally, for $k < s\leq g$ and each positive integer $q$, the inverse
image under $u^{r,s}$ of the ideal sheaf $\lt( \mc{I}^{(r)}_k \rt)^q$
is an invertible ideal sheaf defining the Cartier divisor
\begin{equation}
\sum_{i=r}^{k} q(k+1-i) E^{(s)}_{i}
\end{equation}
\end{lem}

\medskip\noindent
\begin{proof}
By Proposition~\ref{prop-smooth}, we know that $M^{(r)}_r \subset
M^{(r)}$ is a regular embedding.  And by ~\cite[Prop., p. 67]{ACGH},
the codimension of $M^{(r)}_r$, which equals the codimension of
$M^{(0)}_r$ in $M^{(0)}$, equals $(f-r)(g-r)$.  Therefore the relative
canonical divisor of the blowing up $u^{r,r+1}:M^{(r+1)} \rightarrow
M^{(r)}$ is $(f-r)(g-r)E^{(r+1)}_r$.  The first formula follows since
the relative canonical divisor of a composition of birational
morphisms is the sum of the relative canonical divisors of the
separate morphisms.

\medskip\noindent
The second formula follows from Item $(6)$ of
Proposition~\ref{prop-smooth}.  The final formula follows from the
second formula and that fact that
$(u^{r,r-1})^{-1}\mc{I}^{(r-1)}_{r-1}$ equals the invertible ideal
sheaf $\OO_{M^{(r)}}(-E^{(r)}_{r-1})$.
\end{proof}

\medskip\noindent
\begin{cor} \label{cor-mults} 
Let $(S,\mc{G},\mc{F})$ be a datum with $\text{rank}(\mc{G}) = g$ and
$\text{rank}(\mc{F}) = f$.  Suppose that $S$ is smooth, that $g\geq
1$, that $f \geq 2$ and that $0\leq k \leq g-1$.  Consider the pair
$(M^{(r)}(S,\mc{G},\mc{F}), q\cdot M^{(r)}_{k}(S,\mc{G},\mc{F}))$
where $q \geq 0$.  For $i=r,\dots,k$, the log discrepancies are
\begin{equation}
a(E^{(g)}_i;M^{(r)},q\cdot M^{(r)}_{k}) = (f-i)(g-i) - q(k+1-i).
\end{equation}
Define $a = \text{min}\{ (f-i)(g-i) - q(k+1-i)| i=r,\dots,k \}$.  Then
$(M^{(r)}, q\cdot M^{(r)}_{k})$ is log canonical iff $a \geq 0$, in
which case the minimal log discrepancy $\text{mld}(M^{(r)};M^{(r)},
q\cdot M^{(r)}_{k})$ equals $\text{min}(1,a)$.  
In particular, if $q\leq f-g+1$,
then the pair $(M^{(r)},q\cdot M^{(r)}_{g-1})$ is log canonical and
the minimal log discrepancy equals $\text{min}(1,f-g+1-q)$.
\end{cor}

\medskip\noindent
\begin{proof}
The corollary follows from the computations in Lemma~\ref{lem-mults}
and the definition of the log discrepancies and minimal log
discrepancies of a pair (c.f. ~\cite[Defn. 1.1]{EMY}).
\end{proof}

\medskip\noindent
Let $(S,\mc{G},\mc{F})$ be a datum with $\text{rank}(\mc{G}) = g$ and
$\text{rank}(\mc{F}) = f$.  Let $\mc{A}$ and $\mc{A}'$ be locally free
sheaves on $S$ with $\text{rank}(\mc{A}) = a$ and
$\text{rank}(\mc{A}') =a'$ with $a > 0$.  Define $\mc{E}$ to be the
cokernel of the following sheaf map on $M^{(0)}$
\begin{equation}
\text{Id}\otimes \phi: \mc{A} \otimes_{\OO_S} \mc{G} \otimes_{\OO_S}
\OO_{M^{(0)}} \rightarrow \mc{A} \otimes_{\OO_S} \mc{F}
\otimes_{\OO_S} \OO_{M^{(0)}} \oplus \mc{A}'\otimes_{\OO_S} \OO_{M^{(0)}}.
\end{equation}
The sheaf map is zero on the summand $\mc{A}' \otimes_{\OO_S}
\OO_{M^{(0)}}$, so this locally free sheaf will be a direct summand of
$\mc{E}$. 
Denote by the pair $(\pi:C\rightarrow M^{(0)}, \alpha:\pi^* \mc{E}
\rightarrow \mc{Q} )$ the relative Grassmannian cone parametrizing
rank $r$ locally free quotients of $\mc{E}$.

\medskip\noindent
\begin{prop} \label{prop-cone} 
If $a\cdot r \leq f-g$ and $S$ is smooth and geometrically connected,
then $C$ has pure dimension equal to the expected dimension $d =
\text{dim}(C) = \text{dim}(S) + f\cdot g + r\lt( (a(f-g)+a') - r\rt)$
and $C$ is a normal, integral, local complete intersection scheme
which is canonical.
\end{prop}

\medskip\noindent
\begin{proof}
Denote $g' = a\cdot g = \text{rank}(\mc{A} \otimes_{\OO_S} \mc{G})$.
For each $k=0,\dots, g-1$ and $l=0,\dots,a-1$, we have that $B_{a\cdot
k + l} = B_{a \cdot k} = M^{(0)}_k$ which has codimension
$(f-k)(g-k)=(f-g)(g-k) + (g-k)^2$.  And $r(g'-(a\cdot k + l))+1 =
a\cdot r(g-k)+1 - rl$.  By assumption, $f-g \geq a\cdot r$, and for $k
\leq g-1$ we have $(g-k)^2 \geq 1$.  Therefore $(f-g)(g-k) + (g-k)^2
\geq a\cdot r(g-k) +1 -rl$ for all $k=0,\dots, g-1$ and
$l=0,\dots,a-1$.  Therefore, by Lemma~\ref{lem-codims} we have that
$C$ is irreducible of the expected dimension.  And by the proof of
Lemma~\ref{lem-discrep}, $C$ is a local complete intersection scheme.
In particular it is Cohen-Macaulay.  Moreover, $\pi:C \rightarrow
M^{(0)}$ is smooth over $M^{(0)}-M^{(0)}_{g-1}$.  So $C$ is
generically reduced.  Since $C$ is Cohen-Macaulay, it follows that $C$
is everywhere reduced.  A reduced Cohen-Macauly scheme satisfies
Serre's condition $S2$ for normality.  Thus to prove that $C$ is
normal, it suffices to prove that $C$ is regular in codimension $1$.

\medskip\noindent
If $a\cdot r < f-g$, then the same parameter count as above shows that
for all $k=0,\dots, g-1$ and $l=0,\dots, a-1$ we have that $(f-k)(g-k)
\geq r(g'-(a\cdot k + l)) + 2$.  By Lemma~\ref{lem-codims} it follows
that $C$ is regular in codimension $1$ so that $C$ is normal.
Therefore assume that $a\cdot r = f-g$.  Of course $C$ is regular on
the dense open subset $\pi^{-1}(M^{(0)} - M^{(0)}_{g-1})$.  The only
codimension one point not contained in this locus is the generic point
of $\pi^{-1}(M^{(0)}_{g-1} - M^{(0)}_{g-2})$.

\medskip\noindent
As in the proof of Lemma~\ref{lem-codims}, denote by $(\rho:C'
\rightarrow M^{(0)}, \beta: \rho^*\mc{F} \rightarrow \mc{Q}')$ the
Grassmannian bundle parametrizing rank $r$ locally free quotients of
$\mc{A} \otimes_{\OO_S} \mc{F} \otimes_{\OO_S} \OO_{M^{(0)}} \oplus
\mc{A}'\otimes_{\OO_S} \OO_{M^{(0)}}$.  There is a natural closed
immersion $h:C \rightarrow C'$ compatible with projection to
$M^{(0)}$.

\medskip\noindent
Observe that $C' = P \times_S M^{(0)}$ where $\sigma:P \rightarrow S$
is the Grassmannian bundle parametrizing rank $r$ locally free
quotients of $\mc{A} \otimes_{\OO_S} \mc{F} \oplus \mc{A}'$.
Therefore we have a projection $\text{pr}_2:C \rightarrow P$
compatible with projection to $S$.  The question is local, so we may
base-change to an open subset of $S$ over which $\mc{A}$ is trivial.
Choose an ordered basis for $\mc{A}$ so that $\mc{A} \otimes_{\OO_S}
\mc{F}$ is just $\mc{F}^{\oplus a}$.  Let $W \subset P$ denote the
dense open set over which the sheaf map $(\pi')^*(\mc{A}
\otimes_{\OO_S} \mc{F}) \rightarrow \mc{Q}'$ is surjective.  On $W$
there is a smooth, surjective morphism to the Grassmannian $P'$
parametrizing rank $r$ locally free quotients of $\mc{A}
\otimes_{\OO_S} \mc{F}$.  Any rank $r$ quotient space of $\mc{A}
\otimes_{\OO_S} \mc{F}$, i.e. of $\mc{F}^{\oplus a}$ is represented as
the image of a matrix $\mc{F}^{\oplus a} \rightarrow \OO^{\oplus r}$
of the form
\begin{equation}
M = \lt[ \begin{array}{cccc}
v_{1,1} & v_{1,2} & \dots & v_{1,a} \\
\vdots & \vdots  & \ddots & \vdots \\
v_{r,1} & v_{r,2} & \dots & v_{r,a}
\end{array} \rt]
\end{equation}
where the $v_{i,j}$ are sections of $\mc{F}^\vee$.  This matrix
determines a point $x'$ in $P'$.  Let $x\in W$ be any point mapping to
$x'$.  If we think of the fiber of $\text{pr}_2:C \rightarrow P'$ over
$x$ as a subscheme of $M^{(0)}$, then it is just the subscheme of
matrices $L :\mc{G} \rightarrow \mc{F}$ such that $M\circ L^{\oplus
a}$ is the zero matrix, i.e. it is the set of matrices $L$ such that
the kernel of the transpose matrix $L^\dagger$ contains the subspace
\begin{equation}
K = \text{span}\{v_{i,j} | 1\leq i \leq r, 1\leq j \leq a \}.
\end{equation}
This is the space of matrices $L^\dagger$ from $\mc{F}^\vee / K$ to
$\mc{G}^\vee$.  Let $V' \subset P'$ be the dense open subset
parametrizing quotients where $\text{dim}(K) = a\cdot r$ and let $V
\subset W$ be the preimage of $V$.  Over $V$ we conclude that
$\text{pr}_2:C \rightarrow P$ is a vector bundle of rank $(f-a\cdot
r)g$.  In particular, the preimage is a nonempty, smooth scheme,
i.e. $\text{pr}_2^{-1}(V)$ is contained in the smooth locus
$C_{\text{smooth}}$.  But of course, the map $L^\dagger$ may still
have any rank between $0$ and $g$ (recall that $f-a\cdot r \geq g$, so
that the dimension of $\mc{F}^\vee /K$ is greater than the dimension
of $\mc{G}^\vee$.  Therefore this open set intersects the preimage of
every strata $M^{(0)}_k - M^{(0)}_{k-1}$.  Combined with an obvious
homogeneity argument, the open set intersects every fiber of $\pi$.
So the smooth locus $C_{\text{smooth}}$ intersects every fiber of
$\pi$ and we conclude that $C$ is regular in codimension $1$ points.
So, in every case, $C$ is normal.

\medskip\noindent
The $(a\cdot (f-g)+a')^{\text{th}}$ Fitting ideal $\mc{J}$ of
$\text{Id} \otimes \phi$, i.e. the ideal generated by the maximal
minors of the matrix of $\text{Id}\otimes \phi$, is easily seen to be
$\mc{I}^a$ where $\mc{I}$ is the $(f-g)^{\text{th}}$ Fitting ideal of
$\mc{I}$, i.e. the ideal sheaf of $M^{(0)}_{g-1}$.  By
Proposition~\ref{prop-relcan}, the cone $(C,\emptyset)$ is canonical
iff the pair $(M^{(0)}, a\cdot r \cdot M^{(0)}_{g-1} )$ is canonical.
And by Corollary~\ref{cor-mults}, $(M^{(0)},a\cdot r \cdot
M^{(0)}_{g-1})$ is canonical.  Therefore $(C,\emptyset)$ is canonical,
which finishes the proof.
\end{proof}

\begin{rmk} \label{rmk-cone} 
When $a\cdot r > f-g$, the cone $C$ has more than one irreducible
component.  It would be interesting to determine the minimal log
discrepancies of the different irreducible components of $C$, in
particular of the unique irreducible component which dominates
$M^{(0)}$.  The first case is when $a\cdot r=f-g+1$.  In this case it
follows from the proof above that the second irreducible component is
the closure of the preimage of $M^{(0)}_{g-1} - M^{(0)}_{g-2}$, and
that the restriction of $\text{pr}_2$ to this irreducible component is
birationally a vector bundle of rank $(f-a\cdot r)g$ over $P$.  It may
be possible to use this structure to compute the minimal log
discrepancies of the two irreducible components.
\end{rmk}

\section{Adjunction for $(B,B_{g-1})$} \label{sec-spec} 
\medskip\noindent
Let $K$ be a field, not necessarily algebraically closed or of
characteristic zero.  In this section we will work in the category of
$K$-schemes.  We use the results of the last section.  The interested
reader will see how to prove all the analogous results over an
arbitrary base scheme.

\medskip\noindent
In Section~\ref{sec-discrep} we proved that the log discrepancies of a
Grassmannian cone $C \rightarrow B$ are related to the log
discrepancies of the pair $(B,B_{g-1})$.  In this section we describe
a construction which, when combined with results related to inversion
of adjunction, often allow one to prove that $(B,B_{g-1})$ is log
canonical (resp. Kawamata log terminal, canonical) in a neighborhood
of a Cartier divisor $S\subset B$ by proving that $(S,S_{g-1})$ is log
canonical (resp. Kawamata log terminal, canonical).  The construction
is not strictly necessary for our application to the singularities of
Grassmannian cones, but it is an obvious variation of the
constructions in the last sections and it makes the statements of
certain lemmas more natural.  Throughout this section we assume that
$\mc{E}$ is torsion-free of rank $e=f-g > 0$.

\medskip\noindent
Let the pair $(\rho:C' \rightarrow B, \beta:\rho^* \mc{F} \rightarrow
\mc{Q}')$ denote the Grassmannian bundle of rank $g$ locally free
quotients of $\mc{F}$.  Our first construction equates the log
discrepancies of the pair $(B,B_{g-1})$ with the log discrepancies of
a pair $(C',\mc{D}_{\phi})$ for a Cartier divisor $\mc{D}_{\phi}
\subset C'$.  Morally, the construction is a version of the principle
that a pair $(B,Z)$ is log canonical (resp. Kawamata log terminal,
canonical) if for a general hypersurface $H\subset B$ containing $Z$,
the pair $(B,H)$ is log canonical (resp. Kawamata log terminal,
canonical).  Indeed, if we locally trivialize $\mc{F}$ so that $C'
\cong B \times \text{Grass}(g,f)$, then the fibers of $\mc{D}_{\phi}
\rightarrow \text{Grass}(g,f)$ considered as subvarieties of $B$ are
hypersurfaces containing $B_{g-1}$.  Therefore one expects that the
general fiber of the pair $(C',\mc{D}_{\phi})$ is log canonical,
etc. iff $(B,B_{g-1})$ is log canonical, etc.  In this case, it is
even true that $(C',\mc{D}_{\phi})$ is log canonical, etc. iff
$(B,B_{g-1})$ is log canonical, etc.

\medskip\noindent
\begin{notat} \label{notat-Dphi} 
Denote by $\delta:\rho^* \mc{G} \rightarrow \mc{Q}'$ the composition
\begin{equation}
\delta: \rho^* \mc{G} \xrightarrow{\rho^* \phi} \rho^* \mc{F}
\xrightarrow{\beta} \mc{Q}'.
\end{equation}
Both $\rho^* \mc{G}$ and $\mc{Q}'$ have rank $g$, so $\delta$
induces a well-defined morphism of invertible sheaves
$\text{det}(\delta):\rho^*\text{det}(\mc{G}) \rightarrow
\OO_{C'}(1)$.  Denote by $\mc{D}_{\phi} \subset C'$ the zero
scheme of this morphism.
\end{notat}

\begin{lem} \label{lem-D} 
The Cartier divisor $\mc{D}_{\phi}$ is irreducible and generically
reduced, and the projection morphism $\pi:\mc{D}_{\phi} \rightarrow B$
admits a dualizing complex of the form
\begin{equation}
\omega_{\mc{D}_{\phi}/B} = \rho^*\lt( \text{det}(\mc{F})^{\otimes(r-1)}
  \otimes_{\OO_B} \text{det}(\mc{E}) \rt)\otimes_{\OO_{C'}}
  \OO_{C'}(-(f-1))|_{\mc{D}_{\phi}} [r(f-r) - 1].
\end{equation}
If $B$ is Cohen-Macaulay (resp. Gorenstein) then also $\mc{D}_{\phi}$
is Cohen-Macaulay (resp. Gorenstein).  If $B$ is Cohen-Macaulay and
$\text{codim}_B(B_{g-1}) \geq 3$, then $\mc{D}_{\phi}$ is normal.
\end{lem}

\medskip\noindent
\begin{proof} 
The proof is the same sort of argument as in the proofs of
Lemma~\ref{lem-codims} and Lemma~\ref{lem-discrep}.  The details are
left to the reader.  There is one extra detail in the proof of the
last claim: over $B-B_{g-1}$, the divisor $\mc{D}_{\phi}$ is not
necessarily smooth.  However the singular locus of $\mc{D}_{\phi}$ is
the locus where $\text{delta}$ has rank at most $g-2$, and this has
codimension at least $4$ in $\rho^{-1}(B - B_{g-1})$.  So the singular
locus of $\mc{D}_{\phi} \cap \rho^{-1}(B-B_{g-1})$ has codimension at
least $3$ in $\mc{D}_{\phi}$.  If $\text{codim}_B(B_{g-1}) \geq 3$,
then the codimension of $\mc{D}_{\phi} \cap \rho^{-1}(B_{g-1})$ in
$C'$, i.e. the codimension of $\rho^{-1}(B_{g-1})$ in $C'$, is
$\text{codim}_B(B_{g-1}) \geq 3$.  Therefore the codimension of
$\rho^{-1}(B_{g-1})$ in $\mc{D}_{\phi}$ is $\text{codim}_B(B_{g-1}) -
1 \geq 2$.  So the singular locus of $\mc{D}_{\phi}$ has codimension
at least $2$ in $\mc{D}_{\phi}$, and by Serre's criterion we conclude
that $\mc{D}_{\phi}$ is normal.
\end{proof}

\begin{rmk} \label{rmk-D} 
The last condition for $\mc{D}_{\phi}$ to be normal is not a necessary
condition.  It seems certain that $\mc{D}_{\phi}$ is normal provided
that for every codimension $2$ point $\eta$ of $B$ contained in
$B_{g-1}$, $B$ is regular at $\eta$ and $\mc{I}_{g-1}\cdot
\OO_{B,\eta}$ equals the maximal ideal $\mathfrak{p}_\eta$.  In fact
it seems likely that $\mc{D}_{\phi}$ is normal provided that for every
codimension $2$ point $\eta$ of $B$ contained in $B_{g-1}$, $B$ is
regular at $\eta$ and $\mc{I}_{g-1} \cdot \OO_{B,\eta}/
\mathfrak{p}_{\eta}^2$ inside of
$\mathfrak{p}_\eta/\mathfrak{p}_\eta^2$ has dimension at least $1$ as
a $\kappa(\eta)$-vector space.  Examples of such maps are when
$B=\AAA^2$ and $\phi:\OO_{\AAA^2} \rightarrow \OO_{\AAA^2}^{\oplus 2}$
is the map with matrix $(x,y^m)^\dagger$ for $m\geq 1$.
\end{rmk}

\medskip\noindent 
Let $u:\widetilde{B} \rightarrow B$ be a resolution of $\mc{E}$, and
let $\widetilde{\mc{E}}$ and $\widetilde{\phi}: \widetilde{\mc{G}}
\rightarrow u^* \mc{F}$ be as in Notation~\ref{notat-res}.  Inside of
the fiber product $\widetilde{B} \times_B C'$, we have the Cartier
divisor $\mc{D}_{\widetilde{\phi}}$.  The projection morphism
$\text{pr}_2:\widetilde{B} \times_B C' \rightarrow C'$ maps
$\mc{D}_{\widetilde{\phi}}$ onto $\mc{D}_{\phi}$.  Denote by
$w:\mc{D}_{\widetilde{\phi}} \rightarrow \mc{D}_{\phi}$ the induced
morphism.

\begin{lem} \label{lem-relcan2} 
The inverse image Cartier divisor $\text{pr}_2^*\mc{D}_\phi$ has
support in the divisor $\mc{D}_{\widetilde{\phi}} \cup
\text{pr}_1^{-1}(E_1 \cup \dots E_k)$.  Moreover we have
\begin{equation}
K_{\widetilde{B} \times_B C'} + \mc{D}_{\widetilde{\phi}} =
\text{pr}_2^*\lt( K_{C'} + \mc{D}_{\phi} \rt) +   
\sum_{i=1}^k \lt(a(E_i;B,B_{g-1}) -1\rt) \text{pr}_1^* E_i.
\end{equation}
\end{lem}

\begin{proof}
Over $\widetilde{B} - (E_1 \cup \dots \cup E_k)$, it is clear that
$\text{pr}_2$ is an isomorphism.  Therefore we have that
\begin{equation}
K_{\widetilde{B} \times_B C'} + \mc{D}_{\widetilde{\phi}} =
\text{pr}_2^*\lt( K_{C'} + \mc{D}_{\phi} \rt) +   
\sum_{i=1}^k \lt(a_i -1\rt) \text{pr}_1^* E_i
\end{equation}
for some sequence of rational numbers $a_1,\dots,a_k$.  To compute the
integers $a_i$, we first restrict to $\mc{D}_{\widetilde{\phi}}$.  By
adjunction, the restriction of $K_{\widetilde{B} \times_B C'} +
\mc{D}_{\widetilde{\phi}}$ to $\mc{D}_{\widetilde{\phi}}$ is just
$K_{\mc{D}_{\widetilde{\phi}}}$.  And the restriction of
$\text{pr}_2^*(K_{C'} + \mc{D}_{\phi})$ equals $w^*K_{\mc{D}_{\phi}}$
(by which we mean the sum of $K_B$ and $C_1(\omega_{\mc{D}_{\phi}/B})$
in case $\mc{D}_{\phi}$ is not normal).  But applying
Lemma~\ref{lem-D} to both $\mc{D}_{\phi}$ and
$\mc{D}_{\widetilde{\phi}}$, we conclude that
\begin{equation}
\begin{array}{ll}
K_{\mc{D}_{\widetilde{\phi}}} - w^*K_{\mc{D}_{\phi}} = \\
\widetilde{\pi}^* \lt[ K_{\widetilde{B}} - u^*K_B +
C_1( \text{det}(\widetilde{\mc{E}}) ) - C_1( \text{det}(\mc{E}) ) \rt].
\end{array} 
\end{equation}
As proved in Lemma~\ref{lem-relcan}, the divisor on the right is just
$\sum_{i=1}^k \lt( a(E_i;B,B_{g-1}) - 1 \rt)\widetilde{\pi}^*E_i$.  So
we conclude that the restriction to $\mc{D}_{\widetilde{\phi}}$ of
$\sum_{i=1}^k \lt(a_i - 1\rt) \text{pr}_1^* E_i$ equals the restriction of
$\sum_{i=1}^k \lt( a(E_i;B,B_{g-1}) - 1 \rt) \text{pr}_1^* E_i$.

\medskip\noindent
To finish the argument, we need to prove that the pullback
$\text{Pic}(\widetilde{B}) \rightarrow
\text{Pic}(\mc{D}_{\widetilde{\phi}})$ is injective.  We prove this by
breaking the argument up into two possible cases depending on whether
$f-g \leq g$ or $f-g > g$.

\medskip\noindent
Suppose first that $f-g \leq g$.  Form the Grassmannian bundle
$(\sigma:C'' \rightarrow \widetilde{B}, \epsilon:\sigma^*
\widetilde{\mc{G}} \rightarrow \mc{Q}'')$ parametrizing rank $g-(f-g)$
locally free quotients of $\widetilde{\mc{G}}$.  The coproduct of
$\epsilon:\sigma^* \widetilde{\mc{G}} \rightarrow \mc{Q}''$ and
$\sigma^*\widetilde{\phi}: \sigma^* \widetilde{\mc{G}} \rightarrow
\sigma^* u^*\mc{F}$ gives a surjective morphism of sheaves $\sigma^*
u^*\mc{F} \rightarrow \mc{R}$ where $\mc{R}$ is locally free of rank
$\text{rank}(\mc{F}) - \text{rank}(\widetilde{\mc{G}}) +
\text{rank}(\mc{Q}'') = g$.  There is an induced morphism from $C''$
to $\widetilde{B} \times_B C$ which is compatible with the projection
to $\widetilde{B}$.  And the image is contained in
$\mc{D}_{\widetilde{\phi}}$.  Since $\sigma:C'' \rightarrow
\widetilde{B}$ is a Grassmannian bundle, the pullback map on Picard
groups is injective.  And this map factors through the pullback map on
Picard groups from $\widetilde{B}$ to $\mc{D}_{\widetilde{\phi}}$.
Therefore the pullback map $\text{Pic}(\widetilde{B}) \rightarrow
\text{Pic}(\mc{D}_{\widetilde{\phi}})$ is injective.

\medskip\noindent
Finally, suppose that $f-g > g$.  In this case let $(\sigma:C''
\rightarrow \widetilde{B}, \epsilon: \sigma^*(\widetilde{\mc{E}})
\rightarrow \mc{Q}'')$ be the Grassmannian bundle parametrizing rank
$g$ locally free quotients of $\widetilde{\mc{E}}$.  There is an
obvious closed immersion of $C''$ into $\widetilde{B} \times_B C'$.
And the image clearly lies in $\mc{D}_{\widetilde{\phi}}$.  As in the
last paragraph, this implies that the pullback map
$\text{Pic}(\widetilde{B}) \rightarrow
\text{Pic}(\mc{D}_{\widetilde{\phi}})$ is injective.
\end{proof}

\medskip\noindent
The scheme $\mc{D}_{\widetilde{\phi}} \rightarrow \widetilde{B}$ is
typically not smooth, so we do not yet have a log resolution of
$(C',\mc{D}_{\phi})$.  The construction of a log resolution of
$(\widetilde{B} \times_B C', \mc{D}_{\widetilde{\phi}})$ is
essentially equivalent to the construction in Section~\ref{sec-deter}
(and provides some justification for the tedious arguments of that
section).

\medskip\noindent  
\begin{notat} \label{notat-Co} 
Denote by $C^{(0)}$ the fiber product $C^{(0)} = \widetilde{B}
\times_B C'$ and denote by $M^{(0)} =
M^{(0)}(\widetilde{B},\widetilde{\mc{G}}, \OO_{\widetilde{B}}^{\oplus
g})$ the scheme constructed in Section~\ref{sec-deter}.  Denote by
$p^{(0)}:T^{(0)} \rightarrow C^{(0)}$ the $\text{GL}_g$-torsor
parametrizing sheaf isomorphisms $\text{pr}_2^*\mc{Q}' \rightarrow
\OO_{C^{(0)}}^{\oplus g}$ (with the obvious left $\text{GL}_g$-action)
and denote by $\lambda: (p^{(0)})^* \text{pr}_2^* \mc{Q}' \rightarrow
\OO_{T^{(0)}}^{\oplus g}$ the universal isomorphism.  Denote by
$\epsilon$ the composition of $(p^{(0)})^*\delta:(p^{(0)})^*
\text{pr}_1^* \widetilde{\mc{G}} \rightarrow
(p^{(0)})^*\text{pr}_2^*\mc{Q}'$ with $\lambda$.  Denote by
$q^{(0)}:T^{(0)} \rightarrow M^{(0)}$ the morphism of
$\widetilde{B}$-schemes induced by $\epsilon$.  Observe that this
morphism is equivariant for the obvious $\text{GL}_g$-action on
$M^{(0)}$.
\end{notat}

\begin{lem} \label{lem-qsmooth} 
Denote by $g'$ the maximum of $0$ and $2g-f$.
\begin{enumerate}
\item
The image of $q^{(0)}$ equals $M^{(0)} - M^{(0)}_{g'-1}$.
\item
The morphism $q^{(0)}: T^{(0)} \rightarrow \lt( M^{(0)} -
M^{(0)}_{g'-1} \rt)$ factors as an open immersion into a torsor over
$M^{(0)}$ for the vector bundle over $\widetilde{B}$ associated to
$\textit{Hom}_{\OO_{\widetilde{B}}} \lt( \widetilde{\mc{E}},
\OO_{\widetilde{B}}^{\oplus g} \rt)$.  In particular, $q^{(0)}$ is
smooth.
\item
The inverse image scheme $(p^{(0)})^{-1}(\mc{D}_{\widetilde{\phi}})$
equals the inverse image scheme $(q^{(0)})^{-1}(M^{(0)}_{g-1})$.
\end{enumerate}
\end{lem}

\begin{proof}
Item $(1)$ follows on considering the intersection of the subbundle
$\text{pr}_1^*\widetilde{\mc{G}} \subset \text{pr}_1^* u^*\mc{F}$ with
the kernel of $\text{pr}_1^* u^*\mc{F} \rightarrow
\text{pr}_2^*\mc{Q}'$.  The first subbundle has rank $g$ at every
point, and the second has rank $f-g$.  Therefore the maximal possible
intersection is $f-g$ if $f-g \leq g$, and $g$ otherwise. So the
minimal possible rank of $\epsilon$ is $g-(f-g) = 2g-f$ if $2g-f \geq
0$, and $0$ otherwise, i.e. the minimum possible rank is $g'$.  On the
other hand, up to composing with an isomorphism $\text{pr}_2^*\mc{Q}'
\rightarrow \OO_{C^{(0)}}^{\oplus g}$, clearly we can obtain any
morphism $\widetilde{\mc{G}} \rightarrow \OO_{\widetilde{B}}^{\oplus
g}$ as a fiber of $\epsilon$ over a geometric point of $T^{(0)}$.
Therefore $q^{(0)}: T^{(0)} \rightarrow \lt( M^{(0)} - M^{(0)}_{g'-1}
\rt)$ is surjective.

\medskip\noindent
The torsor over $M^{(0)}$ is simply $M^{(0)}(\widetilde{B}, u^*\mc{F},
\OO_{\widetilde{B}}^{\oplus g})$.  The open immersion from $T^{(0)}$
to this scheme is clear.

\medskip\noindent
By construction, $\mc{D}_{\widetilde{\phi}}$ is the scheme determined
by the determinant of $\delta$.  But the pullback $(p^{(0)})^*\delta$
equals $(q^{(0)})^* \phi$ by construction.  Therefore the inverse
image of $\mc{D}_{\widetilde{\phi}}$ under $p^{(0)}$ is precisely the
inverse image of $M^{(0)}_{g-1}$ under $q^{(0)}$.
\end{proof}

\medskip\noindent
\begin{notat} \label{notat-Tk} 
Denote by
\begin{equation}
T^{(g)} \rightarrow T^{(g-1)}  \rightarrow \dots \rightarrow T^{(1)}
\rightarrow T^{(0)}
\end{equation}
the sequence of morphisms which are obtained via base-change by
$q^{(0)}$ from the sequence of morphisms
\begin{equation}
M^{(g)} \rightarrow M^{(g-1)} \rightarrow \dots \rightarrow M^{(1)}
\rightarrow M^{(0)}
\end{equation}
constructed in Lemma~\ref{lem-exist}.
\end{notat}

\medskip\noindent
In particular all of the schemes $T^{(k)}$ have a natural
$\text{GL}_g$-action and the morphisms are all
$\text{GL}_g$-equivariant.  For each $i=1,\dots, g$ the composition
$T^{(i)} \rightarrow T^{(0)}$ is equivalent to the blowing up of an
ideal sheaf $\mc{J}^{(i)}_T$ on $T^{(0)}$.  Moreover, this ideal sheaf
is $\text{GL}_g$-equivariant.  Therefore it is of the form
$(p^{(0)})^{-1} \mc{J}^{(i)}_C$ for some ideal sheaf $\mc{J}{(i)}_C$
on $C^{(0)}$.

\medskip\noindent
\begin{notat} \label{notat-Ck} 
Denote by
\begin{equation}
C^{(g)} \rightarrow C^{(g-1)} \rightarrow \dots \rightarrow C^{(1)}
\rightarrow C^{(0)}
\end{equation}
the sequence of morphisms
which are the blowing ups of the ideal sheaves $\mc{J}{(i)}_C$.  
\end{notat}

\medskip\noindent
For each $0 \leq r < s \leq g$, we have a Cartesian diagram
\begin{equation}
\begin{CD}
T^{(s)} @>>> T^{(r)} \\
@VVV  @VVV \\
C^{(s)} @>>> C^{(r)}
\end{CD}
\end{equation}
where the vertical arrows are $\text{GL}_g$-torsors, and we have a
Cartesian diagram
\begin{equation}
\begin{CD}
T^{(s)} @>>> T^{(r)} \\
@VVV  @VVV \\
M^{(s)} @>>> M^{(r)}
\end{CD}
\end{equation}
where the vertical arrows are open subsets of torsors for a smooth
group scheme.  By Proposition~\ref{prop-smooth} $M^{(g)} \rightarrow
M^{(0)}$ gives a log resolution of the pair $(M^{(0)},
M^{(0)}_{g-1})$. So we conclude that $C^{(g)} \rightarrow C^{(0)}$
gives a log resolution of the pair
$(C^{(0)},\mc{D}_{\widetilde{\phi}})$.  Moreover $C^{(0)} \rightarrow
\widetilde{B}$ is smooth.

\medskip\noindent
\begin{notat} \label{notat-Fk} 
Denote by $F_0,\dots, F_{g-1}$ the Cartier divisors on $C^{(g)}$
corresponding to the divisors $E^{(g)}_0, \dots E^{(g)}_{g-1}$ on
$M^{(g)}$.  Of course $F_i = \emptyset$ for $i < g'$.
\end{notat}

\begin{prop}  \label{prop-res2} 
There exists a log resolution $t:C^{(g)} \rightarrow \widetilde{B}
\times_B C'$ of the pair $\lt( \widetilde{B} \times_B
C',\mc{D}_{\widetilde{\phi}} \rt)$ with exceptional locus $F_{g'} \cup
\dots \cup F_{g-2}$ (where $g' = \text{max}(0,2g-f)$ by definition)
satisfying the following properties:
\begin{enumerate}
\item 
The morphism $\text{pr}_1 \circ t:C^{(g)} \rightarrow \widetilde{B}$
is smooth, and the intersection of every fiber with $F_{g'} \cup \dots
\cup F_{g-2} \cup F_{g-1}$ is a simple normal crossings divisor.
\item 
The morphism $\text{pr}_2 \circ t:C^{(g)} \rightarrow C'$ is a log
resolution of the pair $(C',\mc{D}_\phi)$ with exceptional locus
$F_{g'} \cup \dots \cup F_{g-2} \cup (\text{pr}_2 \circ t)^{-1} \lt(
E_1 \cup \dots \cup E_k \rt)$ and such that the strict transform of
$\mc{D}_{\phi}$ is the divisor $F_{g-1}$.
\item 
We have an equivalence of $\QQ$-Cartier divisors on $C^{(g)}$
\begin{equation}
\begin{array}{l}
K_{C^{(g)}} - (\text{pr}_2 \circ t)^*(K_{C'} + \mc{D}_\phi) = \\
-F_{g-1} + \sum_{j=g'}^{g-2}\lt( (g-1-j)(g-j) -1\rt) F_j + \\
\sum_{i=1}^k
\lt(a(E_i;B,B_{g-1}) - 1 \rt) (\text{pr}_1 \circ t)^* E_i
\end{array}
\end{equation}
\item 
The log discrepancy of $(C',\mc{D}_{\phi})$ equals the minimum of $1$
and the log discrepancy of $(B,B_{g-1})$.  In particular, every
exceptional divisor for $(B,B_{g-1})$ gives rise to an exceptional
divisor for $(C',\mc{D}_{\phi})$.
\item 
The pair $(C',\mc{D}_{\phi})$ is log canonical (resp. purely log
terminal, canonical) iff the pair $(B,B_{g-1})$ is log canonical
(resp. purely log terminal, canonical).
\item 
Assume that $\mc{D}_{\phi}$ is normal.  Then the total discrepancy of
$(\mc{D}_{\phi},\emptyset)$ equals the total discrepancy of
$(B,B_{g-1})$.
\item 
Assume that $\mc{D}_{\phi}$ is normal.  Then
$(\mc{D}_{\phi},\emptyset)$ is log canonical (resp. Kawamata log
terminal, canonical) iff the pair $(B,B_{g-1})$ is log canonical
(resp. purely log terminal, canonical).
\item 
Assume that $\mc{D}_{\phi}$ is normal.  If $(B,B_{g-1})$ is terminal,
then $\mc{D}_{\phi}$ is terminal.  If for every exceptional divisor
$E_i$ we have $a(E_i;B,B_{g-1}) \neq 1$, then the converse also holds.
\end{enumerate}
\end{prop}

\begin{proof}
We may check Item $(1)$ after performing the smooth, surjective
base-change by $T^{(g)} \rightarrow C^{(g)}$.  And again by smooth
surjective base-change, the results on $T^{(g)}$ is equivalent to the
statement that $M^{(g)}$ is smooth over $\widetilde{B}$ and every
fiber intersects $E^{(g)}_{g'} \cup \dots \cup E^{(g)}_{g-1}$ in a
simple normal crossings divisor.  This follows from Item $(1)$ and
Item $(2)$ of Proposition~\ref{prop-smooth}.

\medskip\noindent
As mentioned, $t$ is a log resolution of $\mc{D}_{\widetilde{\phi}}$.
Moreover the divisor $F_{g'} \cup \dots F_{g-2}$ is flat over
$\widetilde{B}$ and intersects every fiber in a simple normal
crossings divisor.  Therefore $F_{g'} \cup \dots F_{g-2} \cup
(\text{pr}_2 \circ t)^{-1} \lt( E_1 \cup \dots \cup E_k \rt)$ is a
simple normal crossings divisor.  It follows that $\text{pr}_2\circ t$
is a log resolution of $(C',\mc{D}_{\phi})$.  This proves Item $(2)$.

\medskip\noindent
Item $(3)$ follows from Corollary~\ref{cor-mults} and
Lemma~\ref{lem-relcan2}.  Item $(4)$ and Item $(5)$ follow immediately
from Item $(3)$.

\medskip\noindent
Next we consider Item $(6)$.  Observe that $F_{g-1} \rightarrow
\mc{D}_{\phi}$ is a resolution of singularities and we have
\begin{equation}
\begin{array}{l}
K_{F_{g-1}} - (w\circ t)^* K_{\mc{D}_{\phi}} = \\
\sum_{j=g'}^{g-2} \lt( (g-1-j)(g-j) - 1 \rt) F_j|_{F_{g-1}} + \\
\sum_{i=1}^k \lt( a(E_i; B,B_{g-1}) - 1 \rt) (\text{pr}_1 \circ t)^*
E_i.
\end{array}
\end{equation}
By the same sort of argument as in Lemma~\ref{lem-relcan2}, all of the
relevant divisor classes are linearly independent on $F_{g-1}$.  Since
the coefficients $(g-1-j)(g-j)$ are at least $2$ for $j
=g',\dots,g-2$, we conclude that the total discrepancy of
$\mc{D}_{\phi}$ equals the total discrepancy of $(B,B_{g-1})$.  As
always, it is possible that some exceptional divisors of $(B,B_{g-1})$
do not give rise to exceptional divisors of $\mc{D}_{\phi}$.  Item
$(7)$ and Item $(8)$ follow immediately from Item $(6)$.
\end{proof}

\medskip\noindent
\begin{cor} \label{cor-spec} 
Let $S \subset B$ be an irreducible Cartier divisor, denote by
$\phi_S$ the restriction of $\phi$ to $S$, and assume that
$\mc{D}_{\phi_S}$ is irreducible.
\begin{enumerate}
\item
Suppose that $S$ is Kawamata log terminal.  Then $(S,S_{g-1})$ is log
canonical iff there exists an open subscheme $U\subset B$ containing
$S$ such that $(U,S+U_{g-1})$ is log canonical.  In particular, if
$(S,S_{g-1})$ is log canonical, then $(U,U_{g-1})$ is log canonical.
\item 
Suppose that $\mc{D}_{\phi_S}$ is irreducible and normal. Then
$(S,S_{g-1})$ is Kawamata log terminal iff there exists an open
subscheme $U\subset B$ containing $S$ such that $(U,S+U_{g-1})$ is
purely log terminal.  In particular, if $(S,S_{g-1})$ is Kawamata log
terminal, then $(U,U_{g-1})$ is Kawamata log terminal.
\item 
Suppose that $B$ is Gorenstein and that $\mc{D}_{\phi_S}$ is
irreducible and normal.  Then $(S,S_{g-1})$ is canonical iff there
exists an open subscheme $U\subset B$ containing $S$ such that
$(U,U_{g-1})$ is canonical.  In particular, if $(S,S_{g-1})$ is
canonical then $(U,U_{g-1})$ is canonical.
\end{enumerate}
\end{cor}

\medskip\noindent
\begin{proof}
For Item $(1)$, first observe that $\rho^{-1}(S) \subset C'$ is also
Kawamata log terminal since $\rho$ is smooth.  By
Proposition~\ref{prop-res2}, $(S,S_{g-1})$ is log canonical iff
$(\rho^{-1}(S), \mc{D}_{\phi_S})$ is log canonical.  By a similar
argument, $(B,S+B_{g-1})$ is log canonical iff
$(C',\rho^{-1}(S)+\mc{D}_{\phi})$ is log canonical.  By
\cite[Thm. 7.5.2]{KollarPairs}, $(\rho^{-1}(S),\mc{D}_{\phi_S})$ is
log canonical iff $(C',\rho^{-1}(S) + \mc{D}_{\phi})$ is log canonical
near $\rho^{-1}(S)$.  Therefore $(S,S_{g-1})$ is log canonical iff
$(B,S+B_{g-1})$ is log canonical near $S$.  So Item $(1)$ holds.

\medskip\noindent
For Item $(2)$ we just combine Item $(7)$ of
Proposition~\ref{prop-res2} with ~\cite[Thm. 7.5.1]{KollarPairs}.  For
Item $(3)$ we just combine Item $(7)$ of Proposition~\ref{prop-res2}
with ~\cite{Stevens88} (see also ~\cite[Thm. 7.9]{KollarPairs}).
\end{proof}

\begin{rmk} \label{rmk-spec} 
If one further assumes that $B$ is a local complete intersection
scheme, then one can also use the results of Ein and Musta\c{t}\v{a}
~\cite{EM} to prove Corollary~\ref{cor-spec} and to relate the minimal
log discrepancy of $(S,S_{g-1})$ to the minimal log discrepancy of
$(B,B_{g-1})$.
\end{rmk}

\section{Deformation to the normal cone} \label{sec-dnc}
\medskip\noindent
Corollary~\ref{cor-spec} allows us to deduce results about
$(B,B_{g-1})$ by analyzing the restriction of $\phi$ to an irreducible
Cartier divisor $S$.  In applications it is also natural to restrict
$\phi$ to an irreducible closed subvariety $Y\subset B$ which is not
necessarily a Cartier divisor; in particular, if $r$ is the smallest
integer such that $B_r \neq \emptyset$, then one natural choice is to
take $Y$ to be an irreducible component of $B_r$.  In this section we
show how to reduce the general case to the case of a Cartier divisor
by using deformation to the normal cone.  We briefly review the
discussion of deformation to the normal cone from ~\cite[Chapter
5]{F}.  All of the unproved assertions regarding deformation to the
normal cone which we use are proved there.

\medskip\noindent
The following setup is a little more general than we will need.  In
this section we consider a morphism of locally free sheaves
$\phi:\mc{G} \rightarrow \mc{F}$, but we do not assume that $f > g$.
Suppose that $Y \subset B$ is a closed subscheme with ideal sheaf
$\mc{J}$.  Denote by $\mc{K}_Y$ and $\mc{E}_Y$ the kernel and cokernel
respectively of the map of $\OO_Y$-modules
\begin{equation}
\phi\otimes \text{Id}: \mc{G} \otimes_{\OO_B} \OO_Y \rightarrow \mc{F}
\otimes_{\OO_B} \OO_Y.
\end{equation}
In particular $\mc{E}_Y$ is simply $\mc{E} \otimes_{\OO_B} \OO_Y$.
Consider the following commutative diagram with exact rows
\begin{equation} \label{eqn-snake1} 
\begin{CD}
0 @>>> \mc{G}\otimes_{\OO_B} \mc{J}/\mc{J}^2 @>>> \mc{G} \otimes_{\OO_B}
 \OO_B/\mc{J}^2 @>>> \mc{G} \otimes_{\OO_B} \OO_Y @>>> 0 \\
& & @V \phi_1 VV @VV\phi_2 V @VV\phi_3 V \\
0 @>>> \mc{F}\otimes_{\OO_B} \mc{J}/\mc{J}^2 @>>> \mc{F} \otimes_{\OO_B}
 \OO_B/\mc{J}^2 @>>> \mc{F} \otimes_{\OO_B} \OO_Y @>>> 0 
\end{CD}
\end{equation}
Each map $\phi_i$ is just $\phi \otimes \text{Id}$.  Since tensor
product is right exact, the cokernel of $\phi_1$ is just $\mc{E}_Y
\otimes_{\OO_Y} \mc{J}/\mc{J}^2$.  And by definition the kernel of $\phi_3$ is
$\mc{K}_Y$.  By the Snake Lemma, there is an induced connecting map
from $\text{Ker}(\phi_3)$ to $\text{Coker}(\phi_1)$.

\medskip\noindent
\begin{defn} \label{defn-conn} 
The \emph{connecting map}, denoted $\theta = \theta_{\phi,Y}: \mc{K}_Y
\rightarrow \mc{E}_Y \otimes_{\OO_Y} \mc{J}/\mc{J}^2$, is the map of
$\OO_Y$-modules induced by the Snake Lemma as above.  The
\emph{induced map} of the connecting map is the map $\theta'_{\phi,Y}:
\textit{Hom}_{\OO_Y}(\mc{E}_Y,\mc{K}_Y) \rightarrow \mc{J}/\mc{J}^2$
induced from $\theta_{\phi,Y}$ by adjointness of Hom and tensor
product.
\end{defn}

\medskip\noindent
\begin{rmk} \label{rmk-conn} 
From now on we will assume that $\mc{E}_Y$ is locally free, which
implies that also $\mc{K}_Y$ is locally free. In our later
applications, we will also always have that $\mc{J}/\mc{J}^2$ is
locally free (but we do not make this assumption in the remainder of
this section).
\end{rmk}

\begin{lem} \label{lem-conn} 
\begin{enumerate}
\item 
For the transpose $\phi^\dagger$, the kernel of $\phi^\dagger|_Y$ is
$\mc{E}_Y^\vee$, the cokernel of $\phi^\dagger|_Y$ is $\mc{K}_Y^\vee$,
and the connecting map $\theta'_{\phi^\dagger,Y}$ is identified with
$\theta'_{\phi,Y}$ under the canonical isomorphism
$$ \textit{Hom}_{\OO_Y}(\mc{K}^\vee_Y,\mc{E}^\vee_Y) \cong
  \textit{Hom}_{\OO_Y}(\mc{E}_Y, \mc{K}_Y).
$$
\item 
Let $\phi':\mc{G}' \rightarrow \mc{F}'$ be a second morphism of
locally free sheaves on $B$ such that $\phi'|_Y$ has locally free
kernel and cokernel $\mc{K}'_Y$ and $\mc{E}'_Y$ and consider
$\phi\oplus \phi': \mc{G} \oplus \mc{G}' \rightarrow \mc{F} \oplus
\mc{F}'$.  The kernel of $(\phi \oplus \phi')|_Y$ is $\mc{K}_Y \oplus
\mc{K}'_Y$, the cokernel of $(\phi \oplus \phi')|_Y$ is $\mc{E}_Y
\oplus \mc{E}'_Y$, and $\theta_{\phi\oplus \phi',Y}$ equals
$\theta_{\phi,Y} \oplus \theta_{\phi',Y}$ via the canonical
isomorphisms.
\item 
Let $\phi'$ be given as above and consider $\phi\otimes \phi':\mc{G}
\otimes_{\OO_B} \mc{G}' \rightarrow \mc{F} \otimes_{\OO_B} \mc{F}'$.
The kernel of $\phi\otimes \phi'$ is the surjective image of $\lt(
\mc{K}_Y \otimes_{\OO_Y} \mc{G}'|_Y \rt) \oplus \lt( \mc{G}|_Y
\otimes_{\OO_Y} \mc{K}'_Y \rt)$, the cokernel of $\phi \otimes \phi'$
is a subsheaf of $\lt( \mc{E}_Y \otimes_{\OO_Y} \mc{F}'|_Y \rt) \oplus
\lt( \mc{F}|_Y \otimes_{\OO_Y} \mc{E}'_Y \rt)$, and
$\theta_{\phi\otimes \phi',Y}$ is the unique morphism compatible with
$( \theta_{\phi,Y} \otimes \phi' ) \oplus (\phi \otimes
\theta_{\phi',Y})$.
\item 
Let $\phi'$ be given as above, and let $\psi_{\mc{G}}:\mc{G}
\rightarrow \mc{G}'$ and $\psi_{\mc{F}}:\mc{F} \rightarrow \mc{F}'$ be
morphisms of $\OO_B$-modules such that we have a commutative diagram:
\begin{equation}
\begin{CD}
\mc{G} @> \phi >> \mc{F} \\
@V \psi_{\mc{G}} VV @VV \psi_{\mc{F}} V \\
\mc{G}' @> \phi'>> \mc{F}' 
\end{CD}
\end{equation}
There are unique morphisms $\psi_{K}:\mc{K}_Y \rightarrow \mc{K}'_Y$
and $\psi_{E}:\mc{E}_Y \rightarrow \mc{E}'_Y$ such that the following
diagram commutes:
\begin{equation}
\begin{CD}
0 @>>> \mc{K}_Y @>>> \mc{G}|_Y @> \phi|_Y >> \mc{F}|_Y @>>> \mc{E}_Y
@>>> 0 \\
& & @V \psi_{\mc{K}} VV @V \psi_{\mc{G}} VV @VV \psi_{\mc{F}} V @VV
\psi_{\mc{E}} V \\
0 @>>> \mc{K}'_Y @>>> \mc{G}'|_Y @> \phi'|_Y >> \mc{F}'|_Y @>>> \mc{E}'_Y
@>>> 0 \\
\end{CD}
\end{equation}
Moreover, the diagram of connecting maps commutes:
\begin{equation}
\begin{CD}
\mc{K}_Y @> \theta_{\phi,Y} >> \mc{E}_Y \otimes_{\OO_Y}
\mc{J}/\mc{J}^2 \\
@V \psi_{\mc{K}} VV @VV \psi_{\mc{E}}\otimes \text{Id} V \\
\mc{K}'_Y @> \theta_{\phi',Y} >> \mc{E}'_Y \otimes_{\OO_Y}
\mc{J}/\mc{J}^2. 
\end{CD}
\end{equation}
\end{enumerate}
\end{lem}

\begin{proof}
Each of these follows by some simple diagram-chasing.  The details are
left to the reader.
\end{proof}

\medskip\noindent
Now we recall the construction of ``deformation to the normal cone''
as discussed in ~\cite[Chapter 5]{F}.  Form the product $B \times
\PP^1$ and consider the closed subscheme $Y \times \{\infty \} \subset
B \times \PP^1$.  The ideal sheaf of this subscheme is
\begin{equation}
\mc{J}' = 
\text{pr}_1^{-1}(\mc{J}) +
\text{pr}_2^{-1}(\OO_{\PP^1}(-\underline{\infty})).
\end{equation}
This decomposition of the ideal sheaf yields a decomposition of the
Rees algebra
\begin{equation}
\begin{array}{l}
\oplus_{n=0}^\infty (\mc{J}')^n/(\mc{J}')^{n+1} \cong   \\
\text{pr}_1^*\lt( \oplus_{n=0}^\infty (\mc{J})^n/ \mc{J}^{n+1} \rt)
\otimes \text{pr}_2^* \lt( \oplus_{n=0}^\infty
\OO_{\PP^1}(-n\underline{\infty})/
\OO_{\PP^1}(-(n+1)\underline{\infty}) \rt).
\end{array}
\end{equation}
The relative Spec of the Rees algebra is the normal cone.  If the
normal cone is the symmetric algebra of a locally free sheaf, the
normal cone is called the normal bundle.  In our case the
decomposition above gives an isomorphism of the normal cone
$C_{Y\times\{ \infty \} }(B\times \PP^1)$ with the fiber product
$\text{pr}_1^* C_{Y} B \times_{Y \times \{ \infty \} } \text{pr}_2^*
N_{\{\infty\} } \PP^1$.

\medskip\noindent Let us denote $C= \text{pr}_1^* C_Y B$ and $C'= C_{Y
\times \{\infty \} } (B \times \PP^1)$.  Of course the normal bundle
$N_{\{ \infty \} } \PP^1$ is just the trivial rank $1$ vector bundle,
which is denoted by $\mathbf{1}$ in ~\cite{F}.  Also our isomorphism
respects the $\mathbb{G}_m$-actions induced by the grading of the
algebras.  Therefore, in the notation of ~\cite{F}, we have an
equivalence of cones $C' \cong C \oplus \mathbf{1}$.

\medskip\noindent
Now let $u:M \rightarrow B \times \PP^1$ be the blowing up of $B
\times \PP^1$ along $Y \times \{ \infty \}$.  Denote by $\varrho: M
\rightarrow \PP^1$ the composition $\text{pr}_2 \circ u$.  This is a
flat morphism.  The preimage of $\AAA^1 = \PP^1 - \{ \infty \}$ is
isomorphic to $B \times \AAA^1$ (compatibly with projection to $B$ and
to $\AAA^1$).  And the Cartier divisor $M_{\infty} =
\varrho^{-1}(\infty)$ is the sum of two effective divisors $B_Y$ and
$\PP(C') = \PP(C \oplus \mathbf{1})$.  Here $B_Y$ is the blowing up of
$B$ along $Y$.  And, as usual, for a cone $K$ the symbol $\PP(K)$
means the relative Proj of the graded algebra associated to $K$.
Denote by $\pi: \PP(C \oplus \mathbf{1}) \rightarrow Y$ the obvious
projection morphism.

\medskip\noindent
The intersection of $B_Y$ and $\PP(C \oplus \mathbf{1})$ is the
exceptional divisor on $B_Y$ and is the ``hyperplane section at
infinity'' $\PP(C)$ in $\PP(C \oplus \mathbf{1})$.  The complement of
the hyperplane section at infinity is identified with the cone $C$
over $Y$.  Finally, there is a closed immersion $\iota: Y \times \PP^1
\rightarrow M$ such that $u\circ \iota: Y \times \PP^1 \rightarrow B
\times \PP^1$ is the obvious closed immersion.  The fiber of $\iota(Y
\times \PP^1)$ over $\infty$ is identifed with the zero section of $C
\subset \PP(C \oplus \mathbf{1})$.

\begin{defn} \label{defn-dnc} 
For a closed subscheme $Y \subset B$, the \emph{deformation to the
normal cone} is the datum $(\varrho: M \rightarrow \PP^1, \iota: Y
\times \PP^1 \hookrightarrow M, B_Y \hookrightarrow M, \PP(C\oplus
\mathbf{1}) \hookrightarrow M)$.  We denote by $\phi_M: \mc{G}_M
\rightarrow \mc{F}_M$ the morphism of locally free sheaves
$u^*\text{pr}_1^* \phi$.
\end{defn}

\medskip\noindent
On $\PP(C\oplus \mathbf{1})$ there is a rank $1$ locally free quotient
$\beta: \pi^*\lt( \text{pr}_1^* \mc{J}/\mc{J}^2 \oplus \OO_Y \rt)
\rightarrow \OO_{\PP(C \oplus \mathbf{1})}(1)$ (which satisfies a
universal property we won't bother stating).  Denote by $\beta_1:
\pi^* \lt( \text{pr}_1^* \mc{J}/\mc{J}^2 \rt) \rightarrow \OO_{\PP(C
\oplus \mathbf{1})}(1)$ and $\beta_2: \OO_{\PP(C\oplus \mathbf{1})}
\rightarrow \OO_{\PP(C \oplus \mathbf{1})}(1)$ the two components of
$\beta$.  Of course the zero scheme of the section $\beta_2$ is
precisely the hyperplane section at infinity $\PP(C) \subset \PP(C
\oplus \mathbf{1})$.  The invertible sheaf $\OO_M(-\PP(C\oplus
\mathbf{1}))|_{\PP(C\oplus \mathbf{1})}$ is canonically isomorphic to
$\OO_{\PP(C \oplus \mathbf{1})}(1)$; the isomorphism is induced by the
isomorphism of ideal sheaves $\text{u}^{-1}\mc{J}' \cong \OO_M(-\PP(C
\oplus \mathbf{1}))$.

\medskip\noindent 
The pullback $\phi_M$ factors through an \emph{elementary transform
up} of $\mc{G}_M$.  To describe this elementary transform, first we
dualize everything.  Consider the adjoint morphism $\phi_M^\dagger:
\mc{F}_M^\vee \rightarrow \mc{G}_M^\vee$.  The restriction of $\phi_M$
to $\PP(C\oplus \mathbf{1})$ is just $\pi^*(\phi|_Y)$.  In particular,
the image of $\lt( \mc{F}_M^\vee \rt)|_{\PP(C \oplus \mathbf{1})}$ is
contained in the kernel of $\pi^*\lt( \mc{G} \otimes_{\OO_B} \OO_Y
\rt)^\vee \rightarrow \pi^*\lt( \mc{K}_Y \rt)^\vee$.  Define the
subsheaf $(\widetilde{\mc{G}} )^\vee \subset \lt( \mc{G}_M \rt)^\vee$
to be the kernel of the surjection
\begin{equation}
\lt( \mc{G}_M \rt)^\vee \rightarrow \lt( \mc{G}_M
\rt)^\vee|_{\PP(C\oplus \mathbf{1})} \cong \pi^*\lt( \mc{G}
\otimes_{\OO_B} \OO_Y \rt)^\vee \rightarrow \pi^* \lt( \mc{K}_Y
\rt)^\vee.
\end{equation}
Then $\phi_M^\dagger$ factors through the subsheaf
$(\widetilde{\mc{G}})^\vee$.  Define $( \widetilde{\phi} )^\dagger:
\mc{F}_M^\vee \rightarrow ( \widetilde{\mc{G}} )^\vee$ to be the
induced map.

\medskip\noindent
\begin{lem} \label{lem-elemtrans} 
Denote by $\widetilde{\mc{G}}$ the dual of
$(\widetilde{\mc{G}})^\vee$, and denote by $\widetilde{\phi}$ the
adjoint of $(\widetilde{\phi})^\dagger$.
\begin{enumerate}
\item
The sheaf $(\widetilde{\mc{G}} )^\vee$ is locally free of rank $g$.
Therefore also $\widetilde{\mc{G}}$ is locally free of rank $g$.
\item
The cokernel of the 
sheaf map $\mc{G}_M
\rightarrow \widetilde{\mc{G}}$ is canonically isomorphic to the
push-forward from $\PP(C\oplus \mathbf{1})$ of the locally free sheaf
$\pi^* \mc{K}_Y \otimes \OO_{\PP(C \oplus \mathbf{1})}(-1)$.  
\item
The restriction of $\mc{G}_M \rightarrow \widetilde{\mc{G}}$ to $\PP(C
\oplus \mathbf{1})$ fits into an exact sequence
\begin{equation}
\begin{array}{lllllll}
0 &  \longrightarrow & \pi^* \mc{K}_Y & \longrightarrow & \pi^*(\mc{G}
 \otimes_{\OO_B} \OO_Y ) & \longrightarrow \\
 & &
\widetilde{\mc{G}}|_{\PP(C \oplus \mathbf{1})} & \longrightarrow &
 \pi^* \mc{K}_Y 
\otimes \OO_{\PP(C \oplus \mathbf{1})}(-1)
& \longrightarrow & 0
\end{array}
\end{equation}
\end{enumerate}
\end{lem}

\begin{proof}
Item $(1)$ is very easy and is just the fact that an \emph{elementary
transform down} along a Cartier divisor gives rise to a locally free
sheaf.  To see Item $(2)$ and Item $(3)$, observe that the restriction
to $\PP(C \oplus \mathbf{1})$ of $(\widetilde{\mc{G}})^\vee
\rightarrow (\mc{G}_M)^\vee$ fits into an exact sequence
\begin{equation}
\begin{array}{lllllll}
0 & \longrightarrow & \textit{Tor}_1^{\OO_M}(\OO_{\PP(C \oplus
  \mathbf{1})}, \pi^*\mc{K}^\vee_Y) & \longrightarrow &
  (\widetilde{\mc{G}})^\vee|_{\PP(C \oplus \mathbf{1})} \\
 & \longrightarrow 
  &  \pi^*( \mc{G} \otimes_{\OO_B} \OO_Y )^\vee & \longrightarrow &
  \pi^* \mc{K}^\vee_Y & \longrightarrow & 0.
\end{array}
\end{equation}
Of course $\textit{Tor}_1^{\OO_M}(\OO_{\PP(C \oplus \mathbf{1})},
\OO_{\PP(C \oplus \mathbf{1})} )$ is just $\OO_M(-\PP(C \oplus
\mathbf{1}))|_{ \PP(C \oplus \mathbf{1}) }$, i.e. $\OO_{\PP(C \oplus
\mathbf{1})}(1)$.  Therefore the left-most term in the exact sequence
above is just $\mc{K}^\vee_Y \otimes \OO_{ \PP(C \oplus \mathbf{1})
}(1)$.  Dualizing this sequence gives Item $(2)$ and Item $(3)$.
\end{proof}

\medskip\noindent 
\begin{notat} \label{notat-gamma} 
The restriction of $\widetilde{\phi}$ to $\PP(C \oplus \mathbf{1})$
induces a morphism of locally free sheaves
\begin{equation}
\widetilde{\mc{G}}|_{\PP(C\oplus
  \mathbf{1})}/\pi^*(\mc{G}\otimes_{\OO_B} \OO_Y) \rightarrow
  \pi^*(\mc{F} \otimes_{\OO_B} \OO_Y)/\pi^*(\mc{G} \otimes_{\OO_B}
  \OO_Y).
\end{equation}
Up to canonical isomorphism, this is the same as a morphism
\begin{equation}
\gamma: \pi^*\mc{K}_Y \otimes \OO_{\PP(C \oplus \mathbf{1})}(-1)
\rightarrow \pi^* \mc{E}_Y.
\end{equation}
And the cokernel of $\widetilde{\phi}$ on $\PP(C \oplus \mathbf{1} )$
equals the cokernel of $\gamma$.
\end{notat}

\medskip\noindent
Our next goal is to show that the map $\gamma$ essentially is the same
as the map $\theta_{\phi,Y}$.

\medskip\noindent
The canonical inclusion $\widetilde{\mc{G}}\otimes \OO_M(-\PP(C \oplus
\mathbf{1}) ) \rightarrow \widetilde{\mc{G}}$ factors through the
kernel of $\widetilde{\mc{G}} \rightarrow \pi^*\mc{K}_Y \otimes
\OO_{\PP(C \oplus \mathbf{1})}(-1)$.  So there is an induced inclusion
$\widetilde{\mc{G}} \otimes \OO_M(-\PP(C \oplus \mathbf{1}))
\rightarrow \mc{G}_M$ whose cokernel is $\pi^*\lt( \mc{G}
\otimes_{\OO_B} \OO_Y/ \mc{K}_Y \rt)$.  In particular, we have a
commutative diagram:
\begin{equation} \label{eqn-cd1}
\begin{CD}
\widetilde{\mc{G}} \otimes \OO_M(-\PP(C \oplus \mathbf{1}) ) @>>>
\mc{G}_M \\
@V \widetilde{\phi} \otimes \text{Id} VV @VV \phi_M V \\
\mc{F}_M \otimes \OO_M(-\PP(C \oplus \mathbf{1}) ) @>>> \mc{F}_M
\end{CD}
\end{equation}

\medskip\noindent
\begin{lem} \label{lem-thetagamma} 
There is a commutative diagram of coherent sheaves:
\begin{equation} 
\begin{CD}
\pi^*\mc{K}_Y @> \pi^*\theta_{\phi,Y} >> \pi^*\mc{E}_Y \otimes
\pi^*(\mc{J}/\mc{J}^2) \\
@V \text{Id} VV @VV \text{Id} \otimes \beta_1 V \\
\pi^*\mc{K}_Y @> \gamma \otimes \text{Id} >> \pi^* \mc{E}_Y \otimes
\OO_{\PP(C \oplus \mathbf{1})}(1)
\end{CD}
\end{equation}
\end{lem}

\medskip\noindent
\begin{proof}
To ease notation in this proof, denote $\PP = \PP(C \oplus
\mathbf{1})$.  Consider the commutative diagram with exact rows
analogous to Equation(~\ref{eqn-snake1}) whose rows are
\begin{equation} \label{eqn-snake2a}
\begin{array}{ccccc}
0 & \longrightarrow & \mc{G}_M \otimes_{\OO_M} \OO_M(-\PP)/ 
\OO_M(-2\PP) & \longrightarrow & \mc{G}_M \otimes_{\OO_M}
\OO_M/ \OO_M(-2\PP) \\
 & \longrightarrow & \mc{G}_M
\otimes_{\OO_M} \OO_{\PP} & \longrightarrow & 0 
\end{array}
\end{equation} 

\begin{equation} \label{eqn-snake2b}
\begin{array}{ccccc}
0 & \longrightarrow & \mc{F}_M \otimes_{\OO_M} \OO_M(-\PP)/ 
\OO_M(-2\PP) & \longrightarrow & \mc{F}_M \otimes_{\OO_M}
\OO_M/ \OO_M(-2\PP) \\
 & \longrightarrow & \mc{F}_M
\otimes_{\OO_M} \OO_{\PP} & \longrightarrow & 0 
\end{array}
\end{equation} 
Associated to this commutative diagram, the snake lemma produces a
connecting map 
$\theta_{\phi_M,\PP}:\pi^* \mc{K}_Y \rightarrow \pi^*\mc{E}_Y \otimes
\OO_{\PP}(1)$.  

\medskip\noindent
Observe that the ideal sheaf $u^{-1}\text{pr}_1^{-1}(\mc{J})$ is
contained in the ideal sheaf $\OO_M(-\PP)$ (moreover when we divide by
the defining equation of $\PP$, the residual ideal sheaf is the ideal
sheaf of the closed immersion $\iota:Y \times \PP^1 \rightarrow M$).
Therefore there is a map from the pullback by $\text{pr}_1 \circ u$ of
the commutative diagram in Equation(~\ref{eqn-snake1}) to the
commutative diagram above.  In particular we have a commutative
diagram of connecting maps
\begin{equation}
\begin{CD}
\pi^* \mc{K}_Y @> \pi^* \theta_{\phi, Y} >> \pi^* \mc{E}_Y \otimes
\pi^*\lt( \mc{J}/\mc{J}^2 \rt) \\
@V \text{Id} VV @VV \text{Id} \otimes \beta V \\
\pi^* \mc{K}_Y @> \theta_{\phi_M, \PP} >> \pi^* \mc{E}_Y \otimes
\OO_{\PP}(1)
\end{CD}
\end{equation}
But now consider the map $\widetilde{\mc{G}}(-\PP) \rightarrow
\mc{G}_M \otimes_{\OO_M} \OO_M/\OO_M(-2\PP)$ constructed above.  The
image of $\widetilde{\mc{G}}(-\PP)$ in the quotient $\mc{G}_M
\otimes_{\OO_M} \OO_{\PP}$ is precisely $\pi^* \mc{K}_Y$.  Moreover
there is the map $\widetilde{\phi} \otimes \text{Id}:
\widetilde{\mc{G}}(-\PP) \rightarrow \mc{F}_M(-\PP)$ and the diagram
in Equation (~\ref{eqn-cd1}) commutes.  Therefore we can use
$\widetilde{\mc{G}}(-\PP)$ to compute the connecting map
$\theta_{\phi_M,\PP}$.  But the construction of $\gamma$ was by
precisely the same construction, i.e.  the connecting map
$\theta_{\phi_M,\PP}$ equals $\gamma\otimes \text{Id}$.  This proves
the lemma.
\end{proof}

\section{The stack of multiple covers of lines} \label{sec-covers}
\medskip\noindent
From this point on, we assume that our field $K$ is algebraically
closed of characteristic $0$.  Let $V$ denote a $K$-vector space of
dimension $n+1$ so that the projective space $\PP(V)$ is isomorphic to
$\PP^n$.  In the next sections, we will apply the analysis above to
$\Kbm{0,0}{\PP(V),e}$, the Kontsevich moduli stack parametrizing
stable maps from unmarked, genus $0$ curves to $\PP(V)$ of degree $e$.
For more details about this stack, see ~\cite{FP}.  The goal is to
prove that for a positive integer $d$ with $d+e \leq n$, for a general
hypersurface $X\subset \PP(V)$ of degree $d$, the closed substack
$\Kbm{0,0}{X,e} \subset \Kbm{0,0}{\PP(V),e}$ has only canonical
singularities.

\medskip\noindent
What does it mean to say that a pair of Deligne-Mumford stacks is
canonical (resp. log canonical, etc.)?  For a pair $(B,Y)$, one can
compute the log discrepancy $a(E;B,Y)$ \'etale locally on $B$, i.e. if
$(f_i:B_i \rightarrow B)$ is an \'etale cover, then $a(E;B,Y) =
\text{min}(a(f_i^* E; B_i, f_i^{-1}Y) | \text{center}(f_i^* E) \neq
\emptyset)$.  There is a standard way of extending any \'etale local
notion for schemes to Deligne-Mumford stacks: the Deligne-Mumford
stack has an \'etale local cover by schemes and the log discrepancies
are defined using this cover by the formula above.

\medskip\noindent
Consider the Kontsevich moduli stack $\Kbm{0,r}{\PP(V),e}$.  Let
$p:\mc{C} \rightarrow \Kbm{0,r}{\PP(V),e}$ denote the universal curve,
and let $f:\mc{C} \rightarrow \PP(V)$ denote the universal map.  For
each integer $d> 0$, on $\Kbm{0,r}{\PP(V),e}$ we have a locally free
sheaf $\mc{P}_d$ of rank $ed+1$ defined by $\mc{P}_d = p_* f^*
\OO_{\PP(V)}(d)$.  The fact that all the higher direct image sheaves
vanish and that $\mc{P}_d$ is locally free of rank $ed+1$ follows from
standard facts about stable maps together with cohomology and base
change.  We are actually only interested in the case that $r=0$, but
we mention that this holds for arbitrary $r$.  There is a canonical
\emph{evaluation morphism} of locally free sheaves on
$\Kbm{0,0}{\PP(V),e}$, $\phi_d^\dagger: H^0(\PP(V),\OO_{\PP(V)}(d))
\otimes_\kappa \OO \rightarrow \mc{P}_d$.

\medskip\noindent
\begin{defn} \label{defn-GtoF} 
For each $e \geq 1$ and $d\geq 0$, we define $\mc{G}_d$ to be the dual
of $\mc{P}_d$, we define $\mc{F}_d$ to be the trivial locally free
sheaf $H^0(\PP(V),\OO_{\PP(V)}(d))^\vee \otimes_\kappa \OO$, and we
define the \emph{co-evaluation morphism} for degree $d$, to be
\begin{equation}
\phi_d: \mc{G} \rightarrow \mc{F}
\end{equation}
which is the adjoint of the evaluation morphism above.  We define
$\mc{E}_d$ to be the cokernel of $\phi_d$.
We define $(\pi_d:C_d \rightarrow \Kbm{0,0}{\PP(V),e}, \alpha_d: \pi^*
\mc{E}_d \rightarrow \mc{Q}_d)$ to be the projective Abelian cone
parametrizing rank $1$ locally free quotients of $\mc{E}_d$.
\end{defn}

\medskip\noindent
When there is no risk of confusion, we will drop the subscripts.  We
will apply our techniques to analyze the singularities of $C$.  We
denote by $(\rho:C' \rightarrow \Kbm{0,0}{\PP(V),e},\beta:
H^0(\PP(V),\OO_{\PP(V)}(d))^\vee \otimes_\kappa \OO \rightarrow
\mc{Q}')$ the projective bundle parametrizing rank $1$ locally free
quotients of the trivial locally free sheaf
$H^0(\PP(V),\OO_{\PP(V)}(d)) \otimes_\kappa \OO$.  But $C'$ is just
the same as the product $\PP H^0(\PP(V),\OO_{\PP(V)}(d)) \times
\Kbm{0,0}{\PP(V),e}$.  As in Section~\ref{sec-discrep}, define $h:C
\rightarrow C'$ to be the tautological closed immersion.  Our interest
in $C$ is the following easy result.  The proof is left to the reader
(but also see ~\cite[Lemma 4.5]{HRS2}).

\medskip\noindent
\begin{lem} \label{lem-PE} 
The Deligne-Mumford stack $C$ parametrizes pairs $([X],[f:C\rightarrow
X])$ where $[X]\in \PP H^0(\PP(V),\OO_{\PP(V)}(d))$ is a hypersurface
of degree $d$ in $\PP(V)$, and where $f:C\rightarrow X$ is a
Kontsevich stable map of genus $0$ and degree $e$ to $X$.  In
particular, for each $[X] \in \PP H^0(\PP(V),\OO_{\PP(V)}(d))$, the
fiber of $C$ over $[X]$ is canonically identified with the Kontsevich
moduli stack $\Kbm{0,0}{X,e}$.
\end{lem}

\medskip\noindent
\begin{notat} \label{notat-Grass} 
Denote by $\mathbb{G} = \mathbb{G}(2,V)$ the Grassmannian variety over
$\kappa$ parametrizing rank $2$ locally free quotients of $V^\vee$,
i.e. parametrizing $2$-dimensional linear subspaces of $V$.  Let
$V^\vee \otimes_\kappa \OO_{\mathbb{G}} \rightarrow S^\vee$ denote the
universal quotient, so that the adjoint $S \hookrightarrow V
\otimes_\kappa \OO_{\mathbb{G}}$ is the universal rank $2$ subbundle.
Denote the quotient of the universal subbundle by $V \otimes_\kappa
\OO_{\mathbb{G}} \rightarrow T$.  Observe that $T^\vee \hookrightarrow
V^\vee \otimes_\kappa \OO_{\mathbb{G}}$ is simply the annihilator of
$S$.  For each $d \geq 0$, there is an induced filtration on
$\text{Sym}^d(V^\vee) \otimes_\kappa \OO_{\mathbb{G}}$
\begin{equation}
\begin{array}{l}
\text{Sym}^d(V^\vee) \otimes_\kappa \OO_{\mathbb{G}} = F^0 \supset F^1
\supset \dots \supset F^d \\
F^i = \text{Sym}^i(T^\vee) \cdot \text{Sym}^{d-i}(V^\vee), \ \ \
F^i/F^{i+1} \cong \text{Sym}^i(T^\vee) \otimes_{\OO_{\mathbb{G}}}
\text{Sym}^{d-i}(S^\vee)
\end{array}
\end{equation}
This filtration is the same as the filtration by order of vanishing
along $S$.
\end{notat}

\medskip\noindent
\begin{defn} \label{defn-covers} 
There is an induced morphism $\PP S \rightarrow \PP(V)$ identifying
$\mathbb{G}$ with the Hilbert scheme of lines in $\PP(V)$.  The
\emph{stack of multiple covers of lines} $Y$ is defined to be the
closed substack of $\mathbb{G} \times \Kbm{0,0}{\PP(V),e}$
parametrizing pairs $([L],[f:C \rightarrow L])$ where $[L]\in
\mathbb{G}$ is a line in $\PP(V)$ and $f:C\rightarrow L$ is a
Kontsevich stable map of genus $0$ and degree $e$.
\end{defn}

\medskip\noindent 
There are several equivalent definitions.  Of course the projection
$\text{pr}_{\mathbb{G}}: Y \rightarrow \mathbb{G}$ is Zariski locally
isomorphic to the product $\mathbb{G} \times \Kbm{0,0}{\PP^1,e}$.  An
easy observation is that the projection $Y \rightarrow
\Kbm{0,0}{\PP(V),e}$ is a closed immersion.

\medskip\noindent
\begin{lem} \label{lem-equivdefn} 
The $1$-morphism $Y \rightarrow \Kbm{0,0}{\PP(V),e}$ is representable
by closed immersions and the image is the rank $2$ locus
$(\Kbm{0,0}{\PP(V),e})_2$ for $\phi_1:\mc{G}_1 \rightarrow \mc{F}_1$.
Moreover, for each $\phi_d$ the reduced substack of the rank $(d+1)$
locus $(\Kbm{0,0}{\PP(V),e})_{d+1}$ equals the image of $Y$.
\end{lem}

\medskip\noindent
\begin{proof}
It is clear that $Y \rightarrow \Kbm{0,0}{\PP(V),e}$ is injective on
geometric points; after all for a pair $([L],[f:C \rightarrow L])$, we
have that $L = f(C)$, so the line $[L]$ is uniquely determined by $f:C
\rightarrow \PP(V)$.  Moreover for a stable map $f:C\rightarrow X$,
the rank of $\phi_1$ restricted to the residue field of $[f]$ is at
least as big as $H^0(f(C),\OO_{\PP(V)}(1)|_{f(C)})$.  For a pure
$1$-dimensional subscheme of $\PP(V)$, this dimension is always at
least $2$, and equals $2$ only if $f(C)$ is a line.  Therefore every
the geometric points of $(\Kbm{0,0}{\PP(V),e})_2$ equal the geometric
points of $Y$.  The same sort of argument shows that for every $d$,
the rank $(d+1)$ locus of $\phi_d$ equals $Y$ on the level of sets of
geometric points.

\medskip\noindent
Moreover, since $(\Kbm{0,0}{\PP(V),e})_1$ is empty, on
$(\Kbm{0,0}{\PP(V),e})_2$ the quotient $\mc{E}$ of $\phi_1$ is a
locally free sheaf of rank $2$.  Then $\mc{E}$ is a quotient of
$H^0(\PP(V),\OO_{\PP(V)}(1))^\vee$ which induces a morphism from
$(\Kbm{0,0}{\PP(V),e})_2$ to $\mathbf{G}$, from which it is easy to
construct an inverse to $Y \rightarrow (\Kbm{0,0}{\PP(V),e})_2$.
\end{proof}

\medskip\noindent
Next we analyze the restriction of $\phi_d$ to $Y$.  It is simpler to
phrase the results for the adjoint $\phi_d^\dagger$, but by Item $(1)$
of Lemma~\ref{lem-conn}, they are both equivalent.

\begin{lem} \label{lem-phi1} 
The kernel of $\phi_1^\dagger|_Y$ equals $\text{pr}_{\mathbb{G}}^*
T^\vee \subset V^\vee \otimes_\kappa \OO_Y$.  The cokernel of
$\phi_1^\dagger|_Y$ is a locally free sheaf $\mc{R}$ of rank $e-1$.
And the induced connecting map $\theta'_{\phi_1^\dagger,Y}:
\textit{Hom}_{\OO_Y} ( \mc{R}, \text{pr}_{\mathbb{G}}^* T^\vee )
\rightarrow \mc{J}/\mc{J}^2$ is an isomorphism of $\OO_Y$-modules.
\end{lem}

\medskip\noindent
\begin{proof}
The first two assertions follow from the proof of
Lemma~\ref{lem-equivdefn}; the details are left to the reader.  The
third assertion can probably be proved directly, but it also follows
from the deformation theory of Kontsevich stable maps developed in
\cite{B} and \cite{BF} (see also ~\cite[Sec. 3]{HS2}).  Since $Y$ is
smooth and since $\Kbm{0,0}{\PP(V),e}$ is smooth, the conormal sheaf
$\mc{J}/\mc{J}^2$ is a locally free $\OO_Y$-module.  We begin by
computing the dual of this locally free sheaf.

\medskip\noindent
As above, let $\pi:\mc{C} \rightarrow Y$ be the universal curve.  Let
$g:\mc{C} \rightarrow \PP(S)$ be the universal map (compatible with
projection to $\mathbb{G}$), and let $\text{pr}_{\PP(V)}:\PP(S)
\rightarrow \PP(V)$ be the obvious projection so that $f =
\text{pr}_{\PP(V)} \circ g$ is the universal map from $\mc{C}$ to
$\PP(V)$.  There is a perfect complex of amplitude $[-1,0]$ on
$\mc{C}$, denoted $L_f$, such that the object $(\mathbb{R}\pi_*
L_f^\vee)[1]$ in the derived category of $Y$ is quasi-isomorphic to
the restriction of the tangent bundle of $\Kbm{0,0}{\PP(V),e}$.
Similarly, there is a perfect complex of amplitude $[-1,0]$, denoted
$L_g$, such that the object $(\mathbb{R}\pi_* L_g^\vee)[1]$ in the
derived category of $Y$ is quasi-isomorphic to the vertical tangent
bundle of the morphism $\text{pr}_{\mathbb{G}}: Y \rightarrow
\mathbb{G}$.  These complexes are as follows:
\begin{equation}
\begin{CD}
& -1 & & 0 \\
L_f: & f^* \Omega^1_{\PP(V)} @> (df)^\dagger >> \Omega^1_{\mc{C}/Y} \\
L_g: & g^* \Omega^1_{\PP(S)/\mathbb{G}} @> (dg)^\dagger >>
\Omega^1_{\mc{C}/Y}
\end{CD}
\end{equation}
Of course the derivative of the morphism $\text{pr}_{\PP(V)}:\PP(S)
\rightarrow \PP(V)$ induces a surjective sheaf map from $f^*
\Omega^1_{\PP(V)}$ to $g^* \Omega^1_{ \PP(S)/\mathbb{G} }$ whose
kernel is just $g^*\OO_{\PP(S)}(-1)\otimes \text{pr}_{\mathbb{G}}^*
T^\vee$.  So we get a distinguished triangle
\begin{equation}
g^*\OO_{\PP(S)}(-1)\otimes \text{pr}_{\mathbb{G}}^* T^\vee [1]
\rightarrow L_f \rightarrow L_g
\end{equation}

\medskip\noindent
Applying the derived functor $\mathbb{R} \textit{Hom}_{ \OO_{\mc{C}} }
( * ,\OO_{\mc{C}} )$ to the distinguished triangle above produces a
distinguished triangle
\begin{equation}
L_g^\vee \rightarrow L_f^\vee \rightarrow \textit{Hom}_{\OO_{\mc{C}}} \lt(
\text{pr}_{\mathbb{G}}^* T^\vee, g^*\OO_{\PP(S)}(1) \rt) [-1]
\end{equation}
Finally we apply $\mathbb{R}\pi_*$ to this distinguished triangle to
get a distinguished triangle
\begin{equation}
(\mathbb{R}\pi_* L_g^\vee)[1] \rightarrow (\mathbb{R}\pi_*
  L_f^\vee)[1] \rightarrow \textit{Hom}_{\OO_Y} \lt(
  \text{pr}_{\mathbb{G}}^* T^\vee, \mc{P}_1 \rt)[0]
\end{equation}
(obviously we are skipping some details which are left to the reader).
So the derivative map from the vertical tangent bundle of
$\text{pr}_{\mathbb{G}}: Y \rightarrow \mathbb{G}$ to the restriction
of the tangent bundle of $\Kbm{0,0}{\PP(V),e}$ has cokernel isomorphic
to $\textit{Hom}_{\OO_Y} \lt( \text{pr}_{\mathbb{G}}^* T^\vee,
\mc{P}_1 \rt)$.

\medskip\noindent
The map $\phi_1^\dagger|_Y$ has image $\text{pr}_{\mathbb{G}}^*
S^\vee$.  Therefore inside of this cokernel we have this subsheaf
$\textit{Hom}_{\OO_Y} \lt( \text{pr}_{\mathbb{G}}^* T^\vee,
\text{pr}_{\mathbb{G}}^* S^\vee \rt)$.  This is just the pullback of
the tangent bundle of $\mathbb{G}$.  The cokernel of this subsheaf is
the normal bundle of $Y \rightarrow \Kbm{0,0}{\PP(V),e}$.  And, since
$\mc{R}$ is $\mc{P}_1/\text{pr}_{\mathbb{G}}^* S^\vee$ by definition,
we have that this cokernel is canonically isomorphic to
$\textit{Hom}_{\OO_Y} \lt(\text{pr}_{\mathbb{G}}^* T^\vee, \mc{R}
\rt)$.  Therefore the dual sheaf, $\mc{J}/\mc{J}^2$, is canonically
isomorphic to $\textit{Hom}_{\OO_Y} \lt( \mc{R},
\text{pr}_{\mathbb{G}}^* T^\vee \rt)$.  Obviously this is very close
to what we need to prove.

\medskip\noindent
By the canonical isomorphism above, the induced connecting map
$\theta'_{\phi_1^\dagger, Y}$ is now an endomorphism of the locally
free sheaf $\textit{Hom}_{\OO_Y} \lt( \mc{R}, \text{pr}_{\mathbb{G}}^*
T^\vee \rt)$.  An endomorphism of a locally free sheaf is an
isomorphism iff the determinant of the endomorphism is invertible.
Since $Y$ is a proper, smooth, connected Deligne-Mumford stack, the
global sections of $\OO_Y$ are just the constants.  So to prove that
the determinant is invertible, it suffices to prove that it is
somewhere nonzero.

\medskip\noindent
Now we are reduced to a simple, slightly tedious computation in local
coordinates.  Choose homogeneous coordinates $Y_0,\dots,Y_n$ on
$\PP(V)$, i.e. $Y_0,\dots,Y_n$ is an ordered basis for $V^\vee =
H^0(\PP(V),\OO_{\PP(V)}(1))$.  Choose homogeneous coordinates
$X_0,X_1$ on $\PP^1$.  Let $\AAA$ be the affine space associated to
the \emph{dual vector space} $W$ of linear transformations
\begin{equation}
W^\vee := \text{Hom}_{\kappa}\lt( 
H^0(\PP^1,\OO_{\PP^1}(e-2)), \text{span}\{Y_2,\dots,Y_n\} \rt).
\end{equation}
A basis for the vector space $W^\vee$, i.e. for the vector space of
linear forms on $\AAA$, is given by the tensors $(X_0^i
X_1^{e-2-i})^\vee \otimes Y_j$ for $j=2,\dots,n$ and $i=0,\dots,e-2$.
Define $F:\PP^1 \times \AAA \rightarrow \PP(V)$ to be the morphism
with $F^* \OO_{\PP(V)}(1) = \text{pr}_{\PP^1}^* \OO_{\PP^1}(e)$ and
where the pullback of homogeneous coordinates is defined by
\begin{equation}
\begin{array}{lll}
F^* Y_0 & = & \text{pr}_{\PP^1}^*X_0^e, \\
F^* Y_1 & = & \text{pr}_{\PP^1}^*X_1^e, \\
F^* Y_j & = & \sum_{i=0}^{e-2} \text{pr}_{\PP^1}^*
(X_0^{i+1}X_1^{e-1-i}) \cdot 
\text{pr}_{\AAA}^*( (X_0^i X_1^{e-2-i})^\vee \otimes Y_j), \
j=2,\dots,n
\end{array}
\end{equation}

\medskip\noindent
The morphism $F$ is a family of stable maps of degree $e$ over $\AAA$
and defines a $1$-morphism $\zeta:\AAA \rightarrow
\Kbm{0,0}{\PP(V),e}$.  The pullback by $\zeta$ of $\mc{P}_1$ is simply
the trivial vector bundle $H^0(\PP^1,\OO_{\PP^1}(e))$ on $\AAA$, and
the pullback of $\phi_1^\dagger$ is simply the map
\begin{equation}
\begin{array}{lll}
\phi_1^\dagger(1\otimes Y_0) & = & 1\otimes X_0^e, \\
\phi_1^\dagger(1\otimes Y_1) & = & 1\otimes X_1^e, \\
\phi_1^\dagger(1\otimes Y_j) & = & \sum_{i=0}^{e-2} \lt( (X_0^i
X_1^{e-2-i})^\vee \otimes Y_j \rt) \otimes X_0^{i+1} X_1^{e-1-i}, \
j=2,\dots, n
\end{array}
\end{equation}
It follows that the rank $2$ locus is the origin $0\in \AAA$.  The
inverse image ideal sheaf $\zeta^{-1}\mc{J}$ is just the ideal of the
origin, i.e. the ideal with generators $(X_0^i X_1^{e-2-i})^\vee
\otimes Y_j$ for $i=0,\dots,e$ and $j=2,\dots,n$.  The kernel of
$\phi_1^\dagger|_0$ is $\text{span}\{ Y_2,\dots, Y_n \}$.  The image
of $\phi_1^\dagger|_0$ is $\text{span} \{ X_0^e,X_1^e \}$, and the
cokernel is $X_0X_1 \cdot H^0(\PP^1,\OO_{\PP^1}(e-2))$.  So the
pullback of $\textit{Hom}_{\OO_Y}(\mc{E}_Y,\mc{K}_Y)$ is the vector
space
$$
\text{Hom}_\kappa \lt( H^0(\PP^1,\OO_{\PP^1}(e-2)),
\text{span}\{Y_2,\dots,Y_n \} \rt).
$$
Chasing through the snake diagram associated to $F$, the pullback of
the map $\theta'_{\phi_1^\dagger,Y}: \textit{Hom}_{\OO_Y}(\mc{E}_Y,
\mc{K}_Y) \rightarrow \mc{J}/\mc{J}^2$ is the map
\begin{equation}
(X_0^i X_1^{e-2-i})^\vee \otimes Y_j \mapsto (X_0^i X_1^{e-2-i})^\vee
  \otimes Y_j
\end{equation}
i.e. it is the identity map.  This proves that
$\theta_{\phi_1^\dagger,Y}$ is an isomorphism when restricted to the
image of $0\in \AAA$.  As mentioned above, this suffices to prove that
$\theta_{\phi_1^\dagger,Y}$ is everywhere an isomorphism.  In fact,
via the canonical isomorphism of $\mc{J}/\mc{J}^2$ with
$\textit{Hom}_{\OO_Y}(\mc{E}_Y,\mc{K}_Y)$ induced by the deformation
theory computation, the endomorphism $\theta_{\phi_1^\dagger,Y}$ is
the identity map.
\end{proof}

\begin{prop} \label{prop-phid} 
Using the isomorphisms from Lemma~\ref{lem-phi1}, we have the
following.
\begin{enumerate}
\item
For each $d\geq 1$, the kernel of $\phi_d^\dagger|_Y$ is
$\text{pr}_{\mathbb{G}}^* F^1 \subset \text{Sym}^d(V^\vee)
\otimes_\kappa \OO_Y$.
\item 
The cokernel of $\phi_d^\dagger|_Y$ is canonically isomorphic to
$\mc{R} \otimes_{\OO_Y} \text{pr}_{\mathbb{G}}^*
\text{Sym}^{d-1}(S^\vee)$.
\item
The induced connecting map
\begin{equation}
\theta'_{\phi_d^\dagger,Y}: \textit{Hom}_{\OO_Y} \lt( 
\mc{R} \otimes_{\OO_Y} \text{pr}_{\mathbb{G}}^*
\text{Sym}^{d-1}(S^\vee) ,
\text{pr}_{\mathbb{G}}^* F^1 \rt) 
\rightarrow \mc{J}/\mc{J}^2
\end{equation}
is the zero map on the subsheaf of homomorphisms from the domain to
$\text{pr}_{\mathbb{G}}^* F^2$.  
\item
Identify $F^1/F^2$ with $T^\vee \otimes_{\OO_{\mathbb{G}}}
\text{Sym}^{d-1}(S^\vee)$, and identify $\mc{J}/\mc{J}^2$ with
$\textit{Hom}_{\OO_Y}\lt( \text{pr}_{\mathbb{G}}^* T^\vee, \mc{R}
\rt)$ using $\theta'_{\phi_1^\dagger,Y}$.  The following map,
``induced'' by the induced connecting map:
\begin{equation}
\begin{array}{r}
\theta_{\phi_d^\dagger,Y}'': \textit{Hom}_{\OO_Y} \lt(
\mc{R} \otimes_{\OO_Y} \text{pr}_{\mathbb{G}}^*
\text{Sym}^{d-1}(S^\vee) ,
\text{pr}_{\mathbb{G}}^* T^\vee \otimes_{\OO_Y} \text{pr}_{\mathbb{G}}^*
\text{Sym}^{d-1}(S^\vee)
\rt) \\
\rightarrow \textit{Hom}_{\OO_Y} \lt( \mc{R}
\text{pr}_{\mathbb{G}}^* T^\vee \rt)
\end{array}
\end{equation}
is equal to the ``obvious'' map obtained by contracting the
$\text{Sym}^{d-1}(S^\vee)$ factors.
\end{enumerate}
\end{prop}

\begin{proof}
The sheaf $\mc{P}_d$ is defined to be $\pi_* f^* \OO_{\PP(V)}(d)$.
Consider the fiber product $\PP(S) \times_{\mathbb{G}} Y$.  There is
an induced map $(g,\pi):\mc{C} \rightarrow \PP(S) \times_{\mathbb{G}}
Y$, and $f^* \OO_{\PP(V)}(d)$ equals $(g,\pi)^* \text{pr}_{\PP(S)}^*
\OO_{\PP(S)}(d)$.  As $\pi = \text{pr}_{Y} \circ (g,\pi)$, we have the
identity
\begin{equation}
\mc{P}_d = (\text{pr}_{Y})_* (g,\pi)_* (g,\pi)^*
\text{pr}_{\PP(S)}^* \OO_{\PP(S)}(d)
\end{equation}
There is a canonical sheaf map $\text{pr}_{\PP(S)}^* \OO_{\PP(S)}(d)
\rightarrow (g,\pi)_* (g,\pi)^* \text{pr}_{\PP(S)}^* \OO_{\PP(S)}(d)$.
Because $(g,\pi)$ is proper and surjective, this sheaf map is
injective.  Since pushforward is left exact, we have an injective
sheaf map from $(\text{pr}_{Y})_* \text{pr}_{\PP(S)}^*
\OO_{\PP(S)}(d)$ to $\mc{P}_d$.  But of course the first sheaf is just
$\text{pr}_{\mathbb{G}}^* \text{Sym}^d(S^\vee)$.  And the injective
sheaf map $\text{Sym}^d(S^\vee) \rightarrow \mc{P}_d$ is the image of
$\phi_d^\dagger|_Y$.  Therefore the kernel of $\phi_d^\dagger|_Y$ is
the pullback of the kernel of $\text{Sym}^d(V^\vee) \otimes_\kappa
\OO_{\mathbb{G}} \rightarrow \text{Sym}^d(S^\vee)$, i.e. the first
filtered subsheaf $F^1$ of $\text{Sym}^d(V^\vee) \otimes_\kappa
\OO_{\mathbb{G}}$.  This proves that the kernel of $\phi_d^\dagger|_Y$
equals $\text{pr}_{\mathbb{G}}^* F^1$.  So Item $(1)$ is verified.

\medskip\noindent
On $\PP(S)$ we have a short exact sequence
\begin{equation}
0 \rightarrow \text{pr}_{\mathbb{G}}^*
\mathbb{S}_{(d-1,1)}(S^\vee)\otimes \OO_{\PP(S)}(-1)
\rightarrow \text{pr}_{\mathbb{G}}^* \text{Sym}^{d-1}(S^\vee)
\rightarrow \OO_{\PP(S)}(d-1) \rightarrow 0
\end{equation}
where $\mathbb{S}_{(d-1,1)}$ is the \emph{Schur functor} as defined in
~\cite[Sec. 6.1]{FH}.  Twist this sequence by $\OO_{\PP(S)}(1)$ and
pullback to $\mc{C}$ by $(g,\pi)$ to get a short exact sequence on
$\mc{C}$.  When we pushforward by $\pi$, there is a long exact
sequence of higher direct image sheaves.  Since $\pi:\mc{C}
\rightarrow Y$ is a flat family of at-worst-nodal curves of genus $0$,
we have that $R^1\pi_* \OO_{\mc{C}}$ is zero.  Therefore the long
exact sequence is really the following short exact sequence:
\begin{equation}
0 \rightarrow \text{pr}_{\mathbb{G}}^* \mathbb{S}_{(d-1,1)}(S^\vee)
\rightarrow \text{pr}_{\mathbb{G}}^* \text{Sym}^{d-1}(S^\vee) \otimes
\mc{P}_1 \rightarrow \mc{P}_d \rightarrow 0
\end{equation}
And of course we have a commutative diagram with exact rows
\begin{equation}
\begin{array}{ccccccccc}
0 & \rightarrow & \text{pr}_{\mathbb{G}}^*
\mathbb{S}_{(d-1,1)}(S^\vee) & \rightarrow &
\text{pr}_{\mathbb{G}}^* \text{Sym}^{d-1}(S^\vee) \otimes
\text{pr}_{\mathbb{G}}^* S^\vee & \rightarrow &
\text{pr}_{\mathbb{G}}^* \text{Sym}^d(S^\vee) & \rightarrow & 0 \\
 & & \text{Id} \downarrow & & \downarrow \text{Id}\otimes
\phi_1^\dagger & & \downarrow \phi_d^\dagger \\
0 & \rightarrow & \text{pr}_{\mathbb{G}}^* \mathbb{S}_{(d-1,1)}(S^\vee)
& \rightarrow & \text{pr}_{\mathbb{G}}^* \text{Sym}^{d-1}(S^\vee) \otimes
\mc{P}_1 & \rightarrow & \mc{P}_d & \rightarrow & 0
\end{array}
\end{equation}
Applying the snake lemma to this diagram, we get an isomorphism of the
cokernel of $\phi_d^\dagger$ with $\mc{R}\otimes
\text{pr}_{\mathbb{G}}^* \text{Sym}^{d-1}(S^\vee)$.  This proves Item
$(2)$.

\medskip\noindent
In case $d=1$, Item $(3)$ follow from Lemma~\ref{lem-phi1}.  Therefore
we may suppose that $d > 1$.  To prove Item $(3)$, we use the fact
that we have a commutative diagram of sheaves on
$\Kbm{0,0}{\PP(V),e}$:
\begin{equation}
\begin{CD}
V^\vee \otimes_\kappa \text{Sym}^{d-1}(V^\vee) \otimes_\kappa \OO @>
\phi_1^\dagger \otimes \text{Id} >> \mc{P}_1 \otimes_\kappa
\text{Sym}^{d-1}(V^\vee) \\
@V \psi_{\mc{G}} VV @VV \psi_{\mc{F}} V \\
\text{Sym}^d(V^\vee) \otimes_\kappa \OO @> \phi_d^\dagger >> \mc{P}_d
\end{CD}
\end{equation}
Associated to $\phi_1^\dagger \otimes \text{Id}$ is the induced
connecting map
\begin{equation}
\theta'_{\phi_1^\dagger \otimes \text{Id},Y}:
\textit{Hom}_{\OO_Y} \lt( \mc{R} \otimes_\kappa \text{Sym}^{d-1}(V^\vee),
\text{pr}_{\mathbb{G}}^* T^\vee \otimes_\kappa
\text{Sym}^{d-1}(V^\vee) \rt)  \rightarrow \mc{J}/\mc{J}^2.
\end{equation}
Of course this is obtained from $\theta'_{\phi_1^\dagger,Y}$ by
contracting the $\text{Sym}^{d-1}(V^\vee)$ factors.
Define $\theta''_{\phi_1^\dagger \otimes \text{Id}, Y}$ to be the
restriction of $\theta'_{\phi_1^\dagger \otimes \text{Id}, Y}$ to the
subsheaf 
\begin{equation}
\textit{Hom}_{\OO_Y} \lt( \mc{R} \otimes_{\OO_Y}
\text{pr}_{\mathbb{G}}^* \text{Sym}^{d-1}(S^\vee),
\text{pr}_{\mathbb{G}}^* T^\vee \otimes_\kappa
\text{Sym}^{d-1}(V^\vee) \rt).
\end{equation}
By Item $(4)$ of Lemma~\ref{lem-conn}, we have a commutative diagram
of induced connecting maps
\begin{equation}
\begin{CD}
\textit{Hom}_{\OO_Y} \lt( \mc{R} \otimes_{\OO_Y} \text{Sym}^{d-1}(S^\vee),
\text{pr}_{\mathbb{G}}^* T^\vee \otimes_\kappa
\text{Sym}^{d-1}(V^\vee) \rt) @> \theta'_{\phi_1^\dagger\otimes
  \text{Id},Y} >> \mc{J}/\mc{J}^2 \\
@V\psi VV @VV \text{Id} V \\
\textit{Hom}_{\OO_Y} \lt( \mc{R} \otimes_{\OO_Y}
\text{pr}_{\mathbb{G}}^* \text{Sym}^{d-1}(S^\vee),
\text{pr}_{\mathbb{G}}^* F^1 \rt) @> \theta'_{\phi_d^\dagger, Y} >>
\mc{J}/\mc{J}^2 
\end{CD}
\end{equation}
Since $\theta'_{\psi_1^\dagger\otimes \text{Id},Y}$ is obtained by
contracting the $\text{Sym}^{d-1}(V^\vee)$ factors,
in particular the kernel of
$\theta_{\phi_1^\dagger \otimes \text{Id}, Y}$ contains the subsheaf
\begin{equation}
\textit{Hom}_{\OO_Y} \lt( \mc{R} \otimes_{\OO_Y}
\text{pr}_{\mathbb{G}}^* \text{Sym}^{d-1}(S^\vee),
\text{pr}_{\mathbb{G}}^* T^\vee \otimes_{\OO_Y} F^1 \rt).
\end{equation}
Therefore the kernel of $\theta'_{\phi_d^\dagger,Y}$ contains the
image under $\psi$ of this subsheaf.  But the image under $\psi$ is
just
\begin{equation}
\textit{Hom}_{\OO_Y} \lt( \mc{R} \otimes_{\OO_Y}
\text{pr}_{\mathbb{G}}^* \text{Sym}^{d-1}(S^\vee),
\text{pr}_{\mathbb{G}}^*F^2 \rt).
\end{equation}
Forming the quotient by this sheaf, we have another commutative
diagram
\begin{equation}
\begin{CD}
\textit{Hom}_{\OO_Y} \lt( \mc{R} \otimes_{\OO_Y} 
\text{Sym}^{d-1}(S^\vee),
\text{pr}_{\mathbb{G}}^* T^\vee \otimes_\kappa
\text{Sym}^{d-1}(V^\vee)  \rt) 
@> \theta'_{\phi_1^\dagger\otimes
  \text{Id},Y} >> \mc{J}/\mc{J}^2 \\
@V\psi'' VV @VV \text{Id} V \\
\textit{Hom}_{\OO_Y} \lt( \mc{R} \otimes_{\OO_Y}
\text{pr}_{\mathbb{G}}^* \text{Sym}^{d-1}(S^\vee),
\text{pr}_{\mathbb{G}}^* T^\vee \otimes_{\OO_Y}
\text{pr}_{\mathbb{G}}^* \text{Sym}^{d-1}(S^\vee) 
\rt) @> \theta''_{\phi_d^\dagger, Y} >>
\mc{J}/\mc{J}^2 
\end{CD}
\end{equation}
Since $\psi''$ is surjective, $\theta''_{\phi_d^\dagger,Y}$ is the
unique morphism making the above diagram commute.  But using the fact
that $\theta'_{\phi_1^\dagger \otimes \text{Id}, Y}$ is obtained from
contracting the $\text{Sym}^{d-1}(V^\vee)$ factors, it is clear that
the diagram also commutes when we replace
$\theta''_{\phi_d^\dagger,Y}$ by the map obtained from
$\theta'_{\phi_1^\dagger,Y}$ by contracting the
$\text{Sym}^{d-1}(S^\vee)$ factors.  Therefore
$\theta''_{\phi_d^\dagger,Y}$ equals the map obtained from
$\theta'_{\phi_1^\dagger,Y}$ by contracting the
$\text{Sym}^{d-1}(S^\vee)$ factors.  This finshes the proof.
\end{proof}

\medskip\noindent
\begin{cor} \label{cor-phid} 
Consider the restriction $\phi_d|_Y:\mc{G}_d|_Y \rightarrow
\mc{F}_d|_Y$.
\begin{enumerate}
\item 
The cokernel $\mc{E}_Y$ is canonically isomorphic to
$\text{pr}_{\mathbb{G}}^* (F^1)^\vee$.
\item 
The kernel $\mc{K}_Y$ is canonically isomorphic to $\lt( \mc{R}
\otimes_{\OO_Y} \text{pr}_{\mathbb{G}}^* \text{Sym}^{d-1}(S^\vee)
\rt)^\vee$.
\item 
The induced connecting map $\theta'_{\phi_d,Y}$ is canonically
isomorphic to the induced connecting map $\theta'_{\phi_d^\dagger,Y}$.
\end{enumerate}
\end{cor}

\medskip\noindent
\begin{proof}
This follows from Item $(1)$ of Lemma~\ref{lem-conn} and
Proposition~\ref{prop-phid}.
\end{proof}

\section{Proof of the main theorem} \label{sec-proof} 

\medskip\noindent 
Let an integer $e \geq 1$ be fixed.  To simplify notation, in this
section we denote $B=\Kbm{0,0}{\PP(V),e}$.  By Lemma~\ref{lem-PE}, for
each $d \geq 1$ the projective Abelian cone $\pi:C \rightarrow B$ of
the coherent sheaf $\mc{E}_d = \text{Coker}(\phi_d)$ is the
Deligne-Mumford stack parametrizing pairs $([X],[f:C \rightarrow X])$
where $X\subset \PP(V)$ is a hypersurface of degree $d$ and
$[f:C\rightarrow X]$ is a point in $\Kbm{0,0}{X,e}$.  In this section
we describe the singularities of the cone $C_d$.  We reiterate that we
are working in characteristic $0$ (though not necessarily over an
algebraically closed field).

\medskip\noindent
We begin with the simplest case $e=1$.  The next two results are
already known, in fact in arbitrary characteristic
~\cite[Thm. V.4.3]{K}.  We only prove that part which we will use.

\medskip\noindent
\begin{prop} \label{prop-e=1} 
If $e=1$, then $B$ is a smooth projective scheme and for all $d\geq 1$
the morphism $\pi_d:C_d \rightarrow B$ is a projective bundle of the
expected dimension.  In particular $C_d$ is a geometrically
irreducible, smooth scheme of the expected dimension.
\end{prop}

\begin{proof}
Since $e=1$, $Y$ equals $B$, which is simply the Grassmannian
$\mathbb{G}$.  Thus $\mc{R}$ is the zero sheaf.  And by
Corollary~\ref{cor-phid}, the cokernel $\mc{E}_d$ is locally free of
the expected dimension.  Therefore $\pi:C \rightarrow B$ is a
projective bundle.
\end{proof}

\begin{thm} \label{thm-e=1} \cite[V.4.3]{K}
If $e=1$ and if $d > 2n-3$, then the projection morphism $h_d:C_d
\rightarrow \PP H^0\lt( \PP(V), \OO_{\PP(V)}(d) \rt)$ is not
surjective, i.e. the general fiber is empty.  If $d \leq 2n-3$, then
the projection morphism $h_d$ is surjective, and the general fiber is
a smooth scheme of the expected dimension $2n-d-3$.  Moreover if $d <
2n-3$ and $(d,n) \neq (2,3)$, then the general fiber is geometrically
connected.
\end{thm}

\medskip\noindent
\begin{proof}
For the full proof, the reader should consult ~\cite{K}.  We will only
prove that part of the theorem which we shall use, namely that the
general fiber is nonempty and smooth for $d \leq n-1$ and that the
general fiber is geometrically connected for $d \leq n-2$.

\medskip\noindent
Suppose that $d\leq n-1$.  Since $C_d$ is irreducible of the expected
dimension, $h_d$ is surjective iff the general fiber has the expected
dimension.  To prove the general fiber has the expected dimension, it
suffices to find one pair $([X],[L])$ consisting of a hypersurface
$X\subset \PP(V)$ and a line $L\subset X$ such that $X$ is smooth
along $L$ and such that $H^1(L,N_{L/X})$ is zero.  Then the Zariski
tangent space of the fiber, i.e. $H^0(L,N_{L/X})$, has the expected
dimension which proves that on a nonempty (hence dense) open subset of
$C$, $h_d$ has the expected fiber dimension.  Choose homogeneous
coordinates $Y_0,Y_1,Y_2,\dots,Y_n$ on $\PP(V)$.  Define $L$ to be the
vanishing locus of $Y_2,\dots,Y_n$.  Define $X\subset \PP(V)$ to be
the hypersurface with defining equation
\begin{equation}
 F:= \sum_{j=2}^{d+1} Y_0^{d+1-j}Y_1^{j-2} Y_j
\end{equation}
At every point of $L$, either the partial derivative $F_2$ is nonzero
or the partial derivative $F_{d+1}$ is nonzero.  Therefore $X$ is
smooth along $L$.  Moreover, by the usual exact sequence
\begin{equation}
0 \rightarrow N_{L/X} \rightarrow \OO_L(1)^{n-1}
\xrightarrow{F_2,\dots,F_{d+1},0,\dots,0} \OO_L(d) \rightarrow 0
\end{equation}
we conclude that $N_{L/X} \cong \OO_L^{d-1}\oplus \OO_L(1)^{n-d-1}$.
Therefore $H^1(L,N_{L/X})$ is zero which proves that $h_d$ is
surjective and the generic fiber has the expected dimension.  Since
$C_d$ is smooth and since $\PP H^0(\PP(V),\OO_{\PP(V)}(d) )$ is
smooth, it follows by generic smoothness that the generic fiber of
$h_d$ is everywhere smooth.

\medskip\noindent
Next we prove that the generic fiber of $h_d$ is connected when $d
\leq n-2$.  In this case, by the same sort of dimension computation as
above, we conclude that for a general hypersurface $X\subset \PP(V)$,
every irreducible component of $\Kbm{0,1}{X,1}$ surjects to $X$ under
the evaluation morphism $\text{ev}:\Kbm{0,1}{X,1} \rightarrow X$.
Therefore, to prove that $\Kbm{0,1}{X,1}$ is connected, and thus that
$\Kbm{0,0}{X,1}$ is connected, it suffices to prove that every fiber
of $\text{ev}$ is connected.  Now for a given $p\in X$, the set of
lines in $\PP(V)$ which contain $p$ is canonically isomorphic to the
projective space $\PP(V/L_p)$ where $L_p \subset V$ is the
one-dimensional vector subspace corresponding to $p$.  Choose
homogeneous coordinates so that $p=[1:0:\dots:0]$.  Then the Taylor
expansion of the defining equation $F$ of $X$ about the point $p$ has
the form
\begin{equation}
F = Y_0^{d-1}F_1(Y_1,\dots, Y_n) + Y_0^{d-2}F_2(Y_1,\dots,Y_n) +
\dots + Y_0F_{d-1}(Y_1,\dots, Y_n) + F_d(Y_1,\dots, Y_n)
\end{equation}
where each $F_i$ is a homogeneous polynomial of degree $i$ in
$Y_1,\dots,Y_n$.  A line passing through $p$ with parametric equation
$[(1-t):tY_1:tY_2:\dots:tY_n]$ is contained in $X$ iff the equations
$F_1(Y_1,\dots,Y_n), \dots , F_d(Y_1,\dots, Y_n)$ are all zero.  It
should be observed that some of the homogeneous equations $F_i$ may be
identically zero.  Nonetheless, the common zero locus of $F_1,\dots,
F_d$ in $\PP(V/L_p)$ is an intersection of at most $d$ hypersurfaces
in a projective space of dimension $n-1$.  In a projective space of
dimension $n-1$, an intersection of at most $n-2$ hypersurfaces is
always connected.  Since $d \leq n-2$, we conclude that the common
zero locus of $F_1,\dots, F_d$ is connected.  Therefore every fiber of
$\text{ev}$ is connected, which proves that the fiber $\Kbm{0,0}{X,1}
= h_d^{-1}([X])$ is connected.
\end{proof}

\medskip\noindent
Now suppose that $e>1$.  Consider the closed immersion $Y
\hookrightarrow B$.  Because the conormal sheaf of $Y$ is locally
free, we denote the normal cone $C$ by the letter $N$ (also to avoid
confusion with the cone $C_d$ which we are studying).  Associated to
the closed immersion $Y \rightarrow B$, we have the deformation to the
normal cone $(\varrho:M \rightarrow \PP^1, \iota:Y \times \PP^1
\rightarrow M, B_Y \rightarrow M, \PP(N\oplus \mathbf{1}) \rightarrow
M)$ as in Definition~\ref{defn-dnc}.  Define $M^o = M - B_Y$.  In
particular, the intersection of $M^o$ with $\PP(N \oplus \mathbf{1})$
is just $N$.

\medskip\noindent
Let $\widetilde{\phi}_d:\widetilde{\mc{G}} \rightarrow \mc{F}_M$ be
the elementary transform up of the pullback of $\phi_d$ as described
in Lemma~\ref{lem-elemtrans}.  Let $\widetilde{\mc{E}}_d$ be the
cokernel of $\widetilde{\phi}_d|_{M^o}$ and let
$\widetilde{\pi}_d:\widetilde{C}_d \rightarrow M^o$ be the projective
Abelian cone parametrizing rank $1$ locally free quotients of
$\widetilde{\mc{E}}_d$.  Observe that over the open subset
$\varrho^{-1}(\AAA^1) = B \times \AAA^1$, $\widetilde{\mc{E}}_d$ is
simply $\text{pr}_B^* \mc{E}_d$ and $\widetilde{C}_d$ is simply $C_d
\times \AAA^1$.

\medskip\noindent
\begin{lem} \label{lem-restrict} 
Let $e \geq 2$.  If $d+e \leq n$, then the fiber product
$\widetilde{C}_d \times_{M^o} N$ is an integral, normal, local
complete intersection scheme of the expected dimension which has
canonical singularities.
\end{lem}

\medskip\noindent
\begin{proof}
By Lemma~\ref{lem-phi1}, the normal bundle $N\rightarrow Y$ is the
vector bundle associated to the locally free sheaf
$\textit{Hom}_{\OO_Y}(\mc{R}^\vee, \text{pr}_{\mathbb{G}}^* T)$,
i.e. $N = M^{(0)}(Y,\mc{R}^\vee, \text{pr}_{\mathbb{G}}^* T)$ in the
notation of Section~\ref{sec-deter}.  For ease of notation, denote by
$\mc{A}$ the locally free sheaf $\text{pr}_{\mathbb{G}}^* \lt(
\text{Sym}^{d-1}(S^\vee) \rt)^\vee$.  This is a locally free sheaf of
rank $d$ on $Y$.

\medskip\noindent
By Corollary~\ref{cor-phid}, the kernel of $\phi_d|_Y$ is the locally
free sheaf
\begin{equation}
\mc{K}_Y = 
\mc{A} \otimes_{\OO_Y} \lt(\mc{R}^\vee \rt) 
\end{equation}
and the cokernel $\mc{E}_Y$ of $\phi_d|_Y$ fits into a short exact
sequence of locally free sheaves
\begin{equation}
0 \rightarrow \mc{A} \otimes_{\OO_Y} \lt( \text{pr}_{\mathbb{G}}^*
T\rt) 
\rightarrow \mc{E}_Y \rightarrow \text{pr}_{\mathbb{G}}^*(F^2)^\vee
\rightarrow 0.
\end{equation}
Denote the first sheaf in this sequence by $\mc{E}'_Y$ and denote the
third sheaf in this sequence by $\mc{A}'$.  By
Notation~\ref{notat-gamma}, the cokernel of $\widetilde{\phi}_d$ on
the Cartier divisor $N\subset M^o$ equals the cokernel of the sheaf
map $\gamma: \pi_Y^* \mc{K}_Y \rightarrow \pi_Y^* \mc{E}_Y$.  Here we
are using the fact that $\OO_{\PP(N\oplus \mathbf{1})}(1)$ is
canonically trivialized on $N$, so we will treat all twists by this
sheaf as twists by the trivial sheaf.  By Lemma~\ref{lem-thetagamma},
the map $\gamma$ is the unique map induced by $\theta_{\phi_d,Y}$.  By
Item $(1)$ of Lemma~\ref{lem-conn}, this means that $\gamma$ is the
transpose of the unique map induced by $\theta_{\phi_d^\dagger,Y}$.
By Item $(3)$ of Lemma~\ref{prop-phid}, the map
$\theta_{\phi_d^\dagger,Y}$ is the map induced by the universal
homomorphism from $\text{pr}_{\mathbb{G}}^*(T^\vee)$ to $\mc{R}$.

\medskip\noindent
Putting all the pieces together, we conclude two things.  First, the
image of $\gamma$ is actually the image of a sheaf map
\begin{equation}
\gamma': \mc{K}_Y \otimes_{\OO_Y} \OO_N \rightarrow \mc{E}'_Y
\otimes_{\OO_Y} \OO_N
\end{equation}
so that we have a short exact sequence (which is split Zariski locally
over Y)
\begin{equation}
0 \rightarrow \text{Coker}(\gamma') \rightarrow \text{Coker}(\gamma)
\rightarrow \text{pr}_Y^* \mc{A}' \rightarrow 0.
\end{equation}
And second, if we denote by $\psi: \mc{R}^\vee \otimes_{\OO_Y} \OO_N
\rightarrow \text{pr}_{\mathbb{G}}^* T \otimes_{\OO_Y} \OO_N$ the
universal sheaf homomorphism on $N$, then the sheaf map $\gamma'$ is
just $\text{Id} \otimes \psi$:
\begin{equation}
\gamma' = \text{Id} \otimes \psi:  
\mc{A} \otimes_{\OO_Y} \lt(\mc{R}^\vee \rt) \otimes_{\OO_Y} \OO_N
\rightarrow 
\mc{A} \otimes_{\OO_Y} \lt( \text{pr}_{\mathbb{G}}^* T\rt)
\otimes_{\OO_Y} \OO_N. 
\end{equation}

\medskip\noindent
Now the rank of $\mc{R}^\vee$ is $e-1$, the rank of
$\text{pr}_{\mathbb{G}}^* T$ is $n-1$ and the rank of $\mc{A}$ is $d$.
And, by assumption, $d \leq (n-1) - (e-1)$.  But this means that,
Zariski locally over $Y$, we are in the situation of
Proposition~\ref{prop-cone}.  So we conclude that $\widetilde{C}_d
\times_{M^o} N$ is a normal, integral, local complete intersection
scheme of the expected dimension which has canonical singularities.
\end{proof}

\medskip\noindent
Now let $W \subset B$ be the maximal open substack over which $C_d$ is
an integral, normal scheme of the expected dimension and with only
canonical singularities. Observe that $W\times \AAA^1 \subset B\times
\AAA^1$ is the maximal open subset over which $\widetilde{C}_d
\times_{M^o} (B\times \AAA^1)$ is an integral, normal scheme of the
expected dimension and which only canonical singularities.  Define $W'
\subset M^o$ to be the maximal open substack over which
$\widetilde{C}_d$ is an integral, normal scheme of the expected
dimension and with only canonical singularities

\medskip\noindent
\begin{lem} \label{lem-VcontY} 
If $e \geq 2$ and $d+e \leq n$, then the open substack $W \subset B$
contains the closed substack $Y \subset B$.
\end{lem}

\medskip\noindent
\begin{proof}
By Lemma~\ref{lem-restrict}, by Corollary~\ref{cor-discrep}, by
Proposition~\ref{prop-relcan}, and by Item $(3)$ of
Corollary~\ref{cor-spec} applied to $N \subset M^o$, the open substack
$W'$ contains $N$.

\medskip\noindent
There is one slight hiccup in checking the hypothesis of Item $(3)$ of
Corollary~\ref{cor-spec}.  In case $n-e \geq 2$, we have that
$\text{codim}_N(N_{g-1}) = (f-(g-1))(g-(g-1)) = f-g+1$ equals $(n-1) -
(e-1) + 1 = n-e+1 \geq 3$, so that Lemma~\ref{lem-D} proves that the
hypothesis of Item $(3)$ of Corollary~\ref{cor-spec} is satisfied.

\medskip\noindent
The one remaining case is when $d=1$ and $e=n-1$, and in this case we
give an ad hoc argument.  In this case the morphism
$\widetilde{\phi}_1$ restricted to $N$ is just the universal sheaf map
$\psi$.  In this special case $\mc{D}_{\phi_N}$ sits inside $N
\times_Y \PP(\text{pr}_{\mathbb{G}}^* T)$.  The projection
$\mc{D}_{\phi_N} \rightarrow \PP(\text{pr}_{\mathbb{G}}^* T)$ is a
Zariski locally trivial bundle.  Given a closed point $p\in Y$ and a
one-dimensional subspace $L\subset \text{pr}_{\mathbb{G}}^* T|_p$, the
fiber over this point, considered as a subvariety of $N =
M^{(0)}(Y,\mc{R}^\vee, \text{pr}_{\mathbb{G}}^* T)$, equals the cone
whose vertex set is $Hom_{\kappa(p)}(\mc{R}^\vee|_p, L)$ and whose
quotient by the vertex set is the set of linear maps in
$Hom_{\kappa(p)}(\mc{R}^\vee|_p, \text{pr}_{\mathbb{G}}^* T|_p/L)$
which are not isomorphisms.  Observe that this second vector space is
essentially just the vector space of square $(e-1)\times(e-1)$
matrices.  In the special case $e=2$, the cone is just a linear space
and so it is smooth.  In case $e\geq 3$ the vertex set has codimension
$(e-1)^2 -1 \geq 3$ in the fiber of $\mc{D}_{\phi_N}$, and the
singular locus of the quotient has codimension $4-1 = 3$:
$Hom_{\kappa(p)}(\mc{R}^\vee|_p, \text{pr}_{\mathbb{G}}^*
T|_p/L)_{e-3}$ has codimension $4$.  Therefore the singular locus has
codimension $3$ in the fiber of $\mc{D}_{\phi_N}$.  So the fiber is
normal, which implies that $\mc{D}_{\phi_N}$ is normal.  Therefore
when $d=1$ and $e=n-1$, the hypothesis of Item $(3)$ of
Corollary~\ref{cor-spec} are again satisfied.

\medskip\noindent
Of course we have $W' \cap \varrho^{-1}(\AAA^1) = W \times \AAA^1$.
Let $p\in Y$ be any point and consider $\iota(\{p\} \times \PP^1)
\subset M$.  By construction of the deformation to the normal cone
from Definition~\ref{defn-dnc}, $\iota(p,\infty)$ is the point on the
zero section of $N\rightarrow Y$ over $p\in Y$.  In particular
$\iota(\{p\} \times \PP^1) \subset M^o$ and $\iota(p,\infty) \in N$.
Therefore $\iota(\{ p \} \times \PP^1)$ intersects $W'$.  So $\iota(\{
p \} \times \AAA^1)$ intersects $W' \cap \varrho^{-1}(\AAA^1)$, i.e
$\{ p \} \times \AAA^1$ intersects $W\times \AAA^1$.  Therefore $p\in
W$.  So we conclude that $Y \subset W$ as was to be proved.
\end{proof}

\begin{thm} \label{thm-main} 
If $e\geq 2$ and if $d+e \leq n$, then $C_d$ is an integral, normal,
local complete intersection stack of the expected dimension which has
canonical singularities.
\end{thm}

\begin{proof}
By Lemma~\ref{lem-VcontY}, the open substack $W$ contains $Y$.  Now
the automorphism group $\text{GL}(V)$ acts on $\PP(V)$ and thus on
$B$.  Moreover the sheaves $\mc{G}_d$ and $\mc{F}_d$ have natural
$\text{GL}(V)$-linearizations and the morphism $\phi_d$ is
$\text{GL}(V)$-equivariant.  Therefore $W$ is a
$\text{GL}(V)$-invariant open substack of $B$.  So to prove that $W =
B$, it suffices to prove that the closure of every
$\text{GL}(V)$-orbit intersects $Y$.

\medskip\noindent
Let $f:D \rightarrow \PP(V)$ be any stable map of genus $0$ and degree
$e$.  Choose a direct sum decomposition $V = V_2 \oplus V_{n-1}$ so
that $\PP(V_{n-1}) \subset \PP(V)$ is disjoint from $f(D)$.  Consider
the $\mathbb{G}_m$-action $m: \mathbb{G}_m \times \PP(V) \rightarrow
\PP(V)$ given by $t\cdot(v,v') = (v,tv')$ where $v\in V_2, v'\in
V_{n-1}$.  This defines an action of a subgroup scheme of
$\text{GL}(V)$, and acting on $[f:D \rightarrow \PP(V)]$ yields a
$1$-morphism $\zeta:\mathbb{G}_m \rightarrow B$.  The limit as
$t\rightarrow 0$ of this action is simply the stable map $g\circ f: D
\rightarrow \PP(V_2) \subset \PP(V)$ where $g:\PP(V) - \PP(V_{n-1})
\rightarrow \PP(V_2)$ is the projection map.  In particular, $g\circ
f$ is a multiple cover of the line $\PP(V_2)$.  Therefore the closure
of the $\text{GL}(V)$-orbit of $[f:D\rightarrow \PP(V)]$ intersects
$Y$ in the point $[g\circ f:D \rightarrow \PP(V)]$.  It follows that
$W$ is all of $B$, i.e. $C_d$ is an integral, normal, local complete
intersection stack of the expected dimension and with only canonical
singularities.
\end{proof}

\begin{cor} \label{cor-main} 
If $e \geq 2$ and if $d+e \leq n$, then for a general hypersurface
$X\subset \PP(V)$ of degree $d$, the Kontsevich moduli space
$\Kbm{0,0}{X,e}$ is an integral, normal, local complete intersection
stack of the expected dimension $(n+1-d)e+(n-3)$ with only canonical
singularities.
\end{cor}

\begin{proof}
By Theorem~\ref{thm-main}, $C_d$ is integral, normal, Gorenstein and
canonical.  Consider the projection $h_d:C_d \rightarrow \PP
H^0(\PP(V), \OO_{\PP(V)}(d) )$.  The pullback of the hyperplane linear
system gives a base-point-free linear system on $C_d$.  By repeated
application of ~\cite[Thm. 1.13]{Reid80} (see also
~\cite[Prop. 7.7]{KollarPairs}), we conclude that the general fiber of
$h_d$ is a reduced, normal, local complete intersection stack with
only canonical singularities.  The one issue that remains is
connectedness, i.e. it is a priori possible that $h_d:C_d \rightarrow
\PP H^0(\PP(V), \OO_{\PP(V)}(d) )$ has a nontrivial Stein
factorization.  But observe that by Corollary~\ref{thm-e=1} (in fact
by the part that we proved there), the restriction of $h_d$ to $Y$ is
surjective and has a trivial Stein factorization.  So $h_d|_Y$ yields
a section of the Stein factorization of $h_d$, which is irreducible
and finite over $\PP H^0(\PP(V), \OO_{\PP(V)}(d) )$.  It follows that
also the Stein factorization of $h_d$ is trivial, i.e. the fibers of
$h_d$ are connected.  So the general fiber of $h_d$ is an integral,
normal, local complete intersection stack of the expected dimension
with only canonical singularities.
\end{proof}

\section{Singularities of $\kbm{0,0}{\PP^n,e}$}
~\label{sec-RSBT} 

\medskip\noindent
In this section we use the Reid--Shepherd-Barron--Tai criterion to
prove that (with a very few exceptions) the coarse moduli spaces
$\kbm{0,0}{\PP^n,e}$ are terminal.  We use the same computations to
show that if one carries out the deformation to the normal cone
construction as in Section~\ref{sec-proof} over the coarse moduli
space, one obtains a family whose general fiber is the coarse moduli
space of $C_d$ and whose special fiber is a normal, $\QQ$-Gorenstein,
canonical variety.  It follows that the ``inversion of adjunction''
conjecture implies that $C_d$ is itself canonical, and therefore
$\kbm{0,0}{X,e}$ is canonical for $X$ general.

\medskip\noindent
Let $\Gamma$ be a finite cyclic group of order $r$ and let $\zeta \in
\text{Hom}_{\text{group}}(\Gamma, \mathbb{G}_m)$ be a generator for
the character group of $\Gamma$.  Let $M$ be a finite dimensional
$\Gamma$-representation over $k$ (we are working over a field of
characteristic $0$, but the following definition makes sense if the
characteristic is prime to $\#\Gamma$).  There is a direct sum
decomposition
\begin{equation}
M = \oplus_{i=0}^{r-1} L_{\zeta^i}^{\oplus a_i}
\end{equation}
where each $L_{\zeta^i}$ is the one-dimensional representation
corresponding to the character $\zeta^i$.  

\begin{defn} \label{defn-invt} 
The \emph{invariant of} $M$ \emph{with respect to} $\zeta$ (after
Reid--Shepherd-Barron--Tai) is
\begin{equation}
\alpha_\zeta(M) = \frac{1}{r} \sum_{i=0}^{r-1} ia_i.
\end{equation}
The \emph{invariant of} $M$ is $\alpha(M) = \min{\alpha_\zeta(M)}$ as
$\zeta$ varies over all generators of the character group.
\end{defn}

The relevance is the following:

\medskip\noindent
\begin{thm} [Reid--Shepherd-Barron--Tai criterion, ~\cite{YPGS}]
    \label{thm-RSBT} 
Let $Y$ be a smooth $k$-variety, let $G$ be a finite subgroup of the
group of $k$-automorphisms of $Y$, and suppose that $G$ acts without
quasi-reflections.  Then the quotient variety $X=Y//G$ is terminal
(resp. canonical) iff for every cyclic subgroup $\Gamma \subset G$ and
every closed point $x\in Y^\Gamma$, the invariant of the Zariski
tangent space to $Y$ at $x$ satisfies
\begin{equation}
\alpha(T_x Y) > 1, \ \ \text{resp.} \ \ \alpha(T_x Y) \geq 1.
\end{equation}
\end{thm}

\medskip\noindent
\begin{cor} \label{cor-RSBT} 
Let $\mc{X}$ be a Deligne-Mumford stack which is smooth over $k$.
Denote by $p:\mc{X} \rightarrow X$ the coarse moduli space, and
suppose that $p$ is an isomorphism on the complement of a closed
subset of codimension $\geq 2$.  Then $X$ is terminal
(resp. canonical) iff for every geometric point $x$ of $\mc{X}$, and
for every cyclic subgroup $\Gamma$ of the stabilizer group of $x$, the
invariant of the Zariski tangent space to $\mc{X}$ at $x$ satisfies
\begin{equation}
\alpha(T_x \mc{X}) > 1, \ \ \text{resp.} \ \ \alpha(T_x Y) \geq 1.
\end{equation}
\end{cor}

\medskip\noindent
\begin{proof}
This is just a rewording of Theorem~\ref{thm-RSBT} into the language
of Deligne-Mumford stacks.
\end{proof}

\medskip\noindent
\begin{cor} \label{cor-DMsmooth} 
Let $\mc{X}$ be a Deligne-Mumford stack which is smooth over $k$.  Let
$f:\mc{Y} \rightarrow \mc{X}$ be a smooth, representable $1$-morphism
of Deligne-Mumford stacks.  Suppose that the morphism $p:\mc{X}
\rightarrow X$ is an isomorphism away from codimension $2$ and that
$X$ is terminal (resp. canonical).  Then the map to the coarse moduli
space of $\mc{Y}$, say $q:\mc{Y} \rightarrow Y$, is an isomorphism
away from codimension $2$ and $Y$ is terminal (resp. canonical).
\end{cor}

\medskip\noindent
\begin{proof}
Let $U\subset \mc{X}$ denote the maximal open substack over which $p$
is an isomorphism.  Then $\mc{X} - U$ has codimension at least $2$ in
$\mc{X}$.  Since $f$ is smooth, in particular it is flat.  Therefore
$f^{-1}(\mc{X} -U)$ has codimension at least $2$ in $\mc{Y}$.  And
$f^{-1}(U)$ is a scheme because $f$ is representable.  Therefore $q$
is an isomorphism when restricted to $f^{-1}(U)$, which shows that $q$
is an isomorphism away from codimension $2$.

\medskip\noindent
Now we apply Corollary~\ref{cor-RSBT}.  Let $y$ be a geometric point
of $\mc{Y}$ and let $x=f(y)$.  Because $f$ is representable, the
homomorphism from the stabilizer group of $y$ to the stabilizer group
of $x$ is injective.  So a cyclic subgroup $\Gamma$ of the stabilizer
group of $y$ is also a cyclic subgroup of the stabilizer group of $x$.
By Corollary~\ref{cor-RSBT}, the invariant of $T_x \mc{X}$ as a
$\Gamma$-representation is greater than $1$ (resp. at least $1$).
Since $f$ is smooth, the differential $df:T_y \mc{Y} \rightarrow T_x
\mc{X}$ is surjective.  Therefore the invariant of $T_y \mc{Y}$ is
greater than or equal to the invariant of $T_x \mc{X}$.  Applying
Corollary~\ref{cor-RSBT} one more time, we conclude that $Y$ is
terminal (resp. canonical).
\end{proof}

\medskip\noindent
\begin{rmk} \label{rmk-DMsmooth} 
Unfortunately, for \emph{nice} representable $1$-morphisms $f:\mc{Y}
\rightarrow \mc{X}$ of smooth Deligne-Mumford stacks which are not
smooth, Corollary~\ref{cor-DMsmooth} often fails.  For instance, if
$\mc{Z} \subset \mc{X}$ is a Zariski closed substack which is smooth
and $f:\mc{Y} \rightarrow \mc{X}$ is the blowing up of $\mc{X}$ along
$\mc{Z}$, it can happen that $p:\mc{X} \rightarrow X$ is an
isomorphism away from codimension $2$, that $q:\mc{Y} \rightarrow Y$
is an isomorphism away from codimension $2$, that $X$ is terminal
(resp. canonical), but $Y$ is not terminal (resp. canonical).  For
instance, consider the action of the group of third roots of unity
$\mathbb{\mu}_3$ on affine $4$-space $\AAA^4$ by $\omega \cdot
(X_1,X_2,X_3,X_4) = (\omega X_1, \omega X_2, \omega X_3, \omega X_4)$.
Let $Z \subset \AAA^4$ be the variety associated to the invariant
ideal $\langle X_1, X_2, X_3 \rangle$.  Let $Y\rightarrow \AAA^4$
denote the blowing up along $Z$.  Then $f:[Y/\mathbb{\mu}_3]
\rightarrow [X/\mathbb{\mu}_3]$ is a $1$-morphism of smooth
Deligne-Mumford stacks satisfying the hypotheses above and
$X//\mathbb{\mu}_3$ is terminal.  But $Y//\mathbb{\mu}_3$ is not even
canonical.
\end{rmk}

\medskip\noindent
Let $\Gamma$ be a finite cyclic group of order $r$ and let $\Delta
\subset \Gamma$ be a subgroup of index $s$.  Let $\gamma:\Gamma
\rightarrow \mathbb{G}_m$ be a generator for the character group of
$\Gamma$.  The restriction of $\gamma$ to $\Delta$ is a generator for
the character group of $\Delta$.  The following lemma is a rewording
of the argument on pp. $33$--$34$ of ~\cite{HM}.

\medskip\noindent
\begin{lem} \label{lem-induced} 
Let $V$ be a finite-dimensional representation of $\Delta$ and let $V
\otimes_{k[\Delta]} k[\Gamma]$ be the induced $\Gamma$-representation.
The relation between the invariant of $V\otimes_{k[\Delta]} k[\Gamma]$
as a $\Gamma$-representation and the invariant of $V$ as a
$\Delta$-representation is:
\begin{equation}
\alpha_{\Gamma,\gamma}\lt( V\otimes_{k[\Delta]} k[\Gamma] \rt) =
\alpha_{\Delta,\gamma|_{\Delta}}(V) + \frac{s-1}{2}\text{dim}_k(V).
\end{equation}
\end{lem}

\medskip\noindent
\begin{proof}
Each side of the equation is additive in $V$, therefore we may reduce
to the case that $V$ is an irreducible representation, i.e. a
character $V=L_{\gamma|_{\Delta}^l}$ for some integer $l=0,\dots
\frac{r}{s}-1$.  Let $\phi$ be a generator for $\Gamma$, so that
$\phi^s$ is a generator for $\Delta$.  Let $\overline{\epsilon}$ be a
nonzero element of $V$.  For each integer $j=0,\dots,s-1$, denote
$m=-l-\frac{j\cdot r}{s}$ and define the element $\epsilon_j \in V
\otimes_{k[\Delta]} k[\Gamma]$ to be:
\begin{equation}
\epsilon_j = \sum_{i=0}^{s-1} \gamma^{m}(\phi^i)
\overline{\epsilon}\otimes \phi^i.
\end{equation}
It is trivial to compute that $\epsilon_j \cdot \phi =
\epsilon_j \cdot \gamma^{-m}(\phi)$.  So $\epsilon_j$ spans an
irreducible subrepresentation of $V \otimes{k[\Delta]} k[\Gamma]$
which is isomorphic to $L_{\gamma^{-m}}$.  

\medskip\noindent
This gives $s$ different irreducible subrepresentations of
$V\otimes_{k[\Delta]} k[\Gamma]$, which is also the dimension as a
$k$-vector space.  So we have an irreducible decomposition:
\begin{equation}
V\otimes_{k[\Delta]} k[\Gamma] \cong \bigoplus_{j=0}^{s-1}
L_{\gamma^{l + j\cdot\frac{r}{s}}}.
\end{equation}
It follows that the invariant of $V\otimes_{k[\Delta]} k[\Gamma]$ as a
$\Gamma$ representation is
\begin{equation}
\begin{array}{r}
\frac{1}{r} \lt[l + (l+\frac{r}{s}) + \dots + (l +
(s-1)\frac{r}{s})\rt] =
l/\lt({\frac{r}{s}}\rt) + \frac{s-1}{2} = \\
 \alpha_{\Delta,\gamma|_{\Delta}}(V)
+ \frac{s-1}{2} \text{dim}_k(V)
\end{array}
\end{equation}
which is what we needed to prove. 
\end{proof}

\medskip\noindent
Recall that $\Kbmo{0,0}{\PP^n,e} \subset \Kbm{0,0}{\PP^n,e}$ is the
open substack parametrizing stable maps with irreducible domain, and
$\kbmo{0,0}{\PP^n,e}$ is the coarse moduli space of
$\Kbmo{0,0}{\PP^n,e}$.

\begin{prop} \label{prop-RSBT1} 
Let $x$ be a geometric point of $\Kbmo{0,0}{\PP^n,e}$ and let $\Gamma$
be a subgroup of the stabilizer group of $x$.  Denote $r = \#\Gamma$.
The invariant of the Zariski tangent space to $\Kbmo{0,0}{\PP^n,e}$ at
$x$ equals
\begin{equation}
\alpha(T_x \Kbmo{0,0}{X,e}) = \frac{e(n+1)}{2} \lt( 1 - \frac{1}{r}
\rt) - 1.
\end{equation}
Except in the cases $(e,n) = (2,1)$ and $(e,n) = (2,2)$, the map
$p:\Kbmo{0,0}{\PP^n,e} \rightarrow \kbmo{0,0}{\PP^n,e}$ is an
isomorphism away from codimension $2$.  Disregarding $(e,n)=(2,1),
(2,2)$, $\kbmo{0,0}{\PP^n,e}$ is terminal except in the the cases
$(e,n) = (3,1), (2,3) $ and in these cases it is canonical.
\end{prop} 

\medskip\noindent
\begin{proof}
By the same $\text{GL}_{n+1}$-invariance argument as in the proof of
Theorem~\ref{thm-main}, it suffices to prove the result when $x$ is a
geometric point of $Y \cap \Kbmo{0,0}{\PP^n,e}$.  At such a point we
can decompose the Zariski tangent space as a direct sum of the Zariski
tangent space to $Y$ and the normal bundle to $Y$.  The Zariski
tangent space to $Y$ further decomposes as the direct sum of the
Zariski tangent space to $\mathbb{G}$ and the vertical tangent bundle
of $\text{pr}_{\mathbb{G}}:Y \rightarrow \mathbb{G}$.  And by
Lemma~\ref{lem-phi1}, the normal bundle to $Y$ is a direct sum of
$n-1$ copies of $\mc{R}^\vee$.  So we need to compute the invariants
(with respect to some $\zeta$) of the vertical tangent bundle of
$\text{pr}_{\mathbb{G}}$ and of $\mc{R}^\vee$.

\medskip\noindent
The vertical tangent bundle of $\text{pr}_{\mathbb{G}}$ is just the
same as the tangent bundle of $\Kbmo{0,0}{\PP^1,e}$, so suppose now
that $n=1$.  Let the geometric point $x$ parametrize a stable map $f:C
\rightarrow \PP^1$.  Choose a generator for $\Gamma$, which will be an
automorphism $\phi:C \rightarrow C$ such that $f\circ \phi =f$ and
such that $\phi^s = \text{Id}$ iff $r$ divides $s$.  It is easy to
show that, up to a choice of homogeneous coordinates, $\phi:C
\rightarrow C$ is just the isomorphism $[X_0:X_1] \mapsto [X_0:\xi
X_1]$ for some primitive $r^{\text{th}}$ root of unity.

\medskip\noindent
Denote by $g:C \rightarrow C_0$ the quotient of $C$ by $\phi$ and let
$h:C_0 \rightarrow \PP^1 $ be the unique morphism such that $f=h\circ
g$.  The Zariski tangent space to $\Kbmo{0,0}{\PP^1,e}$ is just the
vector space of global sections of the torsion sheaf
$f^*T_{\PP^1}/T_C$.  And this fits into an exact sequence:
\begin{equation}
\begin{CD} 
0 @>>> g^* T_{C_0}/T_C @>>> g^* h^* T_{\PP^1}/T_C @>>> g^*\lt(
h^*T_{\PP^1}/T_{C_0} \rt) @>>> 0
\end{CD}
\end{equation}
Now, as a representation of $\Gamma$, $g^* \lt( h^* T_{\PP^1}/T_{C_0}
\rt)$ is isomorphic to the tensor product $\lt( h^* T_{\PP^1}/T_{C_0}
\rt) \otimes_k k[\Gamma]$ where the first factor is a trivial
representation.  In particular, the invariant with respect to any
generator $\zeta$ is just
\begin{equation}
\alpha_\zeta(g^* ( h^* T_{\PP^1}/T_{C_0} ) ) = 2(\frac{e}{r} -1) \cdot
\lt( \frac{0}{r} + \frac{1}{r} + \dots + \frac{r-1}{r} \rt) =
(e-r)\lt( 1 - \frac{1}{r} \rt).
\end{equation}

\medskip\noindent
By direct computation, as a representation of $\Gamma$,
$g^*T_{C_0}/T_C$ is isomorphic to
\begin{equation}
g^* T_{C_0}/T_C \cong L_{\xi^0}^{\oplus 2} \oplus L_{\xi^1}^{\oplus 1}
\oplus L_{\xi^{r-1}}^{\oplus 1} \oplus \bigoplus_{i=1}^{r-2}
L_{\xi^i}^{\oplus 2}.
\end{equation}
It follows that the invariant with respect to any generator $\zeta$ is
just
\begin{equation}
\alpha_\zeta( g^* T_{C_0}/T_C) = 2 \lt(\frac{0}{r} + \frac{1}{r} +
\dots + \frac{r-1}{r} \rt) - \lt( \frac{i}{r} + \frac{r-i}{r} \rt) =
r\lt( 1 - \frac{1}{r} \rt) - 1.
\end{equation}
Here $i$ is the unique integer such that $\{ \xi^1, \xi^{r-1} \} = \{
\zeta^i, \zeta^{r-i} \}$.  Summing up, we have that with respect to
any generator $\zeta$,
\begin{equation}
\alpha_\zeta( T_x \Kbmo{0,0}{\PP^1,e} ) = e\lt( 1 - \frac{1}{r} \rt) -
1.
\end{equation}

\medskip\noindent
Next we compute the invariant of $\mc{R}|_x$.  Denote this vector
space by $\mc{R}_f$.  Now each of $g:C \rightarrow C_0$ and $h:C_0
\rightarrow \PP^1$ is also a stable map of a genus $0$ curve to
$\PP^1$.  So each of these also has a canonically associated vector
space $\mc{R}_g$, respectively $\mc{R}_h$.  As $\Gamma$
representations, $\mc{R}_g$ is just $k[\Gamma]/k{1}$ (by direct
computation), and $\mc{R}_h$ is a trivial representation of dimension
$\frac{e}{r} - 1$.  It is easy to see that the relationship between
these spaces is that $\mc{R}_f$ is isomorphic as a
$\Gamma$-representation to $\mc{R}_h \otimes_k \mc{R}_g \oplus
\mc{R}_g \oplus \mc{R}_h$.  Therefore the invariant of $\mc{R}_f$ is
\begin{equation}
\alpha_\zeta( \mc{R}_f ) = \frac{e}{2}\lt( 1 - \frac{1}{r} \rt).
\end{equation}
This is also the invariant of $\mc{R}_f^\vee$.

\medskip\noindent
As mentioned above, the normal bundle of $Y$ at $x$ is isomorphic as a
$\Gamma$-representation to a direct sum of $n-1$ copies of
$\mc{R}_f^\vee$.  And the vertical tangent space to
$\text{pr}_{\mc{G}}$ is just the same as the tangent space to
$\Kbm{0,0}{\PP^1,e}$.  Therefore we conclude that
\begin{equation}
\alpha_{\zeta}(T_x \Kbmo{0,0}{\PP^n,e}) = \frac{(n-1)e}{2} \lt( 1 -
\frac{1}{r} \rt) + e\lt( 1 - \frac{1}{r} \rt) - 1 = \frac{(n+1)e}{2}
\lt( 1 - \frac{1}{r} \rt) - 1.
\end{equation}

\medskip\noindent
For $(e,n) \neq (2,1), (2,2)$, the invariant is at least $1$, which
shows that the stabilizer group of $x$ acts without quasi-reflections
and the coarse moduli space has canonical singularities.  Moreover,
except in the extra cases $(e,n) = (3,1), (2,3)$, the invariant is
actually larger than $1$ which shows that the coarse moduli space has
terminal singularities.
\end{proof}

\medskip\noindent
\begin{rmk} \label{rmk-RSBT} 
In case $e=2, \ n=1$, every geometric point of $\Kbmo{0,0}{\PP^1,2}$
has nontrivial stabilizer.  In fact the coarse moduli space
$\kbmo{0,0}{\PP^1,2}$ is isomorphic to the complement of a smooth
plane conic in $\PP^2$ (via the \emph{branch morphism},
c.f. ~\cite{FaP}), and $p:\Kbmo{0,0}{\PP^1,2} \rightarrow
\kbmo{0,0}{\PP^1,2}$ is a $\ZZ/2\ZZ$-gerbe.  In case $e=2, \ n=2$, the
coarse moduli space $\kbmo{0,0}{\PP^2,2}$ is smooth and is isomorphic
to an open subset of the blowing up of $\PP^5$ along a Veronese
surface (the open subset is the complement of the proper transform of
the discriminant hypersurface).  In this case the morphism
$p:\Kbmo{0,0}{\PP^2,2} \rightarrow \kbmo{0,0}{\PP^2,2}$ is an
isomorphism on the complement of the exceptional divisor, and over the
exceptional divisor it is a $\ZZ/2\ZZ$-gerbe.
\end{rmk}

\medskip\noindent
\begin{lem} \label{lem-Gor} 
Let $(e,n)$ be a pair of positive integers other than $(2,1)$ and
$(2,2)$.  Then $\kbmo{0,0}{\PP^n,e}$ is Gorenstein unless both $e$ and
$n$ are even, in which case it is not Gorenstein.
\end{lem}

\medskip\noindent
\begin{proof}
Since $(e,n) \neq (2,1), (2,2)$, we know that $p:\Kbmo{0,0}{\PP^n,e}
\rightarrow \kbmo{0,0}{\PP^n,e}$ is an isomorphism away from
codimension $2$.  It follows from ~\cite[Prop. 5.75]{KM} that the
dualizing sheaf of $\kbmo{0,0}{\PP^n,e}$ is the pushforward by $p_*$
of the dualizing sheaf of $\Kbmo{0,0}{\PP^n,e}$.  Given a geometric
point $x$ of $\Kbmo{0,0}{\PP^n,e}$, the dualizing sheaf of
$\kbmo{0,0}{\PP^n,e}$ is invertible at $p(x)$ iff there exists a
section of the dualizing sheaf near $p(x)$ whose pullback to
$\Kbmo{0,0}{\PP^n,e}$ is non-zero at $x$.  Such a section corresponds
to a nonzero element in the one-dimensional vector space
$\text{det}(T_x \Kbmo{0,0}{\PP^n,e})^\vee$ which is invariant under
the action of the stabilizer group of $x$.  Therefore
$\kbmo{0,0}{\PP^n,e}$ is Gorenstein iff for every geometric point $x$
of $\Kbmo{0,0}{\PP^n,e}$ and for every cyclic subgroup $\Gamma$ of the
stabilizer group of $x$, the induced character $\text{det}(T_x
\Kbmo{0,0}{\PP^n,e})$ is trivial.

\medskip\noindent
As in the proof of Proposition~\ref{prop-RSBT1}, it suffices to
compute the character for geometric points $x$ of $Y\cap
\Kbmo{0,0}{\PP^n,e}$.  In the proof of Proposition~\ref{prop-RSBT1} we
computed all the relevant $\Gamma$-representations.  The character of
$g^*\lt( h^*T_{\PP^1}/T_{C_0} \rt)$ is just
$\text{det}(k[\Gamma])^{\otimes 2(\frac{e}{r} -1)}$ where
$2(\frac{e}{r} -1)$ is the dimension of $h^*T_{\PP^1}/T_{C_0}$.
Similarly, the character of $T_x \Kbmo{\PP^1,e}$ is
$\text{det}(k[\Gamma])^{\otimes 2}$ (the missing $L_{\xi^1}$ and
$L_{\xi^{r-1}}$ factors tensor to give a trivial character).  Finally,
the character of $\mc{R}_g$ is $\text{det}(k[\Gamma])$ and the
character of $\mc{R}_f \cong \lt( \mc{R}_h \otimes_k \mc{R}_g \rt)
\oplus \lt( \mc{R}_h \rt) \oplus \lt( \mc{R}_g \rt)$ is
$\text{det}(k[\Gamma])^{\otimes \frac{e}{r} }$.  Altogether, the
character of $T_x \Kbmo{0,0}{\PP^n,e}$ is
$\text{det}(k[\Gamma])^{\otimes (n+1)\frac{e}{r}}$.

\medskip\noindent
If $r$ is odd, then the character $\text{det}(k[\Gamma])$ is trivial
so that the character of $T_x \Kbmo{0,0}{\PP^n,e}$ is trivial.  But if
$r$ is even, the character $\text{det}(k[\Gamma])$ equals
$L_{\zeta^{\frac{r}{2}}}$ for any generator $\zeta$ of the character
group of $\Gamma$.  This is a nontrivial character whose square is
trivial.  For $r$ even, the character of $T_x \Kbmo{0,0}{\PP^n,e}$ is
nontrivial iff $\frac{e}{r}$ is odd and $n+1$ is odd.  Therefore if
$e$ is odd or $n$ is odd, then the character of $T_X
\Kbmo{0,0}{\PP^n,e}$ is trivial for every geometric point,
i.e. $\kbmo{0,0}{\PP^n,e}$ is Gorenstein.

\medskip\noindent
On the other hand, suppose that $n$ and $e$ are both even.  Then for
any line $L \subset \PP^n$ and any reduced degree $2$ divisor on $L$,
the cyclic cover $f:C \rightarrow L$ of degree $e$ branched over that
divisor gives a stable map of degree $e$ whose stabilizer group is
cyclic of order $r=e$.  Therefore the character of $T_x
\Kbmo{0,0}{\PP^n,e}$ is nontrivial, i.e. $\kbmo{0,0}{\PP^n,e}$ is not
Gorenstein.
\end{proof}

\medskip\noindent
\begin{prop} \label{prop-RSBT2} 
Let $(e,n)$ be a pair of positive integers different from $(2,1)$ and
$(2,2)$.  When $e \geq 3$ and $n=1$, the coarse moduli space
$\kbm{0,0}{\PP^1,e}$ is canonical.  When $(e,n)=(2,3)$, the coarse
moduli space $\kbm{0,0}{\PP^3,2}$ is canonical.  In all other cases,
the coarse moduli space $\kbm{0,0}{\PP^n,e}$ is terminal.
\end{prop}

\medskip\noindent
\begin{proof}
The proof follows closely the argument on pp. $33$--$34$ of
~\cite{HM}.  As in the proof of Proposition~\ref{prop-RSBT1}, to prove
that $\kbm{0,0}{\PP^n,e}$ is canonical (resp. terminal), it suffices
to check that for every geometric point $x$ of $Y$ and for every
cyclic subgroup $\Gamma$ of the stabilizer group of $x$, the invariant
of $T_x \Kbm{0,0}{\PP^n,e}$ is bigger than $1$ (resp. at least $1$).
Choose some line $\PP^1 \subset \PP^n$; throughout we will work with
geometric points $x$ of the closed substack $\Kbm{0,0}{\PP^1,e}
\subset \Kbm{0,0}{\PP^n,e}$.  We will prove that the invariant of $T_x
\Kbm{0,0}{\PP^n,e}$ is bigger than $1$ (resp. at least $1$) by
induction on the number of nodes of $C$.  We have already computed the
invariant when the domain $C$ of the stable map $f:C \rightarrow
\PP^1$ is irreducible, i.e. we have proved the proposition when the
number of nodes is zero.  Therefore suppose that the number of nodes
is positive.

\medskip\noindent
Let $\phi$ be a generator for $\Gamma$ and let $\{q,\phi q, \phi^2 q,
\dots, \phi^{s-1} q\}$ be an orbit of $\Gamma$ on $C$ such that each
$\phi^i q$ is a node.  Of course we have that $s$ divides $r$.  We use
the language of ~\cite{BM} regarding stable $A$-graphs (i.e. the dual
graph of $C$ labelled by the degree of $f$) $\tau$ and the
Behrend-Manin moduli stacks $\Kbm{}{\PP^1, \tau}$.  Let $\tau$ be the
dual graph of $f:C \rightarrow \PP^1$.  Let $E_0,E_1,\dots, E_{s-1}$
be the edges of $\tau$ corresponding to the nodes $q, \phi q, \dots,
\phi^{s-1} q$.  Let $\psi: \tau \rightarrow \sigma$ be the maximal
contraction of $\tau$ which does not contract any of the edges $E_0,
\dots, E_{s-1}$, i.e. $\sigma$ is the same as the dual graph of a
curve obtained by smoothing all the nodes of $C$ except $q,\dots,
\phi^{s-1}q$.  Then $f:C \rightarrow \PP^1$ gives a geometric point of
the Behrend-Manin moduli stack $\Kbm{}{\PP^1,\sigma}$, i.e. the moduli
space of stable maps to $\PP^1$ whose dual graph has a contraction to
$\sigma$.

\medskip\noindent
There is a canonical $1$-morphism $\Kbm{}{\PP^1,\sigma} \rightarrow
\Kbm{0,0}{\PP^1,e}$ which is unramified and whose normal sheaf is
locally free.  Therefore the tangent bundle of $\Kbm{}{\PP^1,\sigma}$
at $[f:C \rightarrow X]$ is a vector subspace of the tangent bundle of
$\Kbm{0,0}{\PP^1,e}$ at $[f:C \rightarrow X]$.  Moreover the cokernel,
i.e. the normal bundle, is precisely
\begin{equation}
N_{[f]} = \bigoplus_{i=0}^{s-1} T'_{\phi^i q} \otimes_k T''_{\phi^i q}
\end{equation}
where $T'_{\phi^i q}$ and $T''_{\phi^i q}$ are the tangent spaces of
the two branches of $C$ at $\phi^i q$ (there isn't any canonical
ordering of the two branches; the notation $T'$ and $T''$ is just for
convenience).  

\medskip\noindent
Now suppose there exists a nonzero section $\epsilon \in N_{[f]}$
which is $\Gamma$-invariant, then we can find a section
$\widetilde{\epsilon}$ of $T_[f] \Kbm{0,0}{\PP^1,e}$ which is
$\Gamma$-invariant and which maps to $\epsilon$.  Over an \'etale
neighborhood of the image of $[f]$ in $\kbm{0,0}{\PP^1,e}$, the stack
$\Kbm{0,0}{\PP^1,e}$ is a finite group quotient $[M/G]$ where $G$ is
the stabilizer group of $f$ and $M$ is a smooth scheme.  In particular
the invariant locus $M^\Gamma \subset M$ is a closed subscheme which
is smooth and whose Zariski tangent space at $[f]$ is the
$\Gamma$-invariant subspace of $T_{[f]} \Kbm{0,0}{\PP^1,e}$.  In
particular we can find a smooth, connected curve $B\subset M^\Gamma$
such that $B$ contains the point $[f]$ and the tangent space to $B$ at
$[f]$ equals $\text{span}(\widetilde{\epsilon})$.  But since
$\widetilde{\epsilon}$ has nonzero image in $N_{[f]}$, the curve $B$
is not contained in the image of $\Kbm{}{\PP^1,\sigma}$.  So the
general point of $B$ parametrizes a stable map with fewer nodes than
$f:C\rightarrow \PP^1$.  On the other hand, the invariant of the
Zariski tangent space is constant in connected families.  So by the
induction hypothesis, the invariant of $f$ equals the invariant of a
general point of $B$ is greater than $1$ (resp. at least $1$).

\medskip\noindent
By the last paragraph, we may now assume that for every node $q$ of
$C$, the $\Gamma$-invariant subspace of $\bigoplus_{i=0}^{s-1}
T'_{\phi^i q} \otimes T''_{\phi^i q}$ is trivial.  Let $\Delta \subset
\Gamma$ be the subgroup generated by $\phi^s$.  There is an action of
$\Delta$ on $T'_q \otimes T''_q$, and the $\Gamma$-representation
$N_{[f]}$ is simply the induced representation $\lt( T'_q \otimes
T''_q \rt) \otimes_{k[\Delta]} k[\Gamma]$.  By
Lemma~\ref{lem-induced}, it follows that the invariant of $N_{[f]}$ as
a $\Gamma$-representation is simply
\begin{equation}
\alpha_\gamma(N_{[f]}) = l\cdot \frac{s}{r} + \frac{s-1}{2}
\end{equation}
where the character $T'_q\otimes T''_q$ of $\Delta$ is
$\gamma|_{\Delta}^l$ for $l=0,\dots, \frac{r}{s}-1$.

\medskip\noindent
If $s \geq 3$, then already the invariant of $N_{[f]}$ is greater than
$1$.  So we may assume that $s=1$ or $s=2$.  If $s=2$, then the
invariant of $N_{[f]}$ is $\frac{1}{2} + l \cdot \frac{2}{r}$.  If
$s=1$, then the invariant is $l \cdot \frac{1}{r}$.  The only
possibilities that don't give an invariant larger than $1$ are:
\begin{enumerate}
\item every node of $C$ is fixed by $\phi$, or
\item there is precisely one pair of nodes $q, \phi q$ not fixed by
  $\phi$.
\end{enumerate}
We will consider these two possibilities in turn.

\medskip\noindent
Suppose that every node of $C$ is fixed by $\phi$.  Then every
irreducible component of $C$ is left invariant by $\phi$.  Since
$\phi$ is nontrivial, there is an irreducible component $C_i$ of $C$
such that $\phi|_{C_i}$ is nontrivial.  Let $C_1$ be an irreducible
component so that the restriction $\phi|_{C_1}$ has maximal order
$r_1$ (i.e. $\phi|_{C_1}^m = \text{Id}$ iff $r_1$ divides $m$).  The
irreducible component $C_1$ contains at least one node $q$ of $C$ and
it does not contain more than two nodes $q$ of $C$ since the only
automorphism of $\PP^1$ which fixes three points is the identity.  In
particular, $C_1$ is not contracted by $f$.

\medskip\noindent
Let $\tau$ be the stable $A$-graph of $C$ and let $\psi:\tau
\rightarrow \sigma$ be the maximal contraction which does not contract
the edges corresponding to nodes on $C_1$.  There is a morphism
$\Kbm{}{\PP^1,\sigma} \rightarrow \Kbm{0,0}{\PP^1,e}$ and the normal
bundle, as mentioned above, is the direct sum over nodes $q$ on $C_1$
of $T'_q \otimes T''_q$.  Let $\xi: \sigma \hookrightarrow \tau'$ be
the combinatorial morphism which is the inclusion of the maximal
sub-$A$-graph of $\tau$ whose only vertex is $v_1$, the vertex of the
irreducible component of $C_1$.  More simply, $\tau'$ is the graph
with the single vertex $v_1$ corresponding to $C_1$ and with one tail
for each node $q$ of $C$ contained in $C_1$ (i.e. either one or two
tails depending on whether $C_1$ contains one or two nodes of $C$).
Let $f_1:(C_1,q) \rightarrow \PP^1$ or $f_1:(C_1,q,q') \rightarrow
\PP^1$ be the marked stable map which is the restriction of $f$ to the
irreducible component $C_1$ marked by the nodes of $C$ contained on
$C_1$.

\medskip\noindent
We have a commutative diagram of $1$-morphisms:
\begin{equation}
\begin{CD}
\Kbm{}{\sigma, \PP^1} @>>> \Kbm{}{\sigma,\PP^n} \\
@V \Kbm{}{\xi, \PP^1} VV @VV \Kbm{}{\xi, \PP^n} V \\
\Kbm{}{\tau', \PP^1} @>>> \Kbm{}{\tau', \PP^n}
\end{CD}
\end{equation}
The horizontal arrows are closed immersions and the vertical arrows
are smooth.  So the invariant of the tangent space $T_{[f]}
\Kbm{}{\sigma,\PP^n}$ is greater 
than or equal to the invariant of the tangent space $T_{[f_1]}
\Kbm{}{\tau',\PP^n}$.  Let $e_1$ be the degree of $f_1:C_1 \rightarrow
\PP^1$.  As a $\Gamma$-representation, $T_{[f_1]} \Kbm{}{\tau', \PP^n}$
is the direct sum of $T_{[f_1]} \Kbmo{0,0}{\PP^n,e_1}$ with the tangent
space $T_q C_1$ (or $T_q C_1 \oplus T_{q'} C_1$ if $C_1$ contains two
nodes).  By Proposition~\ref{prop-RSBT1}, the invariant of $T_{[f_1]}
\Kbmo{0,0}{\PP^n,e_1}$ is $\frac{e_1(n+1)}{2}\lt( 1- \frac{1}{r_1}
\rt) - 1$.  Except for the four cases $(e_1,n) = (2,1), (2,2), (2,3),
(3,2)$, this invariant is already greater than $1$.  So the invariant
of $T_[f] \Kbm{0,0}{\PP^n,e}$ is also greater than $1$.  So we
consider these four cases.

\medskip\noindent
Now if $(e_1,n) = (2,3)$ or $(3,2)$ the invariant
$\frac{e_1(n+1)}{2}\lt( 1- \frac{1}{r_1} \rt) - 1$ equals $1$.  And
then the invariant of $T_q C$ is either $\frac{1}{2}$ for $(2,3)$ or
$\frac{1}{3}$ or $\frac{2}{3}$ for $(3,2)$.  Therefore the invariant
of $\Kbm{}{\tau', \PP^n}$ is greater than $1$.  So the invariant of
$T_{[f]} \Kbm{0,0}{\PP^n,e}$ is also greater than $1$.

\medskip\noindent
If $(e_1,n) = (2,2)$, the invariant of $T_{[f_1]} \Kbmo{0,0}{\PP^2,2}$
is $\frac{1}{2}$.  The invariant of $T_q C_1$ is also $\frac{1}{2}$.
So if there are two nodes on $C_1$, the invariant is already greater
than $1$.  But if there is only one node, so far the invariant only
equals $1$.  But we also have the invariant of $T'_q \otimes T''_q$,
which is positive.  So the invariant of $T_{[f]} \Kbm{0,0}{\PP^n,e}$
is greater than $1$.

\medskip\noindent
Finally, suppose $(e_1,n) = (2,1)$.  The invariant of $T_{[f_1]}
\Kbmo{0,0}{\PP^1,2}$ is zero.  The invariant of $T_q C_1$ is
$\frac{1}{2}$.  If there are two nodes on $C_1$, then the invariant is
$1$ and then the invariants of $T'_q\otimes T''_q$ and $T'_{q'}
\otimes T''_{q'}$ will make the total invariant of $T_{[f]}
\Kbm{0,0}{\PP^1,e}$ positive.  Therefore assume there is only one
node.  Then the invariant of $T_{[f_1]} \Kbm{}{\PP^1,\sigma}$ only
equals $\frac{1}{2}$.  Consider the $\Gamma$-representation $T'_q
\otimes T''_q$.  Let $C_2$ denote the irreducible component of $C$
which intersects $C_1$ at $q$.  Because of our assumption that the
order $r_1$ is the maximal among all orders of $\phi|_{C_i}$, either
$\phi$ acts trivially on $C_2$ or the order of $\phi|_{C_2}$ is $2$.
But in the second case, both $T'_q$ and $T''_q$ give characters of
$\Gamma$ which are $\gamma^{\frac{r}{2}}$, the unique character of
order $2$.  So the tensor product is the trivial character.  This
violates our assumption that for every node there are no non-zero
$\Gamma$-invariant sections of $\bigoplus_{i=0}^{s-1} T'_{\phi^i q}
\otimes T''_{\phi^i q}$.  Therefore $\phi$ acts trivially on $C_2$ and
the invariant of $T'_q \otimes T''_q$ is $\frac{1}{2}$.  So the
invariant of $T_{[f]} \Kbm{0,0}{\PP^1,e}$ is at least $\frac{1}{2} +
\frac{1}{2} = 1$.  Unfortunately, this is all we can conclude -- it is
easy to write down degree $e \geq 3$ covers of $\PP^1$ with reducible
domain where the invariant is equal to $1$.  So for $n=1$ and $e\geq
3$, we can only conclude that the invariant is $\geq 1$,
i.e. $\Kbm{0,0}{\PP^1,e}$ has canonical singularities (of course we
still have to dispense with the case that there are two nodes $q,\phi
q$ interchanged by $\phi$!).

\medskip\noindent
This finishes the analysis when $\phi$ fixes all nodes.  Now we
suppose that there is exactly one pair of nodes $\{ q, \phi q\}$ which
are not fixed by $\phi$.  The node $q$ disconnects $C$ into a union of
two connected subcurves $D_q$ and $C_2$.  Let $C_2$ be the subcurve
which contains $\phi q$.  The node $\phi q$ disconnects $C_2$ into a
union of two connected subcurves $D_{\phi q}$ and $C_1$.  Let $C_1$ be
the subcurve which contains $q$.  So $C_1$ is the maximal connected
subcurve of $C$ containing $q$ and $\phi q$ on which both $q$ and
$\phi q$ are nonsingular points.  Observe that $\phi(D_q) = D_{\phi
q}$, $\phi(D_{\phi q}) = D_q$ and both $D_q$ and $D_{\phi q}$ are
smooth.  Let $C_q \subset C_1$ be the irreducible component which
contains $q$ and let $C_{\phi q} \subset C_1$ be the irreducible
component which contains $\phi q$.  Observe that $\phi(C_q) = C_{\phi
q}$ and $\phi(C_{\phi q}) = C_q$.  There are two possibilities
depending on whether $C_q$ (and thus $C_{\phi q}$) is contracted by
$f$ or not.

\medskip\noindent
First consider the possibility that $C_q$ is contracted by $f$ and
suppose that $C_q \neq C_{\phi q}$. Then $C_q$ contains at least three
nodes, $q$ and two other nodes.  By assumption, each of the two other
nodes is fixed by $\phi$.  Also we have $\phi(C_q) = \phi(C_{\phi
q})$.  Since $C_q$ is not equal to $C_{\phi q}$, then $C_q \cap
C_{\phi q}$ is at most one node.  But then the second of the other
nodes cannot be fixed by $\phi$, which is a contradiction.  So we
conclude that if $C_q$ is contracted by $f$, then $C_q = C_{\phi q}$
and $C_q$ contains at least one other node $r$ of $C_1$ which is one
of the two fixed points of $\phi|_{C_q}$.  Also, since $\phi|_{C_q}$
contains the orbit $\{q, \phi q\}$ of order $2$, we conclude that
$\phi|_{C_q}$ has order $2$.

\medskip\noindent
The node $r$ disconnects $C_1$ into $C_q$ and a connected subcurve
$C_0$.  Suppose that $\phi|_{C_0}$ is the identity.  Then the
$\Gamma$-representation $T'_r \otimes T''_r$ has invariant
$\frac{1}{2}$.  Combined with the invariant $\frac{1}{2} + \frac{l
\cdot r}{2}$ coming from the nodes $\{ q, \phi q\}$, the total
invariant is greater than $1$ so that the invariant of $T_{[f]}
\Kbm{0,0}{\PP^n,e}$ is greater than $1$.  Next suppose that
$\phi|_{C_0}$ is not the identity.  Then by the same analysis as in
the case that $\phi$ fixes all nodes of $C$, we conclude that the
invariant coming from $C_0$ is at least $\frac{1}{2}$ (except when
$(e_1,n) = (2,1)$, the invariant is at least $1$).  Combined with the
invariant $\frac{1}{2} + \frac{l \cdot r} {2}$ coming from the nodes
$\{ q, \phi q\}$, the total invariant is greater than $1$ so that the
invariant of $T_{[f]} \Kbm{0,0}{\PP^n,e}$ is greater than $1$.

\medskip\noindent
We are reduced to the case when $f$ does not contract $C_q$.  By the
same analysis as above, we conclude that $D_q$ and $D_{\phi q}$ are
irreducible and are not contracted by $f$.  Let $\tau$ be the stable
$A$-graph of $f$ and let $\psi:\tau \rightarrow \sigma$ be the maximal
contraction which does not contract the edges corresponding to the
nodes $q$ and $\phi q$.  Then $\sigma$ has three vertices: $v_0$
corresponding to the connected subcurve $D_q$, $v_1$ corresponding to
the connected subcurve $C_1$ and $v_2$ corresponding to the connected
subcurve $D_{\phi q}$.  The stable map $f:C \rightarrow \PP^1$
determines a point of the Behrend-Manin stack, $\Kbm{}{\PP^n,
\sigma}$.  The Zariski tangent space of $\Kbm{}{\PP^n, \sigma}$ at
$[f]$ is a $\Gamma$-sub-representation of the Zariski tangent space of
$\Kbm{0,0}{\PP^n,e}$ and the cokernel is $(T'_q \otimes T''_q) \oplus
(T'_{\phi q} \otimes T''_{\phi q})$.

\medskip\noindent
Let $\xi:\sigma \hookrightarrow \tau'$ be the maximal
\emph{disconnected} subgraph of $\sigma$ which contains the vertices
$v_0$ and $v_2$, i.e. the stable $A$-graph with vertices $v_0$ and
$v_2$ and one flag attached to each vertex corresponding to the marked
point $q$ and $\phi q$ respectively.  There is an associated
$1$-morphism $\Kbm{}{\PP^n,\sigma} \rightarrow \Kbm{}{\PP^n,\tau'}$.
Because $C_1$ is not contracted by $f$, this $1$-morphism is smooth.
In particular, the invariant of $T_{[f]} \Kbm{}{\PP^n,\sigma}$ is at
least as large as the invariant of $T_{[f]} \Kbm{}{\PP^n,\tau'}$.  Of
course $\Kbm{}{\PP^n,\tau'}$ is simply a product of the factor from
$v_0$ and $v_2$, $\Kbm{0,1}{\PP^n,e_1} \times \Kbm{0,1}{\PP^n,e_1}$
where $e_1$ is the degree of $f|_{D_q}:D_q \rightarrow \PP^n$.  The
Zariski tangent space is correspondingly a direct sum of the two
factors coming from $v_0$ and $v_2$.  The automorphism $\phi$ permutes
the two factors and $\phi^2$ acts as an automorphism of each factor.

\medskip\noindent
Let $\gamma$ be a generator for the character group of $\Gamma =
\langle \phi \rangle$.  Then also $\gamma$ is a generator for the
character group of $\langle \phi^2 \rangle$.  The rank of $T
\Kbm{0,1}{\PP^n,e_1}$ equals $(n+1)e_1 + (n-3) + 1$.  Consider the
invariant $\alpha'_\gamma$ of this $\langle \phi^2
\rangle$-representation with respect to $\gamma$.  One contribution
comes from the marked point $q$, which is positive.  So the invariant
is positive.  By Lemma~\ref{lem-induced}, the invariant of $T
\Kbm{}{\PP^n,\tau'}$ as a $\Gamma$-representation with respect to
$\gamma$ is
\begin{equation}
\alpha_\gamma( T \Kbm{}{\PP^n,\tau'}) \geq \alpha'_\gamma + \lt(
(n+1)e_1 + (n-3) + 1 \rt)\frac{1}{2}
\end{equation}
The right-hand-side of the equation is a minimum when $n=1$ and $e_1 =
1$, in which case it is still larger than $\frac{1}{2}$ (since
$\alpha'_\gamma$ is positive).  So the invariant of
$\Kbm{}{\PP^n,\sigma}$ is larger than $\frac{1}{2}$.  And the
invariant of $\lt( T'_q \otimes T''_q \rt) \oplus \lt( T'_{\phi q}
\otimes T''_{\phi q} \rt)$ is larger than $\frac{1}{2}$.  Therefore
the invariant of $T_{[f]} \Kbm{0,0}{\PP^n,e}$ is larger than $1$.
This finishes the proof.
\end{proof}

\medskip\noindent
Now we begin the analysis of the coarse moduli space of $C_d$.  As in
Section~\ref{sec-proof}, let $(\varrho:M \rightarrow \PP^1, \iota: Y
\times \PP^1 \rightarrow M, B_y \rightarrow M, \PP(N\oplus \mathbf{1})
\rightarrow M)$ denote the deformation to the normal cone associated
to the inclusion $Y \hookrightarrow B$, where $B=\Kbm{0,0}{\PP^n,e}$.
Let $\widetilde{\pi}_d: \widetilde{C}_d \rightarrow M^o$ denote the
projective Abelian cone.

\medskip\noindent
\begin{lem} \label{lem-can} 
If $e\geq 3$ and if $d+ e \leq n$, then the map to the coarse moduli
space
\begin{equation}
\widetilde{C}_d \times_{M^o} N \rightarrow (\widetilde{C}_d
\times_{M^o} N )_{\text{coarse}} 
\end{equation}
is an isomorphism away from codimension $2$, and the coarse moduli
space $( \widetilde{C}_d \times_{M^o} N )_{\text{coarse}}$ is a
normal, $\QQ$-Gorenstein variety with only canonical singularities.
\end{lem}

\medskip\noindent
\begin{proof}
Of course $\widetilde{C}_d \times_{M^o} N$ is normal and Gorenstein,
therefore the coarse moduli space is normal and $\QQ$-Gorenstein.  To
see that the coarse moduli map is an isomorphism away from codimension
$2$ and that the coarse moduli space is canonical, observe first that
we have constructed a resolution of $\widetilde{C}_d \times_{M^o} N$
as a Deligne-Mumford stack.  Recall the resolution is constructed as
follows.  First of all, the projection $N \rightarrow Y$ is
$M^{(0)}(Y,\mc{R}^\vee, \text{pr}_{\mathbb{G}}^* T)$.  By
Proposition~\ref{prop-smooth}, we have a $1$-morphism of stacks
$u^{e-1,0}:M^{(e-1)} \rightarrow N$ such that $M^{(e-1)} \rightarrow
Y$ is representable and smooth.  There is a projective bundle $C_d'$
over $M^{(e-1)}$ and a morphism $C_d' \rightarrow \widetilde{C}_d
\times_{M^o} N$ which is a resolution of singularities.  Observe that
also $C_d' \rightarrow Y$ is representable and smooth.  By
Corollary~\ref{cor-mults}, the relative canonical divisor of $C_d'
\rightarrow \widetilde{C}_d \times_{M^o} N$ is effective.

\medskip\noindent
Now consider the coarse moduli spaces $Y \rightarrow
Y_{\text{coarse}}$, $C_d' \rightarrow C_{d,\text{coarse}}'$ and
$\widetilde{C}_d \times_{M^o} N \rightarrow (\widetilde{C}_d
\times_{M^o} N)_\text{coarse}$.  By Proposition~\ref{prop-RSBT2}, the
morphism $Y \rightarrow Y_{\text{coarse}}$ is an isomorphism away from
codimension $2$ and $Y_{\text{coarse}}$ is canonical (this corresponds
to the case $(e,1)$ where $e\geq 3$).  By
Corollary~\ref{cor-DMsmooth}, also $C_d' \rightarrow
C_{d,\text{coarse}}'$ is an isomorphism away from codimension $2$ and
$C_{d,\text{coarse}}'$ is canonical.  There is an open substack $U
\subset \widetilde{C}_d \times_{M^o} N$ such that $C_d' \rightarrow
\widetilde{C}_d \times_{M^o} N$ is an isomorphism over $U$ and such
that the complement of $U$ has codimension at least $2$.  And the
morphism $U \rightarrow U_{\text{coarse}}$ is an isomorphism away from
codimension $2$.  Therefore also the morphism
\begin{equation}
\widetilde{C}_d \times_{M^o} N \rightarrow (\widetilde{C}_d
\times_{M^o} N )_{\text{coarse}} 
\end{equation}
is an isomorphism away from codimension $2$.  

\medskip\noindent
Because $C_d' \rightarrow C_{d,\text{coarse}}'$ is an isomorphism away
from codimension $2$, the relative canonical divisor of
$C_{d,\text{coarse}}' \rightarrow (\widetilde{C}_d \times_{M^o}
N)_\text{coarse}$ equals the image of the canonical divisor of $C_d'
\rightarrow \widetilde{C}_d \times_{M^o} N$.  Therefore the relative
canonical divisor is effective.  But also $C_{d,\text{coarse}}'$ is
canonical.  It follows that also $(\widetilde{C}_d \times_{M^o}
N)_\text{coarse}$ is canonical.  This finishes the proof.
\end{proof}

\medskip\noindent
\begin{rmk} \label{rmk-can} 
When $e=2$ and $d + 3 \leq n$ (note that this inequality is worse than
the usual $d+ e \leq n$), the second part of the lemma also holds by a
slightly more ad hoc argument.  In this case $Y$ is a $\ZZ/2\ZZ$-gerbe
over its coarse moduli space.  And $N = M^{(0)}$ is a vector bundle of
$1\times (n-1)$-matrices.  So the only stratum to blow up to form
$M^{(1)}$ is the zero section.  When we do this, we have that the
$\ZZ/2\ZZ$-invariant locus of $M^{(1)}$ is the whole exceptional
divisor $E$.  A simple computation shows that the $1$-morphism $C_d'
\rightarrow M^{(1)}$ preserves all stabilizer groups of geometric
points (i.e. the induced homomorphisms of stabilizer groups are
isomorphisms).  Therefore the morphism $C_d' \rightarrow
C_{d,\text{coarse}}'$ is a morphism to a smooth variety ramified of
ramification index $1$ along the preimage of $E$, i.e. it does not
satisfy the first part of the lemma.  However, it is straightforward
to compute that the relative canonical divisor of
$C_{d,\text{coarse}}' \rightarrow (\widetilde{C}_d \times_{M^o}
N)_\text{coarse}$ is $\frac{n-3-d}{2} E_{\text{coarse}}$.  Therefore
$(\widetilde{C}_d \times_{M^o} N)_\text{coarse}$ is canonical when $d+
3 \leq n$.
\end{rmk}  

\section{Conjectures about $\kbm{0,0}{X,e}$} \label{sec-conj}

\medskip\noindent
\begin{conj}[Inversion of Adjunction, Conj. 7.3 ~\cite{KollarPairs}]
  \label{conj-ioa} 
Let $X$ be a normal variety, $S$ a normal Cartier divisor and $B=\sum
b_i B_i$ a $\QQ$-divisor.  Assume that $K_X + S + B$ is $\QQ$-Cartier.
Then
\begin{equation}
\text{totaldiscrep}(S, B|S) = \text{discrep}(\text{Center} \cap S \neq
\emptyset, X, S+B),
\end{equation}
where the notation on the right means that we compute the discrepancy
using only those divisors whose center on $X$ intersects $S$.
\end{conj}

\medskip\noindent
\begin{conj} \label{conj-1} 
If $e \geq 3$ and $d+e \leq n$, the coarse moduli space
$C_{d,\text{coarse}}$ is a normal, $\QQ$-Gorenstein variety with only
canonical singularities.  If $e=2$ and $d+3 \leq n$, the coarse moduli
space $C_{d,\text{coarse}}$ is a normal, $\QQ$-Gorenstein variety with
only canonical singularities.
\end{conj}

\medskip\noindent
\begin{conj} \label{conj-2} 
If $e \geq 3$, $d+e \leq n$, and if $X\subset \PP^n$ is a general
hypersurface of degree $d$, then $\kbm{0,0}{X,e}$ is a normal,
$\QQ$-Gorenstein variety with only canonical singularities.  If $e=2$,
$d+3 \leq n$ and if $X \subset \PP^n$ is a general hypersurface of
degree $d$, then $\kbm{0,0}{X,2}$ is a normal, $\QQ$-Gorenstein
variety with only canonical singularities.
\end{conj}

\medskip\noindent
\begin{prop} \label{prop-ioac1} 
Suppose that $d+e \leq n$.
\begin{enumerate}
\item
If $e \geq 2$ then the coarse moduli space $( C_d)_{\text{coarse}}$ is
normal, $\QQ$-Gorenstein and Kawamata log terminal.
\item
If $e \geq 3$ and $d+e \leq n$, then the coarse moduli map
\begin{equation}
C_d \rightarrow ( C_d )_{\text{coarse}}
\end{equation}
is an isomorphism away from codimension $2$. 
\item
Conjecture~\ref{conj-ioa} implies Conjecture~\ref{conj-1}.
\end{enumerate}
\end{prop}

\medskip\noindent
\begin{proof}
First of all, the stack $C_d$ is normal and Gorenstein with canonical
singularities by Theorem~\ref{thm-main}.  So the coarse moduli space
$C_{d,\text{coarse}}$ is normal and $\QQ$-Gorenstein.  And by
~\cite[Prop. 3.16]{KollarPairs}, $C_{d,\text{coarse}}$ is Kawamata log
terminal.  This proves Item ($1$).

\medskip\noindent
Denote by $Z\subset C_d$ the closed substack where the map $C_d
\rightarrow (C_d )_{\text{coarse}}$ is not an isomorphism.  Denote by
$\widetilde{Z} \subset \widetilde{C}_d$ the closed substack where the
map $\widetilde{C}_d \rightarrow ( \widetilde{C}_d )_{\text{coarse}}$
is not an isomorphism.  Clearly $\widetilde{Z} \cap \rho^{-1}(\AAA^1)
= Z \times \AAA^1$.  Of course $Z$ is invariant under the action of
$\text{GL}_{n+1}$.  By the same argument as in the proof of
Theorem~\ref{thm-main}, every irreducible component of $\widetilde{Z}$
has non-empty intersection with the fiber over $\infty$,
i.e. $\widetilde{C}_d \times_{M^o} N$.  By Lemma~\ref{lem-can}, if
$e\geq 3$ then every irreducible component of $\widetilde{Z} \cap
(\widetilde{C}_d \times_{M^o} N)$ has codimension at least $2$.  By
Krull's Hauptidealsatz, we conclude that every irreducible component
of $\widetilde{Z}$ has codimension at least $2$.  Therefore every
irreducible component of $Z \subset C_d$ has codimension at least $2$.
This proves Item ($2$).

\medskip\noindent
Finally we prove Item ($3$).  As in the proof of
Theorem~\ref{thm-main}, let $W \subset \kbm{0,0}{\PP^n,e}$ be the
largest open substack over which $C_{d,\text{coarse}}$ is canonical.
This is a $\text{GL}_{n+1}$-invariant open set, so to prove that $W$
is all of $\kbm{0,0}{\PP^n,e}$, it suffices to prove that $W$ contains
the image of $Y$.

\medskip\noindent
Let $\widetilde{C}_d \rightarrow M^o$ be as in
Section~\ref{sec-proof}.  Let $M^o \rightarrow M^o_{\text{coarse}}$
and $\widetilde{C}_d \rightarrow \widetilde{C}_{d,\text{coarse}}$ be
the coarse moduli spaces.  Let $W'\subset M^o_{\text{coarse}}$ be the
largest open subset over which $\widetilde{C}_{d,\text{coarse}}$ is
canonical.  Of course $W' \cap \varrho^{-1}(\AAA^1) = W \times
\AAA^1$.  Now by Lemma~\ref{lem-can}, the Cartier divisor
$(\widetilde{C}_d \times_{M^o} N)_{\text{coarse}}$ in
$\widetilde{C}_{d,\text{coarse}}$ is normal and canonical.  Assuming
Conjecture~\ref{conj-ioa} is true, we conclude that there is an open
subvariety of $\widetilde{C}_{d,\text{coarse}}$ containing
$(\widetilde{C}_d \times_{M^o} N)_{\text{coarse}}$ which is canonical,
i.e. $W'$ contains $\varrho^{-1}(\infty)$.  By the same argument as in
the proof of Theorem~\ref{thm-main}, we conclude that $W$ contains the
image of $Y$, i.e. $W$ is all of $\kbm{0,0}{\PP^n,e}$.  So
Conjecture~\ref{conj-1} is true.
\end{proof}

\medskip\noindent
\begin{prop} \label{prop-c1c2} 
Let $d+e \leq n$.  Let $X \subset \PP^n$ be a general hypersurface of
degree $d$.
\begin{enumerate}
\item
If $e \geq 2$, the coarse moduli space $\kbm{0,0}{X,e}$ is normal,
$\QQ$-Gorenstein and Kawamata log terminal.
\item
If $e \geq 3$, then the coarse moduli map
\begin{equation}
\Kbm{0,0}{X,e} \rightarrow \kbm{0,0}{X,e}
\end{equation}
is an isomorphism away from codimension $2$.
\item
Conjecture~\ref{conj-1} implies Conjecture~\ref{conj-2}.
\end{enumerate}
\end{prop}

\medskip\noindent
\begin{proof}
This is the same argument as in the proof of Corollary~\ref{cor-main}.
\end{proof}

\medskip\noindent
\begin{rmk} \label{rmk-c1c2} 
In Item ($2$), if $e=2$ then the coarse moduli map fails to be an
isomorphism precisely on the locus $Y \cap \Kbm{0,0}{X,2}$.  By direct
computation, for $X$ general this locus has dimension $2n-d-1$.  And
$\Kbm{0,0}{X,2}$ has dimension $3n-d-2$.  Therefore, if $d +3\leq n$,
then the coarse moduli map $\Kbm{0,0}{X,2} \rightarrow \kbm{0,0}{X,2}$
is an isomorphism away from codimension $2$.
\end{rmk}

\section{The canonical class on $\Kbm{0,r}{X,e}$} \label{sec-canKBM}

\medskip\noindent 
Let $X\subset \PP^n$ be a complete intersection of $c$ hypersurfaces
of degrees $\mathbf{d} = (d_1,\dots,d_c)$.  Associated to the
inclusion morphism, there is a $1$-morphism of Kontsevich moduli
spaces $\Kbm{0,r}{X,e} \rightarrow \Kbm{0,r}{\PP^n,e}$.  This
$1$-morphism is representable and is a closed immersion.  The image is
the zero locus of a section $\sigma$ of a locally free sheaf
$\mc{P}_{\mathbf{d}}$ in the small \'etale site of
$\Kbm{0,r}{\PP^n,e}$.

\medskip\noindent
In the case that $\sigma$ is a \emph{regular} section, we may express
the dualizing sheaf $\omega'$ on $\Kbm{0,r}{X,e}$ as the pullback from
$\Kbm{0,r}{\PP^n,e}$ of the tensor product $\omega\otimes
\text{det}(\mc{P}_{\mathbf{d}})$ where $\omega$ is the dualizing sheaf
on $\Kbm{0,r}{\PP^n,e}$.  Pandharipande has computed the $\QQ$-Picard
group of $\Kbm{0,r}{\PP^n,e}$ in ~\cite{QDiv}.  And he has computed
the $\QQ$-divisor class of $\omega$ in ~\cite{Pand97}.  The purpose of
this section is to compute the class $\text{det}(\mc{P}_{\mathbf{d}})$
in terms of the standard generators of the $\QQ$-Picard group, and
thereby compute the $\QQ$-divisor class of $\omega'$ in the case that
$\sigma$ is a regular section.

\medskip\noindent
Let $p\colon \mathcal{C}\rightarrow \Kbm{0,r}{\PP^n,e}$ denote the
universal curve.  Let $f\colon \mathcal{C} \rightarrow \PP^n$ denote
the universal morphism.  For each integer $d$, we can form the
pullback sheaf $f^*\OO_{\PP^r}(d)$.  We define $\mathcal{E}_d$ to be
the pushforward $p_*(f^*\OO_{\PP^r}(d)$.  More generally, given an
ordered sequence $\mathbf{d} = (d_1,\dots,d_c)$ of integers, we define
$\OO_{\PP^n}(\mathbf{d})$ to be $\OO_{\PP^n}(d_1)\oplus \dots \oplus
\OO_{\PP^n}(d_c)$ and we define $\mc{P}_{\mathbf{d}}$ to be
$\mc{P}_{d_1}\oplus \dots \oplus \mc{P}_{d_c}$.

\medskip\noindent
There is a pullback map on global sections:
\begin{equation}
H^0\lt( \PP^n, \OO_{\PP^n}(d) \rt) \rightarrow H^0\lt( \mathcal{C},
f^*\OO_{\PP^n}(d) \rt) \rightarrow H^0\lt( \Kbm{0,r}{\PP^n,e},
\mc{P}_d \rt).
\end{equation}
We denote the composite map by $f^*$.  More generally, given an
ordered sequence $\mathbf{d}$, we have a pullback map on global
sections:
\begin{equation}
f^*\colon H^0\lt( \PP^n, \OO_{\PP^n}(\mathbf{d}) \rt) \rightarrow
H^0\lt( \Kbm{0,r}{\PP^n,e}, \mc{P}_{\mathbf{d}} \rt).
\end{equation}

\medskip\noindent
\begin{lem} \label{applem1} 
If $d\geq 0$, then $\mathcal{P}_d$ is a locally free sheaf of rank
$de+1$ in the small \'etale site of $\Kbm{0,r}{\PP^n,e}$ and
$R^ip_*\lt( f^*\OO_{\PP^n}(d) \rt)$ is zero for $i>0$.  More
generally, if $\textbf{d} = (d_1,\dots,d_c)$ and $d_1,\dots,d_c \geq
0$, then $\mc{P}_{\textbf{d}}$ is a locally free sheaf of rank
$|\textbf{d}|e+c$, where $|\mathbf{d}| = d_1 + \dots + d_c$ and
$R^ip_* \lt( f^*\OO_{\PP^n}(\textbf{d}) \rt)$ is zero for $i>0$.
\end{lem}

\medskip\noindent
\begin{proof}
This has been proved in other places, in particular it is proved as
part of the the proof of ~\cite[Lemma 4.5]{HRS2}.
\end{proof}

\medskip\noindent 
Now let $d_1,\dots,d_c$ be a sequence of positive integers, and let
$s=(s_1,\dots, s_c)$ be a global section of $H^0\lt( \PP^n,
\OO_{\PP^n}(\mathbf{d}) \rt)$.  Let $X\subset \PP^n$ be the zero locus
of $s$.  Let $\sigma$ denote the pullback section $f^*s\in H^0\lt(
\Kbm{0,r}{\PP^n,e}, \mc{P}_{\mathbf{d}}\rt).$

\medskip\noindent
\begin{lem}\label{applem2} 
The zero locus of $\sigma$ as a closed substack of
$\Kbm{0,r}{\PP^n,e}$ is the image of the closed immersion
$\Kbm{0,r}{X,e} \rightarrow \Kbm{0,r}{\PP^n,e}$.
\end{lem}

\medskip\noindent
\begin{proof}
This is also proved as part of the proof of ~\cite[Lemma 4.5]{HRS2}.
\end{proof}

\medskip\noindent
Since $\mc{P}_{\mathbf{d}}$ is a locally free sheaf of rank
$e|\mathbf{d}| + c$, and since $\Kbm{0,r}{\PP^n,e}$ is smooth, we
conclude that $\sigma$ is a regular section iff the codimension of
$\Kbm{0,r}{X,e}$ in $\Kbm{0,r}{\PP^n,e}$ equals $e|\mathbf{d}|+c$.  In
this case, it follows from the generalized version of the adjunction
theorem that the dualizing sheaf $\omega'$ on $\Kbm{0,r}{X,e}$ is the
pullback of the sheaf $\omega\otimes \text{det}(\mc{P}_{\mathbf{d}})$,
where $\omega$ is the dualizing sheaf on $\Kbm{0,r}{\PP^n,e}$.

\medskip\noindent
In ~\cite{QDiv}, Pandharipande described the $\QQ$-Picard group
$\text{Pic}\lt(\Kbm{0,r}{\PP^n,e}\rt) \otimes \QQ$.  For simplicity,
we assume that $n>1$ and $e>0$ and also that $(n,e) \neq (2,2)$.  The
divisor class $\mathcal{H}$ is defined as the image of a positive
generator $h^2\in \text{CH}^2(\PP^n)$ under the composition
\begin{equation}
\text{CH}^2(\PP^n) \xrightarrow{f^*} \text{CH}^2(\mathcal{C})\otimes
\QQ \xrightarrow{p_*} \text{CH}^1(\Kbm{0,r}{\PP^n,e})\otimes \QQ.
\end{equation}
For each of the $r$ sections $g_i:\Kbm{0,r}{\PP^n,e} \rightarrow
\mathcal{C}$, the divisor class $\mathcal{L}_i$ is defined as the
image of a positive generator $h\in \text{CH}^1(\PP^n)$ under pullback
by $f\circ g_i$.  Finally, for each \emph{weighted partition}
$P=(A\cup B,e_A,e_B)$ of $(\{1,\dots,r\},e)$ there is the class
$\Delta_P$ of the corresponding boundary stratum of
$\Kbm{0,r}{\PP^n,e}$.  A weighted partition is a datum where $A\cup B$
is a partition of $\{1,\dots, r\}$, where $e_A + e_B = e$ with
$e_A,e_B \geq 0$, and where we demand that $|A| \geq 2$ (resp. $|B|
\geq 2$) if $e_A = 0$ (resp. $e_B = 0$).  The corresponding boundary
stratum is the closure of the locally closed substack parametrizing
stable maps whose dual graph is of type $(A\cup B, e_A, e_B)$.
Pandharipande's result is that the $\QQ$-Picard group is a
$\QQ$-vector space with basis
\begin{equation}
\{\mathcal{H}\} \cup \{ \mathcal{L}_i | i=1,\dots, r\} \cup \{
\Delta_P | P = (A\cup B, e_A, e_B) \}.
\end{equation}
In the case that $r=0$, for $i=0,\dots,\lt[\frac{e}{2}\rt]$ we denote
by $\mathcal{D}_{i,0}$ the $\QQ$-divisor class $\Delta_{P}$ where $P
=(\emptyset \cup \emptyset, i, e-i)$.  And for $r>0$, for
$i=0,\dots,\lt[ \frac{e}{2}\rt]$ and $j=0,\dots, r$ we define
$W_{i,j}$ to be the set of weighted partitions $\{(A\cup B, e_A,e_B)|
|A| = j, e_A = i\}$.  We denote by $\mathcal{D}_{i,j}$ the
$\QQ$-divisor class
\begin{equation}
\mathcal{D}_{i,j} = \sum_{P\in W_{i,j}} \Delta_P.
\end{equation}

\medskip\noindent
In ~\cite{Pand97}, Pandharipande computed the $\QQ$-divisor class of
the dualizing sheaf $\omega$ in terms of the basis above.

\medskip\noindent
\begin{prop}[Pandharipande, Prop. 2 ~\cite{Pand97}] \label{Pand1}
The dualizing sheaf $\omega$ on $\Kbm{0,0}{\PP^n,e}$ has $\QQ$-divisor
class
\begin{equation}
\omega = \frac{1}{2e}\lt[ -(n+1)(e+1)\mathcal{H}+
\sum_{i=1}^{\lt[ \frac{e}{2}\rt] } \lt(( n+1)(e-i)i - 4e\rt)
\mathcal{D}_{i,0} \rt]
\end{equation}
\end{prop}

\medskip\noindent
\begin{prop}[Pandharipande, Prop. 3 ~\cite{Pand97}] \label{Pand2}
The first Chern class of the dualizing sheaf $\omega$ on
$\Kbm{0,r}{\PP^n,e}$ has $\QQ$-divisor class
\begin{eqnarray*}
C_1(\omega) = \frac{1}{2e^2}\lt[ -(n+1)(e+1)e + 2r\rt] \mathcal{H}
-\frac{1}{2e} \sum_{p=1}^{n} \mathcal{L}_p + \\
\frac{1}{2e^2}\sum_{i=0}^{\lt[\frac{e}{2}\rt]} \sum_{j=0}^{r} \lt[
(n+1)e(e-i)i + 2e^2j - 4eij + 2ri^2 - 4e^2 \rt] \mathcal{D}_{i,j}.
\end{eqnarray*}
\end{prop}

\medskip\noindent
It remains to compute the $\QQ$-divisor class of the first Chern class
$C_1(\mc{P}_{\textbf{d}})$.  We begin by computing for each integer
$d\geq 0$, the first Chern class $C_1(\mc{P}_d)$.  We may compute this
using the Grothendieck-Riemann-Roch formula ~\cite[Thm. 15.2]{F}.
Observe that $p$ is a representable morphism between smooth
Deligne-Mumford stacks with projective coarse moduli space.  So one
can deduce Grothendieck-Riemann-Roch for $p$ from
Grothendieck-Riemann-Roch for the coarse moduli spaces using
~\cite{Vis}.  Alternatively, one can use the Grothendieck-Riemann-Roch
theorem of Toen ~\cite{Toen}.

\medskip\noindent
By Lemma~\ref{applem1}, the element in $K$-theory,
$Rp_{!}[f^*\OO_{\PP^n}(d)]$ equals $[\mc{P}_d]$.  So, by the
Grothendieck-Riemann-Roch formula, we have
\begin{equation}
\text{ch}[\mc{P}_d] = p_*\lt( f^*\text{ch}[\OO_{\PP^n}(d)]\cap
\text{todd}(p) \rt).
\end{equation}

\medskip\noindent
Let us denote by $h\in \text{CH}^1(\PP^n)$ the first Chern class of
$\OO_{\PP^n}(1)$.  Then, up to terms in $\text{CH}^3(\PP^n)$, we have
the formula
\begin{equation}
\text{ch}[\OO_{\PP^n}(d)] = 1 + dh + \frac{d^2}{2}h^2 + \dots
\end{equation}
By ~\cite[Section 3.E]{HarrMorr}, up to terms in
$\text{CH}^3(\mathcal{C})\otimes \QQ$, we have the formula
\begin{equation}
\text{todd}(p) = 1 - \frac{1}{2}C_1(\omega_p) + \frac{1}{12}\lt( \eta
+ C_1(\omega_p)^2 \rt) + \dots 
\end{equation}
where $\eta$ is the $\QQ$-divisor class of the ramification locus of
$p$.  By ~\cite[Lemma 2.1.2]{QDiv}, $p_*\lt( \eta + C_1(\omega_p)^2
\rt)$ equals zero.  Therefore, up to terms in
$\text{CH}^2(\Kbm{0,r}{\PP^n,e})\otimes \QQ$, we have the formula
\begin{eqnarray*}
 p_*\lt( f^*\text{ch}[\OO_{\PP^n}(d)]\cap \text{todd}(p) \rt) = \ \ \
 \ \ \ \\
p_*\lt(df^*(h) - \frac{1}{2}C_1(\omega_p) \rt) +
\frac{d^2}{2}p_*f^*(h^2) -\frac{d}{2} p_*\lt(f^*(h)\cap C_1(\omega_p)
\rt).
\end{eqnarray*}
Clearly $p_*(f^*(h))$ is just $e$ and $p_*(C_1(\omega_p))$ is just
$-2$.  By definition, $p_*f^*(h^2)$ is the divisor class
$\mathcal{H}$.  

\medskip\noindent
\begin{lem}\label{appChernpre} 
In the $\QQ$-Picard group of $\Kbm{0,r}{\PP^n,e}$, we have the formula
\begin{equation}
p_* \lt( f^*(h)\cap C_1(\omega_p) \rt) = \frac{1}{d}\lt[ -\mathcal{H}
+ \sum_{i=1}^{\lt[ \frac{e}{2} \rt]}\sum_{j=0}^{r} (e-i)i
\mathcal{D}_{i,j} \rt].
\end{equation}
\end{lem}

\medskip\noindent
\begin{proof}
Denote by $\alpha$ the difference of the right-hand side of the
equation from the left-hand side.  So the proposition may be rephrased
as saying that $\alpha$ is zero in the $\QQ$-Picard group.

\medskip\noindent
The method of proof is the same as in ~\cite[Section 1.2]{QDiv}.
Consider the class $\mathcal{S}$ of all pairs $(B,\zeta)$ where $B$ is
a smooth complete curve, $\zeta:B \rightarrow \Kbm{0,r}{\PP^n,e}$ is a
$1$-morphism and such that $(B,\zeta)$ satisfies
\begin{enumerate}
\item 
for the pullback of the universal curve, $p_\zeta:\mathcal{C}_\zeta
\rightarrow B$, $\mathcal{C}_\zeta$ is a smooth surface,
\item 
the general fiber of $p_\zeta$ is a smooth, rational curve,
\item 
every singular fiber of $p_\zeta$ has exactly two irreducible
components,
\item 
blowing down one irreducible component in each singular fiber yields a
ruled surface over $B$.
\end{enumerate}
In ~\cite{QDiv}, it is proved that for any nonzero divisor class
$\beta$ in the $\QQ$-Picard group of $\Kbm{0,r}{\PP^n,e}$, there is a
pair $(B,\zeta)$ in $\mathcal{S}$ such that $\zeta^*(\beta)$ has
nonzero degree on $B$.  So to prove the proposition, it suffices to
prove that for every pair $(B,\zeta)$ in $\mathcal{S}$,
$\zeta^*(\alpha)$ has degree zero.

\medskip\noindent
Suppose $(B,\zeta)$ is in $\mathcal{S}$.  Let $(E_1\cup
E'_1,\dots,E_m\cup E'_m)$ denote the irreducible components of the
singular fibers of $p_\zeta$.  Let $s:B\rightarrow \mathcal{C}_\zeta$
denote a section of $p_\zeta$ which does not intersect any of
$E_1,\dots, E_m$ (by item (4), such a section exists), and let $S$
denote $s(B)$.  Let $F$ denote any smooth fiber of $p_\zeta$.  In the
group of numerical equivalence classes, $N^1(\mathcal{C}_\zeta)$, the
classes $[S],[F], [E_1],\dots, [E_m]$ give a basis for
$N^1(\mathcal{C}_\zeta)$ as a free $\ZZ$-module.

\medskip\noindent
We define $k = -\text{deg}([S]\cap [S])$.  One computes that
\begin{equation}
\lt\{ \begin{array}{ll}
\text{deg}([F]\cap [F])     & = \ 0 \\ 
\text{deg}([F]\cap [S])     & = \ 1 \\
\text{deg}([F]\cap [E_i])   & = \ 0 \\
\text{deg}([S]\cap [E_i])   & = \ 0 \\
\text{deg}([E_i]\cap [E_i]) & = -1 \\
\text{deg}([E_i]\cap [E_j]) & = \ 0,\ \  i\neq j
\end{array} \rt.
\end{equation}
By the adjunction formula, we have that $\lt( \zeta^*\omega_p \otimes
\OO_{\mathcal{C}_\zeta}(E_i) \rt)|_{E_i} \cong \omega_{E_i}$.
Therefore, $\text{deg}\lt( \zeta^*C_1(\omega_p) + [E_i]\rt) \cap [E_i]
= -2$, i.e. $\text{deg}(\zeta^*C_1(\omega_p)\cap [E_i]) = -1$.  
Similarly, we have that $\lt( \zeta^*\omega_p \otimes
\OO_{\mathcal{C}_\zeta}(F) \rt)|_{F} \cong \omega_F$.  Therefore
$\text{deg}\lt( \zeta^*C_1(\omega_p) + [F] \rt)\cap [F] = -2$,
i.e. $\text{deg}( \zeta^*C_1(\omega_p) \cap [F]) = -2$.  Finally, by
adjunction we have that $\lt( \zeta^*\omega_p \otimes
\OO_{\mathcal{C}_\zeta}(S) \rt)|_S$ is isomorphic to the relative
dualizing sheaf of $p_\zeta|_{S}:S \rightarrow B$.  But this is an
isomorphism, so the relative dualizing sheaf is just $\OO_S$.
Therefore we have $\text{deg}\lt( \zeta^*C_1(\omega_p) + [S]\rt) \cap
[S] = 0$, i.e. $\text{deg}( \zeta^*C_1(\omega_p) \cap [S]) = k$.
Putting this all together, we have that the numerical equivalence class
of $\zeta^*C_1(\omega_p)$ equals
\begin{equation}
\zeta^* C_1(\omega_p) = -2[S] -k[F] + \sum_{i=1}^m [E_i].
\end{equation}

\medskip\noindent
Now define $l = \text{deg}(\zeta^*f^*(h))\cap [S]$.  For each
$i=1,\dots, m$, define $e_i = \text{deg}(\zeta^*f^*(h)\cap [E_i])$.
By a similar computation as above, we have that the numerical
equivalence class of $\zeta^* f^*(h)$ equals
\begin{equation}
\zeta^* f^*(h) = e[S] + (l + ek)[F] - \sum_{i=1}^m e_i[E_i].
\end{equation}
So we compute that $\text{deg}(\zeta^*C_1(\omega_p)\cap \zeta^*f^*(h))
= -2l -ek + \sum_{i=1}^m e_i$.  Similarly, we compute that
$\text{deg}(\zeta^*f^*(h) \cap \zeta^*f^*(h)) = 2el+e^2k -\sum_{i=1}^m
ee_i + \sum_{i=1}^m e_i(e-e_i)$,
i.e. $-e\text{deg}(\zeta^*C_1(\omega_p) \cap \zeta^*f^*(h)) +
\sum_{i=1}^m e_i(e-e_i)$.  Finally, observe that we have the formula
\begin{equation}
\text{deg}\lt( \zeta^* \sum_{i=1}^{\lt[\frac{e}{2}\rt]} \sum_{j=0}^{r}
i(e-i) \mathcal{D}_{i,j} \rt) = \sum_{i=0}^m e_i(e-e_i).
\end{equation}
So we conclude that 
\begin{equation}
\text{deg}\zeta^* p_*\lt(C_1(\omega_p) \cap f^*(h) \rt) =
-\frac{1}{e}\text{deg}(\zeta^* \mathcal{H}) + \frac{1}{e}
\sum_{i=1}^{\lt[\frac{e}{2} \rt]} \sum_{j=0}^r i(e-i)
\text{deg}(\zeta^* \mathcal{D}_{i,j})
\end{equation}
just as required.
\end{proof}

\medskip\noindent
\begin{prop}\label{appChern1} 
On $\Kbm{0,r}{\PP^n,e}$, the $\QQ$-divisor class of the first Chern
class of $\mc{P}_d = p_*(f^*\OO_{\PP^n}(d))$ equals
\begin{equation}
C_1(\mc{P}_d) = \frac{d}{2e}\lt[ (ed+1)\mathcal{H} -
\sum_{i=1}^{\lt[\frac{e}{2}\rt]} \sum_{j=0}^{r}
i(e-i)\mathcal{D}_{i,j} \rt].
\end{equation}
More generally, for $\textbf{d} = (d_1,\dots,d_c)$, the $\QQ$-divisor
class of the first Chern class of $\mc{P}_{\textbf{d}} =
p_*(f^*\OO_{\PP^n}(\textbf{d}))$ equals
\begin{equation}
C_1(\mc{P}_{\textbf{d}}) = \frac{1}{2e}
\lt( \prod_{k=1}^c (ed_k + 1) \rt) \lt[
\lt(\sum_{k=1}^c d_k\rt) \mathcal{H} + \lt(\sum_{k=1}^c
\frac{d_k}{ed_k+1} \rt) \sum_{i=1}^{\lt[\frac{e}{2}\rt]} \sum_{j=0}^c
i(e-i) \mathcal{D}_{i,j} \rt].
\end{equation}
\end{prop}

\medskip\noindent
\begin{proof}
Substituting the result from Lemma~\ref{appChernpre} into the
Grothendieck-Riemann-Roch formula yields, up to terms in
$\text{CH}^2(\Kbm{0,r}{\PP^n,e}) \otimes \QQ$,
\begin{equation}
\text{ch}[\mc{P}_d] = (ed+1) + \frac{d^2}{2}\mathcal{H} +
\frac{d}{2e} \lt( \mathcal{H} - \sum_{i=1}^{\lt[ \frac{e}{2} \rt]}
\sum_{j=0}^r i(e-i) \mathcal{D}_{i,j} \rt) + \dots
\end{equation}
Since $\text{ch}[\mc{P}_d] = \text{rank}(\mc{P}_d) + C_1(\mc{P}_d) +
\dots$, the first part of the proposition follows.

\medskip\noindent
Since $\mc{P}_{\mathbf{d}} \cong \oplus_{k=1}^c \mc{P}_{d_i}$, we have
the formula
\begin{equation}
C_1(\mc{P}_{\mathbf{d}}) = \lt( \prod_{k=1}^c
\text{rank}(\mc{P}_{d_k}) \rt) \sum_{k=1}^c
\frac{C_1(\mc{P}_{d_k})}{\text{rank}(\mc{P}_{d_k})}. 
\end{equation}
Substituting the first part of the proposition gives the second part
of the proposition.
\end{proof}

\medskip\noindent
The following corollaries follow immediately from Proposition
~\ref{appChern1}.  We state them as separate corollaries for
notational convenience.

\medskip\noindent
\begin{cor}\label{appCherncor1} 
Let $s\in H^0(\PP^n,\OO_{\PP^n}(d))$ be a section with zero locus
$X\subset \PP^n$.  Consider the locally free sheaf $\mc{P}_d$ on
$\Kbm{0,0}{\PP^n,e}$.  Suppose that the section $\sigma = f^*s$ of
$\mc{P}_d$ is a regular section, i.e. $\Kbm{0,0}{X,e}$ has the
``expected'' codimension $ed+1$ in $\Kbm{0,0}{\PP^n,e}$.  Then the
$\QQ$-divisor class of the first Chern class of the dualizing sheaf on
$\Kbm{0,0}{X,e}$ equals
\begin{equation}
\frac{1}{2e}\lt[ \lt( (d^2 - n -1)e - (n+1-d) \rt) \mathcal{H} +
\sum_{i=1}^{\lt[\frac{e}{2} \rt]} \lt((n+1-d)i(e-i) - 4e\rt)
\mathcal{D}_{i,0} \rt].
\end{equation}
\end{cor} 

\medskip\noindent
\begin{cor}\label{appCherncor2} 
Let $\mathbf{d} = (d_1,\dots, d_c)$.  Let $s\in
H^0(\PP^n,\OO_{\PP^n}(\mathbf{d}))$ be a section with zero locus
$X\subset \PP^n$.  Consider the locally free sheaf
$\mc{P}_{\mathbf{d}}$ on $\Kbm{0,r}{\PP^n,e}$.  Suppose that the
section $\sigma=f^*s$ of $\mc{P}_{\mathbf{d}}$ is a regular section,
i.e. $\Kbm{0,r}{X,e}$ has the ``expected'' codimension $e|\mathbf{d}|
+ c$ in $\Kbm{0,r}{\PP^n,e}$.  Then the $\QQ$-divisor class of the
first Chern class of the dualizing sheaf on $\Kbm{0,r}{X,e}$ equals
\begin{eqnarray*}
\frac{1}{2e^2}\lt[e \lt( \prod_{k=1}^c ed_k+1 \rt) \lt( \sum_{k=1}^c
d_k \rt) + 2r - (n+1)(e+1) \rt] \mathcal{H} -\frac{2}{e} \sum_{j=1}^r
\mathcal{L}_j + \\ 
 \sum_{j=0}^r
 j\mathcal{D}_{0,j} + 
\frac{1}{2e^2}\sum_{i=1}^{\lt[ \frac{e}{2} \rt]} \sum_{j=0}^r
 \lt[ \lt(\prod_{k=1}^c ed_k +1\rt)\lt(\sum_{k=1}^c \frac{ed_k}{ed_k+1} \rt)
i(e-i) + \rt. \\ 
 \lt. (n+1)ei(e-i) + 2e^2j - 4eij + 2ri^2 - 4e^2 \rt] \mathcal{D}_{i,j}
\end{eqnarray*} 
\end{cor}

\medskip\noindent
\begin{lem}\label{applembigB} 
Let $X\subset \PP^n$ be a projective scheme.
\begin{enumerate}
\item 
If every geometric generic point of $\Kbm{0,0}{X,e}$ parametrizes a
stable map which maps birationally to its image, then the pullback of
the $\QQ$-divisor class $\mathcal{H}$ in the $\QQ$-Picard group of
$\Kbm{0,0}{X,e}$ is big.  Moreover, when we pull this divisor class
back to the seminormalization of $\Kbm{0,0}{X,e}$ it is Cartier and
base-point-free.
\item
If every geometric generic point of $\Kbm{0,0}{X,e}$ parametrizes a
smooth rational curve in $X$ which is $a$-normal, then the pullback of
the $\QQ$-divisor class $C_1(\mc{P}_a)$ is an effective Cartier
divisor.
\item
If every geometric generic point of $\Kbm{0,0}{X,e}$ parametrizes a
stable map with irreducible domain, then for $i=1,\dots, \lt[
\frac{e}{2} \rt]$, the pullback of $\mathcal{D}_{i,0}$ is an effective
$\QQ$-Cartier divisor.
\end{enumerate}
\end{lem}

\medskip\noindent
\begin{proof}
To prove Item ($1$), we replace $\Kbm{0,0}{X,e}$ by its
seminormalization $\Kbm{0,0}{X,e}_{\text{sn}}$.  Consider the
universal curve $p:\mc{C} \rightarrow \Kbm{0,0}{X,e}_{\text{sn}}$ and
the universal morphism $f:\mc{C} \rightarrow X$.  Form the closed
image subscheme $C \subset \Kbm{0,0}{X,e}_{\text{sn}} \times X$ of
$(p,f)$.  Now $C$ is a \emph{well defined family of algebraic cycles}
in the sense of ~\cite[Defn. I.3.10]{K}.  By ~\cite[Thm. I.3.21]{K},
there is a Chow variety $\text{Chow}_{1,e}(X)$ and an induced morphism
$\Kbm{0,0}{X,e}_{\text{sn}} \rightarrow \text{Chow}_{1,e}(X)$.  By the
construction in ~\cite[Section I.3.23]{K}, there is an ample
invertible sheaf on $\text{Chow}_{1,e}(X)$ such that $\mc{H}$ is the
pullback of this ample invertible sheaf.  Moreover, by our assumption
that every geometric generic point of $\Kbm{0,0}{X,e}$ parametrizes a
stable map which maps birationally to its image, the morphism above is
generically finite.  Therefore $\mc{H}$ is base-point-free and big.
This proves Item ($1$).

\medskip\noindent
Let $W\subset H^0(\PP^n, \OO_{\PP^n}(a))$ be a general vector subspace
of dimension $ae+1$.  There is an induced map $W \otimes_\CC
\OO_{\Kbm{0,0}{X,e}} \rightarrow \mc{E}_a$ and the first Chern class
is simply the locus where this map fails to be an isomorphism.  By the
assumption that every geometric generic point of $\Kbm{0,0}{X,e}$
parametrizes a stable map which is $a$-normal, there is no irreducible
component of $\Kbm{0,0}{X,e}$ which is contained in this locus.
Therefore this locus is a Cartier divisor which is effective.

\medskip\noindent
Finally, by construction the boundary divisors are effective
$\QQ$-Cartier divisors on $\Kbm{0,0}{\PP^n,e}$.  If no irreducible
component of $\Kbm{0,0}{X,e}$ is contained in the boundary, then the
pullbacks of the boundary divisors are effective $\QQ$-Cartier
divisors on $\Kbm{0,0}{X,e}$.
\end{proof}

\medskip\noindent
\begin{cor} \label{cor-Kbig} 
Let $X\subset \PP^n$ be a general hypersurface of degree $d$.
\begin{enumerate}
\item
If $d < \text{min}(n-3,\frac{n+1}{2})$ and $d^2 \geq n+2$, then for
$e>>0$ the canonical divisor of $\Kbm{0,0}{X,e}$ is big.
\item
If $d < \text{min}(n-6,\frac{n+1}{2})$ and $d^2 + d \geq 2n + 2$, then
for every $e>0$ the canonical divisor of $\Kbm{0,0}{X,e}$ is big.
\end{enumerate}
In particular, if also $e \geq 3$ and $d+e \leq n$ or $e=2$ and $d+3
\leq n$, then Conjecture~\ref{conj-2} implies that $\kbm{0,0}{X,e}$ is
of general type.
\end{cor}

\medskip\noindent
\begin{proof}
When $d<\frac{n+1}{2}$, then it follows from ~\cite[Prop. 7.4]{HRS2}
that $\Kbm{0,0}{X,e}$ satisfies the hypotheses of Item ($1$) and Item
($3$) of Lemma~\ref{applembigB}.  Combining this with the formula from
Corollary~\ref{appCherncor1}, we get Item ($1$) and Item ($2$).

\medskip\noindent 
Finally, by Proposition~\ref{prop-c1c2} and Remark~\ref{rmk-c1c2},
when $e\geq 3$ and $d+e \leq n$ or $e=2$ and $d+3 \leq n$, then the
coarse moduli map $\Kbm{0,0}{X,e} \rightarrow \kbm{0,0}{X,e}$ is an
isomorphism away from codimension $2$ so that the canonical bundle of
$\kbm{0,0}{X,e}$ equals the canonical bundle of $\Kbm{0,0}{X,e}$.
Therefore Conjecture~\ref{conj-2} implies that $\kbm{0,0}{X,e}$ is of
general type.
\end{proof}

\bibliography{my}
\bibliographystyle{abbrv}

\end{document}